\documentclass[12pt,reqno]{amsart}   
\usepackage{amsmath,amsfonts,amssymb,graphicx,amscd,amsthm,amsbsy,epsf}
\usepackage{comment}
\numberwithin{equation}{section}
\newtheorem{thm}{Theorem}[section]
\newtheorem{lem}[thm]{Lemma}
\newtheorem{prop}[thm]{Proposition}
\newtheorem{cor}[thm]{Corollary}
\theoremstyle{definition}
\newtheorem{defn}{Definition}[section]
\newtheorem{remk}{Remark}[section]

\newtheorem{property}{property}[section]

\DeclareMathOperator*{\Imag}{Im}
\DeclareMathOperator*{\cov}{cov}
\DeclareMathOperator*{\var}{var}
\DeclareMathOperator*{\diam}{diam}
\DeclareMathOperator*{\spec}{Spec}
\def\uno{\text{\large{\bf{1}}}}
\def\Pr{{\mathbb P}}
\def\Ep{{\mathbb E}}

\def\real{\Bbb R}

\def\beqn{\begin{eqnarray*}}
\def\eeqn{\end{eqnarray*}}
\newcommand{\beqnl}[1]{\begin{eqnarray} \label{#1}}
\newcommand{\eeqnl}{\end{eqnarray}}

\begin{document}

\title[CLT
 for dynamical systems with weak random 
 noise]{ Renormalization and 
Central limit theorem for critical dynamical systems
 with weak external noise }

\author[O. D\'{i}az--Espinosa and 
R. de la Llave]{Oliver D\'{i}az--Espinosa \and Rafael de la Llave}
\address{Department of mathematics, The University of Texas at Austin,
Austin, TX 78712}
\email{odiaz@math.utexas.edu}
\email{llave@math.utexas.edu}

\begin{abstract} 
We  study of the effect of weak  noise on critical one 
dimensional  maps; that is, maps with a renormalization
theory. 

 We establish a one dimensional central limit theorem for 
weak noises and  obtain Berry--Esseen
estimates for the rate of this convergence.

We analyze in detail maps at the accumulation of period doubling
and critical circle maps with golden mean rotation number.
Using renormalization group methods,
we derive  scaling relations for several features of the effective noise
after long times. We use these scaling relations to
show that the central limit theorem
for weak noise holds in both examples.

We note that, for the results presented here, it is 
essential that the maps have parabolic behavior. They are
false for hyperbolic orbits. 

\end{abstract}


\maketitle

\rightline{\emph{To the memory of D. Khmelev}}

\section{Introduction}\label{intro}
The goal of  this paper is to develop a rigorous renormalization theory 
for weak 
noise  superimposed to one dimensional systems whose orbits have 
some self-similar  structure.

Some examples we
consider in detail are period doubling and
critical circle maps with golden mean rotation
number.

To be more precise,  we
consider  systems of the form 
\begin{eqnarray}
x_{n+1}=f(x_n)+\sigma\xi_{n+1}
\label{ransys}
\end{eqnarray}
 where $f$ is a map of a 
one dimensional space ($\mathbb R$, \, $\mathbb T^1$ or 
$I=[-1,1]$) into itself,   $(\xi_n)$ is a sequence
 of real valued  independent mean zero  random variables of comparable
sizes, and   $\sigma >0$ is a small parameter --called noise level--
 that controls the size of the  noise.

We  will study the scaling limit
 of the effective noise of \eqref{ransys} for small noise level, 
and  large number of iterations of the system. 
For $f$ either 
a map at the accumulation of period doubling
or a critical circle  with golden mean rotation 
and the noises satisfying some mild conditions
(existence of moments and the like),
 we will show that the scaling limits of the noise 
resemble a Gaussian in an appropriate sense. 
That is,
 if we fix a small value of the  noise level $\sigma$ and 
then,  look at the effective noise after a long time, 
the distribution of effective noise, normalized to have variance one,
 will have  a distribution  close to a Gaussian.
See Theorems \ref{main} and \ref{berry-esseen}
 for 
precise statements.  We will also show that statistical
 properties of the noise,
namely Wick-ordered moments or \emph{cumulants} satisfy 
some scaling properties. Observe that Theorems 
\ref{main} and \ref{berry-esseen} are 
similar to the classical central limit 
and the Berry-Esseen theorems, even if the powers
which appear as normalization factors are different from 
those in the classical theorems.

 The papers of 
\cite{crut:sca,shra:sca2} considered  heuristically a
renormalization theory for weak Gaussian noise perturbing
one dimensional maps at the accumulation of period doubling.
The main result in those papers  was that after appropriately 
rescaling space and time, the effective noise of this renormalized system
satisfies some scaling relations.
The paper  \cite{vul:feiguniv}
developed a rigorous thermodynamic formalism for critical maps 
with period doubling. Among many other results, 
\cite{vul:feiguniv},
study the effect of noise on the ergodic theory of
these maps and 
 showed that for systems at the accumulation of
period doubling with weak noise, there is a
 stationary measure depending on the magnitude of noise 
that converges to  the invariant measure
in the attractor.
A very different rigorous renormalization theory for systems with noise is
obtained in \cite{colles:renno}.

The results we present here can be considered as a rigorous version of 
the theory of \cite{crut:sca,shra:sca2}. The theory developed here
also applies to noise of
arbitrary shape  and shows that the scaling limit is Gaussian. 
The main idea is that one can also renormalize 
other statistical properties of the map (cumulants) and, by 
analyzing the different rates of convergence of all these different 
renormalization operators (these are what 
we call the convexity properties of the spectral radii
(See Theorem \ref{Comparison}), we obtain that the effective noise is 
approximately Gaussian. The argument we present 
uses relatively little properties of the renormalization operator. 
Basically we just need that there is some convergence to a scaling limit. 
We also apply similar ideas to the case of circle maps. 
The results there are very similar. 

Besides the results for the renormalization, we show that the previous 
results about the renormalization group give information about 
 the behavior of the noise along a
whole orbit. 

There is some overlap between the results in this paper and 
some of the results in \cite{vul:feiguniv}. The main emphasis on 
the paper \cite{vul:feiguniv} is on statistical properties for a
fixed level of noise. In this paper,
we emphasize the behavior on single orbits for longer times 
but with weaker noises.  The paper \cite{vul:feiguniv}
uses mainly the thermodynamic formalism 
and this paper relies mainly on transfer operators. Nevertheless, there 
are some relations between the two points of view. Of 
course, the relation between thermodynamics and transfer 
operators goes back to the beginning of thermodynamic 
formalism. See \cite{mayer}. In 
\cite[p. 31]{vul:feiguniv} the authors introduce 
one of the cumulant operators 
we use and find a relation between its spectral radius and 
thermodynamic properties. We think it would be possible 
and interesting 
to develop thermodynamic formalism analogues of the convexity properties
of transfer operators obtained in Theorem~\ref{Comparison}. 
It would also be interesting to develop 
thermodynamic formalism analogues  of the arguments developed 
in Sections \ref{sec:feig_clt}, \ref{sec:circ_clt}
which allow to study the behavior along a whole orbit from 
the study of renormalization operators.

This  paper will be organized as follows. In Section \ref{results}
we state a general central limit theorem (Theorems \ref{main})
 and a result on 
Berry--Esseen estimates  (Theorem \ref{berry-esseen})for the 
 convergence in the central limit theorem for one dimensional
dynamical systems with weak noise. The main hypothesis of these results
is that some combinations of derivatives grow at certain rate
(see \eqref{lyapcon}). This condition is reminiscent of the 
classical Lindeberg--Lyapunov central limit theorem 
sums of independent random variables. 

One important class of systems that satisfy 
the  condition \eqref{lyapcon} of Theorem \ref{main}
 is that of  fixed points of 
renormalization operators. Specific results for systems at the 
accumulation of period doubling (Theorem \ref{fclt})
and for critical circle maps 
with golden mean rotation number(Theorem \ref{fclt}) will be stated in 
sections \ref{feigresult} and \ref{circresult}
respectively.

The rest of the paper is devoted to providing proofs of
the results above. In Section \ref{sec:cltgen},
Theorems \ref{main} and \ref{berry-esseen} are proved in detail.
The method of the proof is 
to use the Lindeberg--Lyapunov central limit theorem for 
a linearized approximation of the effective noise, and then to
control the error terms.
Some examples will also be
discussed in Section \ref{examples}. 

In Section \ref{perioddoubling} we study in detail
 unimodal  maps at the accumulation of period--doubling 
 \cite{feig:univ,coeck:itemap,coecklan:univ}, and prove
Theorem \ref{fclt}.
In section \ref{critcircmap} we study critical maps of the circle 
with golden mean rotation number and prove Theorem \ref{circclt}.

The techniques used in Sections \ref{perioddoubling} and \ref{critcircmap}
consists on introducing some auxiliary linear operators -- which we called
Lindeberg--Lyapunov
operators (see Sections \ref{linlyapoperatorfeig}  and
\ref{deffcomopcirc}) -- that describe the statistical 
properties of the renormalized noise. The crucial part of 
the argument is to show that the spectral radii of these operators
satisfy some convexity properties (Theorem~\ref{Comparison}), 
so that in the scaling limit, the properties of the noise 
are that of a Gaussian. This implies 
that the sufficient conditions of  Theorem \ref{main} hold for 
a sequence of times.  A separate argument developed in 
Section~\ref{sec:feig_clt}  shows that, from the 
knowledge at these exponentially separated times, one can obtain 
information of the orbit starting at zero for all times. From that, 
one can obtain information for orbits starting at any iterate of 
the critical point. 

\section{Statement of  results}\label{results}
Throughout this paper, 
we will make the following assumptions
\begin{itemize}
\item[A1.] $f:M\mapsto M$ is a $C^2$ map where
  $M=\mathbb R$, $I\equiv[-1,1]$ or 
$\mathbb T^1$.
\item[A2.] If  $M=I$, we will further assume that
 there is a number  $a>0$ such
that $f\in C^2([-1-a,1+a])$ and $f([-1-a,1+a])\subset I$.
\end{itemize}
 For any function
$f\in C(M)$, we will denote by $\|f\|_{C_0}=\sup_{x\in M}|f(x)|$.

Let  $(\xi_n)$ be  a sequence of independent random
 variables (defined in some probability space
 $(\Omega,\mathcal F, \Pr)$),
with  $p>2$ finite moments.
We will assume that 
\begin{itemize}
\item[A3($p$).] $\Ep[\xi_n]=0$ for all $n$, and that
$$c\leq \|\xi_n\|_2\leq \|\xi_n\|_p\leq C$$
for some constants $c,\,C$.
 Here $\|\xi_n\|_s=(\Ep[|\xi_n|^s])^{1/s}$, for any $s>0$.
\end{itemize}

A direct consequence  of A3($p$) and independence is that
\beqn
c\leq \left\|\max_{1\leq j\leq n}|\xi_j|\right\|_p\leq n^{1/p}\,C
\eeqn
This inequality will be useful 
to give sufficient conditions on the noise level to obtain 
a central limit theorem (see  Theorem \ref{main} below).
\subsection{General results for one--dimensional dynamical
systems with external weak noise}\label{sysextwnoise}
Let $x\in M$ and $\sigma>0$  be fixed. We consider the system
\beqnl{system}
x_m= f(x_{m-1})+\sigma \xi_m,\quad x_0=x
\eeqnl
 $\sigma$ will be referred to as the {\it ``noise level''}.

For $M=\mathbb R$ or
 $M=\mathbb T^1$,  define the  process $x_n(x,\sigma)$
 as the value at time $n$ of \eqref{system}.

For $M=I$, we  define  the process $x_n(x,\sigma)$  by 
\eqref{system} provided
that
 $$\{f(x_{j-1})+\sigma\xi_j\}_{j=1}^n\subset [-1-a,1+a]$$
Observe that
$$\bigcup_{j=1}^n\left\{ |f(x_{j-1})+\sigma\xi_j|>1+a\right\}\subset
\left\{\sigma \max_{1\leq j\leq n}|\xi_j|>a\right\}:=C_n,$$
hence, by  Chebyshev's inequality we get
$$\Pr[C_n]\leq\sigma\frac{\Ep[\max_{1\leq j\leq n}|\xi_j|]}{a}$$
Hence, for a small value of $\sigma$, the event 
$C_n$ occurs with low probability.
Therefore, in order to define $x_n(x,\sigma)$,
 it suffices to condition on the event 
$$\Omega\setminus C_n=
\left\{\sigma \max_{1\leq j\leq n}|\xi_j|\leq a\right\}$$

In the scaling limit, we will consider a  sequence of noise levels
$\{\sigma_n\}$ converging to $0$ so that $\Omega\setminus C_n$ is close
 to the whole  space $\Omega$. Notice that  if $(\xi_j)_{j\in\mathbb N}$
 is supported  on a compact  interval, then by
taking $\sigma_n$ small enough,  the events  $C_n$ will
be  empty sets.

\subsubsection{General central limit theorem for
one dimensional dynamical systems with random weak noise}
Our first result is a general central limit theorem for one--dimensional
dynamical systems with weak random noise (Theorem \ref{main}).
This result is based upon the classical Lindeberg--Lyapunov
central limit theorem for sums of random variables 
\cite[p. 44]{bill:convpr}.
This is reflected in the sufficient condition \eqref{lyapcon} in the
statement of Theorem \ref{main}.

We introduce the following notation
\begin{defn}\label{lyapunov_func}
Let $f$ be a map on $M$  satisfying A1, A2, and let
$(\xi_n)$ be a sequence of random variables with $p>2$ finite moments,
that satisfies A3($p$).
\begin{itemize}
\item[1.]
The Lyapunov functions  $\Lambda_s(x,n)$, for  $ s\geq 0$,
 and $\widehat{\Lambda}$ are defined by 
\begin{eqnarray}
\Lambda_s(x,n)&=&
\sum^n_{j=1} \left|
\left(f^{n-j}\right)^\prime\circ f^j(x)
\right|^s\label{lambdap}\\
\widehat{\Lambda}(x,n)&=& \max_{0\leq i\leq n} \sum^i_{j=0}\left|
\left(f^{i-j}\right)^\prime\circ f^j(x)\right| \label{lambdahat}
\end{eqnarray}
When needed, we will use the notation $\widehat{\Lambda}^f$ and 
$\Lambda^f_s$ to emphasize the dependence on $f$.
\item[2.] Let  $x\in M$ and   $\sigma>0$ be fixed. The 
linearized effective noise is defined as 
\beqnl{c1}
L_n(x)=\sum_{j=1}^n \left(f^{n-j}\right)^\prime\circ f^j(x)\xi_j
\eeqnl
\end{itemize}
\end{defn}
\begin{remk} 
It is very important to observe that
for each $s\geq0$, the  Lyapunov function $\Lambda_s(x,n)$ satisfies
\beqnl{lambdarelation}
\Lambda_s(x,n+m)=\left|(f^m)^\prime\circ f^n(x)\right|^s\Lambda_s(x,n) 
+ \Lambda_s(f^n(x),m)
\eeqnl
We will use \eqref{lambdarelation}  in the study of central limit 
theorems for systems near the accumulation of period doubling,
Section \ref{perioddoubling}, and for  critical maps of the circle,
Section \ref{critcircmap}.
\end{remk}
\begin{remk}\label{sum_moments}
The Lyapunov functions \eqref{lambdap} are used to estimate
the sums of the moments or order $s$ of the terms in \eqref{c1}.
In particular, notice that

\begin{equation}
\var[L_n(x)]=\sum_{j=1}^n\left(\left(f^{n-j}\right)^\prime
\circ f^j(x)\right)^2
\Ep[\xi^2_j]
\end{equation}

The assumption  A3($p$) and H\"{o}lder's inequality imply
that  there are constants $c$, $C$ such that
for any $0<s\leq p$,
\beqnl{order}
c\Lambda_s(x,n)\leq \sum_{j=1}^{n}
 \left|\left(f^{n-j}\right)^\prime\circ f^j(x)\right|^s\Ep[|\xi_j|^s]
\leq
C\Lambda_s(x,n)
\eeqnl
In particular, if  A3($p = 2$)  implies 
\beqnl{ordervar}
c\Lambda_2(x,n)\leq \var[L_n(x)]
\leq
C\Lambda_2(x,n)
\eeqnl
\end{remk}

The main result that we obtain for orbits of a dynamical 
system is that if a deterministic condition on the orbit 
(expressed in terms of Lyapunov functions) holds, then, 
the noise perturbing this orbit satisfies 
a central limit theorem. 

\begin{thm}\label{main}
Let $f$, $M$ be a function satisfying  A1 and A2, and let
$(\xi_n)$ be a sequence of independent 
random variables, with $p>2$ finite
moments, that satisfies A3($p$).
Suppose that  for some $x\in M$ there is an increasing  sequence
of positive integers $n_k$ such that
\beqnl{lyapcon}
\lim\limits_{k\rightarrow\infty}\frac{\Lambda_p(x,n_k)}
{\left(\Lambda_2(x,n_k)\right)^{p/2}}
= 0
\eeqnl

Let   $(\sigma_k)_k$ be a sequence of positive numbers. 
Assume furthermore either of the two conditions
\begin{itemize} 
\item[H1]
The noise satisfies A3$(p)$ with $ p > 2$.

And the sequence $(\sigma_k)$
satisfies:

\beqnl{weaknoise}
\lim\limits_{k\rightarrow\infty}
\frac{\|f''\|_{C_0}\|\max_{1\leq j\leq n_k}|\xi_j|\|_p^2
(\widehat{\Lambda}(x,n_k))^6 
\sigma_k}{\sqrt{\Lambda_2(x,n_k)}}
= 0
\eeqnl
\item[H2]
The noise satisfies A3$(p)$ with $ p \ge 4$.

And the sequence $(\sigma_k)$
satisfies:
\beqnl{weaknoise2}
\lim_{k\rightarrow\infty}
\frac{\|f''\|_{C_0}\|\max_{1\leq j\leq n_k}|\xi_j|\|_p^2
(\widehat{\Lambda}(x,n_k))^3 
\sigma_k}{\sqrt{\Lambda_2(x,n_k)}}
=0
\eeqnl
\end{itemize}

Then, there exists a sequence 
of events $B_k\in\mathcal F$ such that
\begin{itemize}
\item[M1.] $\lim_{k\rightarrow\infty}\Pr[B_k]=1$
\item[M2.] The processes defined by
\begin{eqnarray}
w_{n_k}(x,\sigma_k)&=&\frac{(x_{n_k}(x,\sigma_k)-f^{n_k}(x))}
{\sigma_k\sqrt{\var[L_{n_k}(x)]}}
\label{normprocess}\\
\widetilde{w}_{n_k}(x,\sigma_k)
&=&\frac{(x_{n_k}(x,\sigma_k)-f^{n_k}(x))\uno_{B_k}}
{\sqrt{\var[(x_{n_k}(x,\sigma_k)-f^{n_k}(x))\uno_{B_k}]}}
\label{normprocess2}
\end{eqnarray}
converge in distribution  to a standard Gaussian as  $k\rightarrow\infty$.
\end{itemize}

Furthermore, if  the sequence $\xi_n$ is supported on a 
compact set then,  we can choose $B_k=\Omega$ 
for all $k$.
\end{thm}

The sets $\Omega\setminus B_k$, which we call \emph{outliers}, 
 are events where large fluctuations of noise occur.
We will refer to
\beqnl{efn}
(x_{n_k}(x,\sigma_k)-f^{n_k}(x))\uno_{B_k}
\eeqnl
as the \emph{effective noise}. 

Condition \eqref{lyapcon} in  Theorem \ref{main} is closely 
related to the Lyapunov condition of the classical limit theorem 
for sums of independent random variables  applied to the
linearized effective noise $L_{n_k}(x)$ defined by \eqref{c1}.

We will show in Section \ref{proof:clt_gen} that
if the sequence of noise levels $(\sigma_k)$ satisfy
\eqref{weaknoise} then, the events where large fluctuations of noise occur
(outliers) have  low probability. As a consequence, we will have that the
linearized effective noise \eqref{c1} is, with large probability, 
 a very good approximation to
the effective noise defined by \eqref{efn}.

\begin{remk}
Note that there are two variants of the results \eqref{normprocess}
and \eqref{normprocess2} in Theorem~\ref{main} as well as 
two variants on the hypothesis.

The difference in the conclusions is that the effective noise is normalized 
in two  -- in principle different -- ways. The version \eqref{normprocess}
normalizes the effective noise by its variance and the version 
normalizes by the variance of the linear 
approximation (see Remark \ref{sum_moments}. 

It could, in principle happen that the difference between the 
linear approximation and true process converged to zero 
in probability but that had a significance contribution 
to the variance. This pathology can be excluded by assuming that 
the noise has sufficient moments and that the noise is weak enough. 
These two hypothesis can be traded off. In the first version of 
the hypothesis [H1], we use only $p > 2$ moments and a somewhat 
stronger smallness conditions in the size of noise  \eqref{weaknoise}. 
In the second version of the hypothesis [H2],
 we assume $p \ge 4$ moments, but the smallness
conditions in the noise are weaker. 

The subsequent Theorems~\ref{berry-esseen} will  be true under either 
of the hypothesis.
\end{remk}

\begin{remk} 
As we will see later,  in the proof, the conditions on $\sigma_k$
are just  upper bounds. If we consider two sequences 
$\tilde \sigma_k \le \sigma_k$ and $\sigma_k$ satisfies 
either of \eqref{weaknoise} or \eqref{weaknoise2}, then
$\tilde \sigma_k$ satisfies the same conditions. 

Furthermore, 
the bounds that we obtain for 
the proximity  $w_{n_k} (x,\tilde \sigma_k)$
to the standard Gaussian are smaller than 
the bounds that we obtain on the proximity of  $w_{n_k} (x,\sigma_k)$
to the standard Gaussian. 
\end{remk}

\subsubsection{Berry--Esseen estimates}
The relevance of the outliers will become more evident
in our second result, Theorem \ref{berry-esseen} below, 
which provides the  rate of convergence to Gaussian
in  Theorem \ref{main}. 

We will  use $\Phi(z)$  to denote the
distribution function of the standard Gaussian measure on 
$(\mathbb R, \mathcal B(\mathbb R))$, that is
$$\Phi(z)=\tfrac{1}{\sqrt{2\pi}}\int^z_{-\infty} e^{-t^2/2}\,dt$$

\begin{thm}\label{berry-esseen}
Let $f$, $(\xi_n)_n$ be as in Theorem \ref{main} and let $s=\min(p,3)$.
Assume that 
condition \eqref{lyapcon} holds at   some $x\in M$
If $\sigma_k$ is a sequence
of positive numbers such that
\beqnl{weakernoise}
\frac{(\widehat{\Lambda}(x,n_k))^3}{\sqrt{\var[L_{n_k}(x)]}}
\|f''\|_{C_0}\left\|\max_{1\leq j \leq n_k}|\xi_j|\right\|^2_s
\sigma_k \leq 
\left(\frac{\Lambda_s(x,n_k)}
{\left(\Lambda_2(x,n_k)\right)^{s/2}}\right)^2
\eeqnl
then, we have that 
\beqnl{BE}
\sup_{z\in\mathbb R}\left|\Pr[w_{n_k}(x)\uno_{B_k}\leq z]-
\Phi(z) \right|
\leq A \frac{\Lambda_s(x,n_k)}{(\Lambda_2(x,n_k))^{s/2}}
\eeqnl
where  $A>0$  depends only on $x$.
\end{thm}

\begin{remk}
When $M=\mathbb R$, or $\mathbb T^1$ and $f(x)=x+c$ for some constant 
$c$, Theorems \ref{main} and \ref{berry-esseen} coincide with
the    classical central limit  and Berry--Esseen theorems
 for  sums of independent random variables. Since $f''\equiv0$, the result
 holds regardless of what values $\sigma_k$ takes, and the outliers
are empty sets.
\end{remk}
\begin{remk}
In section  \ref{sec:cltgen} we will consider the map $f(x)=2x$.
This system does not satisfy  \eqref{lyapcon}, and 
indeed the conclusion of Theorem \ref{main} fails. 
Systems with enough hyperbolicity satisfy other types of central limit
theorems for weak noise \cite{ kifer:randym}, or  even in the absence
 noise  \cite{liv:clt}, \cite{fmnicoltorok}. 
Those results are very different from the ones 
we consider in this paper.
\end{remk}
\subsubsection{Sketch of the proof of Theorems \ref{main} and 
\ref{berry-esseen}.}
The proof of these results is obtained in Section~\ref{sec:cltgen} 
by showing that:
\begin{itemize}
\item[1)] The linear approximation to the process satisfies a central 
limit theorem (or a Berry-Esseen theorem)
\item[2)]  Under smallness conditions on the noise level 
 (see \eqref{weaknoise}), the linear 
approximation is much larger than the Taylor reminder, 
so that we can transfer the Gaussian behavior from one to the other. 
\end{itemize}

The main  source of difficulties in the proof
are situations when 
the noise is much larger than expected from the statistical properties
of the linear approximation (outliers). These events, of course have 
small probability and, therefore, do not affect the convergence
in probability. However,
it could  happen in principle that
they change the variance. We will see in Section \ref{sec:varcompar}
that the variance of the effective noise and that of its linear
approximation are asymptotically equal. As a consequence, we will
have the effective noise normalized by its variance converges in
distribution to the standard Gaussian.

The existence of moments of high order will provide an improvement
on the choice of the noise level. In particular, we will see that
if the noise has compact support, large fluctuations of the
effective noise never occur; that is, the outliers are empty sets.

The procedure of  cutting off outliers in Theorem \ref{main}
is similar to the process of elimination of \emph{``Large fields''}
that occurs in the rigorous study of renormalization group in
\cite{gakup,gakokup}.


In this paper, we consider two examples of maps that have a 
renormalization theory, and for which Theorems \ref{main} and
\ref{berry-esseen} apply. Namely,  systems at the accumulation 
of period doubling  and critical circle maps with 
golden mean rotation number.

\subsection{Results for systems at the accumulation of period 
doubling}\label{feigresult}
In section \ref{symmunimodalclt} we consider systems of the 
form \eqref{ransys}
where $f$ is a $2k$--order  analytic unimodal map of the interval
$I$ onto itself. That is  $f(0)=1$, $f^{(j)}(0)=0$,
for $1\leq 2k-1$, $f^{(2k)}(0)\neq0$, and  
$xf'(x)< 0$ for $x\neq0$. Here, $f^{(j)}$ denotes the $j$--th order 
derivative of $f$.  

The period doubling renormalization 
group operator $T$ acting on the space of unimodal maps is defined by
$$Tf (x)= f^2(\lambda_f x)/\lambda_f$$
where $\lambda_f=f(1)$.
We refer to Section \eqref{symmunimodalclt}
for a precise definition of unimodal maps and the period-doubling
renormalization operator.

 For each $k$, there is a set of analytic functions 
$\mathcal W_s(g_k)$, such for maps $f\in\mathcal W_s(g_k)$, we have that 
 $T^nf$ converges to a universal function $g_k$ which is
a fixed point of $T$ (see 
\cite{feig:univ,lan:casproof,Coullet-Tresser, eps,sull,mart,mel,JacobsonS02}).

We will show in Section \ref{perioddoubling} that
\begin{thm}\label{fclt}
Let $f\in\mathcal W_s(g_k)$ then, 
\begin{itemize}
\item[F1.] For any $x=f^l(0)$, $l\in\mathbb N$  and any $p>2$,
\beqn
\lim_{n\rightarrow\infty}\frac{\Lambda_p(x,n)}{\{\Lambda_2(x,n)\}^{p/2}}
=0
\eeqn
\end{itemize}
Let $(\xi_n)$ be
a sequence of independent random variables which have
$p>2$ finite moments, and that  satisfies A3($p$).
\begin{itemize}
\item[F2.] Let $w_n(x,\sigma)$ be 
$$w_n(x,\sigma)=
\frac{\left(x_n(x,\sigma_n)-f^n(x)\right)}
{\sigma \sqrt{\sum_{j=1}^n\left( \left(f^{n-j}\right)^\prime
\circ f^j(x)\right)^2
\Ep[\xi^2_j]}}$$
Then, there is a constant  $\gamma>0$,
such that if
$$\lim_n\sigma_nn^{\gamma+1}=0$$
then, for each  $x=f^l(0)$, $l\in\mathbb N$  there are events 
$\{B_n(x)\}_n\subset \mathcal F$ with 
$$\lim_n\Pr[B_n(x)]=1$$
 such that  $w_n(x,\sigma_n)\uno_{B_n(x)}$ and hence
 $w_n(x,\sigma_n)$
converge in distribution to the standard Gaussian.
\item[F3.] Furthermore, 
 there are constants $\alpha>0$
 and $\theta>0$ depending on $p$ such that  
if  $\sigma_n\leq n^{-\theta}$, then
\beqnl{befeig}
\sup_{z\in \mathbb R}|\Pr[w_n(x,\sigma_n)\uno_{B_n(x)}\leq z]-
\Phi(z)|\leq C_x n^{-\alpha}
\eeqnl
\end{itemize}
\end{thm}
Some explicit values for $\gamma$, $\alpha$, and $\theta$ that
follow from the arguments of the proof are given
in  Section \ref{sec:feig_clt}.

In \cite{DiazL}, empirical values for $\gamma$ are obtained 
for the quadratic Feigenbaum fixed point. The values obtained
there suggest that the sequence of level of noise $\sigma_n$
decays as a power $\gamma_*$ of the number of iterations, that 
is not very different from the one we obtain in
 Section \ref{sec:feig_clt}.

\subsection{Results for critical circle maps}\label{circresult}
In Section \ref{critcircmap} we study critical maps of the circle 
with golden mean rotation number \cite{Lanf,EpsFixed,deFariaI}. 
That is, we consider strictly 
increasing analytic maps $f$ such that $f(x+1)=f(x)+1$, 
 $f^{(j)}(0)=0$  for $j=1\ldots,2k$,  $f^{(2k+1)}(x)\neq0$, and
 $\lim_n (f^n(x)-x)/n=(\sqrt{5}-1)/2 \equiv \beta$. 

\subsubsection{Central limit theorem for Fibonacci times}
Recall that the sequence of Fibonacci numbers
$(Q_n)$, is defined by $Q_0=1=Q_1$, $Q_{n+1}=Q_n+Q_{n-1}$.
Any integer $n$ admits a unique Fibonacci decomposition
$$n=Q_{m_0}+\ldots + Q_{m_{r_n}}$$
where $m_0>\ldots>m_{r_n}>0$ are non-consecutive integers,
(i. e. $m_i \ge m_{i -1} + 2$). 
Notice that $r_n+1$ is the number of terms in the Fibonacci
representation of $n$, and that $r_n\leq m_0\leq [\log_{\beta^{-1}}n]+$
(Here $[\quad]$ stands for the integer part function).

\begin{thm}\label{circclt}
Let $f$ be a critical circle map.
If $\{n_k\}$ is a increasing sequence of integers such that
\beqnl{circlacunar}
\lim_{k\rightarrow\infty}\frac{r_{n_k}}{\log_{\beta^{-1}}n_k}=0
\eeqnl
then,
\begin{itemize}
\item[C1.] For all $x=f^l(0)$, $l\in\mathbb N$, and any  $p>2$, 
\beqn
\lim_{k\rightarrow\infty}\frac{\Lambda_p(x,n_k)}
{\{\Lambda_2(x,n_k)\}^{p/2}}
=0
\eeqn
\end{itemize}
Let  $(\xi_n)$ be a sequence of random independent 
variables that have
$p>2$ moments, and that satisfies A3($p$).
\begin{itemize}
\item[C2.] Let $w_n(x,\sigma)$ be the process defined by
$$w_n(x,n)=
\frac{\left(x_{n}(x,\sigma)-f^{n}(x)\right)}
{\sigma \sqrt{\var[L_{n}(x)]}}
$$
Then, there is a   constant $\gamma>0$ such that if
$$\lim_k\sigma_k n_k^{\gamma+1}=0$$ 
there are events $\{B_k(x)\}\subset\mathcal F$ with
$$\lim_k\Pr[B_k(x)]=1$$
such that $w_{n_k}(x,\sigma_k)\uno_{B_k(x)}$ and hence $w_{n_k}(x,\sigma_k)$
 converge in distribution to the standard Gaussian.
\item[C3.] Furthermore,  there  are constants $\tau>0$ and
 $\upsilon>0$ depending on $p$ such that
 if $\sigma_k\leq {n_k}^{-\tau}$, then
$$\sup_{z\in \mathbb R}|\Pr[w_{n_k}(x,\sigma_k)\uno_{B_k(x)}\leq z]
-\Phi(z)|
\leq D_x\,n_k^{-\upsilon}$$
\end{itemize}
\end{thm}
\subsubsection{Central limit theorem along the whole sequence of times}
In Section \ref{techcircsec} it is shown that the Lyapunov condition 
\eqref{lyapcon}, for orbits starting in $\{f^l(0):l\in\mathbb N\}$,
holds along the whole sequence of times provided that 
some numerical condition \eqref{hypothesiscirc}
is satisfied (see Proposition 
\ref{kurto0ncirc} and Theorem \ref{lyapconorbit0circ}).
The technical condition \eqref{hypothesiscirc} is some relation between
properties of fixed points of a renormalization operators and
the spectral radii of some auxiliary operators. It could be 
verified by some finite computation or implied by some monotonicity
properties of the fixed point theorem, 

Under the technical condition
\eqref{hypothesiscirc},
the central limit theorem holds along the whole sequence
of times. The proof presented gives results that do not depend
on the condition \eqref{hypothesiscirc}.  Namely, we show 
that there as a central limit theorem along sequences of 
numbers which can be expressed as sum of  sufficiently ``few'' Fibonacci 
numbers in terms of their size.  The hypothesis \eqref{hypothesiscirc} 
implies that all the numbers satisfy this property.

\begin{thm}\label{circclt2}
Let $f$ be a critical circle map. Under the numerical condition 
\eqref{hypothesiscirc} (see Section \ref{techcircsec})
\begin{itemize}
\item[C4.] If  $x=f^l(0)$, $l\in\mathbb N$ and   $p>2$, then
\beqn
\lim_{n\rightarrow\infty}\frac{\Lambda_p(x,n)}{\{\Lambda_2(x,n)\}^{p/2}}
=0
\eeqn
\end{itemize}
Let  $(\xi_n)$ be a sequence of random variables which have
$p>2$ moments, and that satisfies A3. Then,
\begin{itemize}
\item[C5.] There is a   constant $\delta>0$ such that if
$$\lim_k\sigma_n n^{\gamma+1}=0$$ 
there are events $\{B_n(x)\}\subset\mathcal F$ with
$$\lim_n\Pr[B_n(x)]=1$$
such that $w_n(x,\sigma_k)\uno_{B_n(x)}$ and hence $w_n(x,\sigma_n)$
 converge in distribution to the standard Gaussian.
\item[C3.] Furthermore, there  are constants $\tau>0$ and
 $\upsilon>0$ depending on $p$
such that if $\sigma_n\leq n^{-\tau}$, then
$$\sup_{z\in \mathbb R}|\Pr[w_n(x,\sigma_n)\uno_{B_n(x)}\leq z]
-\Phi(z)|
\leq D_x\,n^{-\upsilon}$$
\end{itemize}
\end{thm}

\subsubsection{Discussion of the results for critical maps
(Theorems \ref{fclt}, \ref{circclt} and \ref{circclt2})}
The proof of Theorems \ref{fclt} and \ref{circclt} are obtained
using renormalization methods.

Roughly, renormalization gives us control of the effects of noise on
small scales around the critical point for a fixed
increasing sequence of times (powers of 2 in the case of 
period doubling and Fibonacci numbers in the case of circle 
maps with golden mean rotation).
This gives us, rather straightforwardly, a central limit theorem 
when the orbit of zero is observed along these sequences of times.

To obtain a central limit theorem along the sequence of all 
times, we use the fact that an arbitrary 
number can be written as sum of these good numbers. 
We argue by approximation. We observe that the sequence of times 
accessible to renormalization is also the  
sequence of times at which the orbit of 
zero comes close to the origin. Hence, we can write the orbit of 
zero as a sum of approximate Gaussians.

The argument we present has some delicate steps. We need to 
balance how close is the approximation to the Gaussian 
(how fast is the convergence to the CLT) with how fast is the 
recurrence at the indicated times. 

In our approach, to prove a CLT along all times, 
we have to compute and compare  the two effects.
 This comparison depends on
quantitative  properties of 
the fixed point of the renormalization group and some of
the auxiliary operators.

In the period doubling case, the
properties required by our approach  can be established 
and proved by conceptual methods 
(convexity and the like)  from the 
properties of 
the fixed point.

In the case of circle maps however, 
our methods require a property (see \eqref{hypothesiscirc})
which
seems to be true numerically, but which we do not know how 
to verify using only analytical methods.

The analysis presented above raises the possibility that, for 
some systems, the weak noise limit could have a 
CLT along some sequences but not along other ones. Of course, it 
is possible that there are other methods of proof that do not require
such comparisons. We think that it would be interesting either 
to develop a proof that does not require these conditions
or to present an example of a system whose weak noise limit converges 
to Gaussian along a sequence of times but not others.  

\section{Proof of Theorem \ref{main} and 
Theorem \ref{berry-esseen}}\label{sec:cltgen}
In this section, we prove the  two general theorems about
one dimensional dynamical systems with weak random noise,
namely Theorem \ref{main} (a central limit theorem)
 and Theorem \ref{berry-esseen} (a Berry--Esseen theorem).

First, we  consider a linear approximation of the system
\eqref{ransys} and show that it satisfies a central limit
theorem and a  Berry--Esseen theorem, see section \ref{lineartheory}.
Then,  in section \ref{nonlineartheory} we make a comparison between 
the linear approximation process and $x_n(x,\sigma)$,
see  Lemma \ref{nonlinlem}.
The proofs 
of Theorems \ref{main} and 
\ref{berry-esseen} are given in Sections
 \ref{proof:clt_gen} and  \ref{proof:BE_gen} 
respectively. We show that
  for $\sigma$ small enough, see \eqref{weaknoise}, 
the linear approximation process
is a good approximation to $x_n(x,\sigma)$, except perhaps
in sets of decreasing probability, which we call outliers.
Since these sets have probability going to zero, they 
do not affect the convergence in probability. Showing that these 
outliers do not affect the variance requires some extra arguments,
which we present later. 

\subsection{Linear approximation of the effective noise}\label{lineartheory}
For $x\in M$ and $\sigma>0$  fixed, using 
Taylor expansion, we decompose the  process $x_n(x,\sigma)$
as 
\beqnl{tay}
x_n(x,\sigma) = f^n(x) +\sigma L_n(x) +
 \sigma^2Q_n(x,\sigma),
\eeqnl
where the  linear term $L_n$ is the sum of independent random 
variables defined by \eqref{c1}. 
The linear approximation process $y_k(x,\sigma)$
\beqnl{vareqn}
y_n(x,\sigma)= f^n(x) + \sigma L_n(x)
\eeqnl
 satisfies the following central limit theorem.
\begin{lem}\label{linlem}
Let $f$ be a function satisfying  A1 and A2, and 
$(\xi_n)$ be  a sequence of independent random variables with $p>2$
moments,  satisfying A3.
If condition \eqref{lyapcon} holds for some point $x\in M$ then,
\beqnl{linclt}
l_{n_k}(x)\equiv
\frac{L_{n_k}(x)}{\sqrt{\var[L_n(x)]}}
\eeqnl
converges in distribution to the standard Gaussian
as $k\rightarrow\infty$. Moreover, there is a universal
constant $C$ such that
\beqnl{esseen}
\sup_{z\in R}\left|\Pr[l_{n_k}(x)\leq z]-\Phi(z) \right|
\leq C 
\frac{\Lambda_{\min(p,3)}(x,n_k)}
{\left(\Lambda_2(x,n_k)\right)^{\min(p,3)/2}},
\eeqnl
\end{lem}
\proof
From \eqref{order} and \eqref{ordervar}, there is $c>0$ such that
\beqn
c^{-1}\frac{\Lambda_p(x,n_k)}{\sqrt{\Lambda_2(x,n_k)}}\leq
\frac{\sum_{j=1}^{n_k}\left|\left(f^{n_k-j}\right)^\prime\circ
f^j(x)\right|^p}{\sqrt{\var[L_{n_k}(x)]}}\leq
c\frac{\Lambda_p(x,n_k)}{\sqrt{\Lambda_2(x,n_k)}}
\eeqn
Hence, the  convergence of $l_{n_k}$ to the standard
Gaussian follows from the
 classical Lindeberg--Lyapunov
central limit theorem \cite[p. 44]{bill:convpr}, since
 condition \eqref{lyapcon} is 
equivalent in this case to the Lyapunov condition for sums of
independent random variables.

 The second assertion, \eqref{esseen},
 follows from Berry--Esseen's Theorem  for  sums 
of independent random variables \cite[p. 115]{petrov:sumran}
\endproof
\begin{remk}
Notice that the linear dependence of  $y_n(x,\sigma)$ on the random
variables $\xi_n$ implies that the convergence 
in Lemma \ref{linlem} is 
independent of the noise level $\sigma$.
The size of the noise $\sigma$ will be 
important in the  control of the  non-linear term $\sigma^2 Q_n(x,\sigma)$.
\end{remk}
\subsection{Nonlinear theory} \label{nonlineartheory}
In this section, we prove a result (Lemma \ref{nonlinlem})
that will help us make a  comparison between  
 $x_n(x,\sigma)$ and $y_n(x,\sigma)$.  This result is  analog to a
 well known result on  variational equations for ODE's  \cite{hart:ode}.
The method we use in the proof of
Lemma \ref{nonlinlem} is similar to the proof of the Shadowing
 Lemma in \cite{shubm,katok3p,katok3pen}. 

In the rest of this section, we will use the norm
$\|{\bf x}\|:=\max_{1\leq j\leq m}|x_j|$ for vectors ${\bf x}$
in Euclidean space $\mathbb R^m$, and the corresponding
induced norms for linear  and bilinear operators.
\begin{lem}\label{nonlinlem}
Suppose $f\in C^2(M)$, $x_0\in M$, and  let  
$\bf{\Delta}\in{\mathbb R}^{N+1}$ with
$\Delta_0=0$. Consider the sequences ${\bf{\bar{x}}}$, ${\bf x}$
and ${\bf c}$ in ${\mathbb R}^{N+1}$ defined by
\begin{itemize}
\item[\emph{a)}]$\bar{x}_0=x_0$,  $\bar{x}_{i+1}=f(\bar{x}_i)$,
\item[\emph{b)}] $ x_{i+1}=f(x_i) + \Delta_{i+1}$, 
\item[\emph{c)}] $c_{i+1}=f'(\bar{x}_i)c_i + \Delta_{i+1}$ with $c_0=0$,
\end{itemize}
for $i=0,\ldots,N-1$.\\
 Assume that
\beqnl{cond4cont}
\|{\bf{\Delta}}\|\|f''\|_{C_0}\{\widehat{\Lambda}(x_0,N)\}^2 \leq \frac14
\eeqnl
Then, 
\beqnl{nonlminl}
\begin{split}
\|{\bf{\bar{x} + c}} -{\bf x}\|& \leq 
\|{\bf{c}}\|^2\widehat{\Lambda}(x_0,N)\|f''\|_{C_0}\\
& \leq
\|{\bf{\Delta}}\|^2(\widehat{\Lambda}(x_0,N))^3 \|f''\|_{C_0} 
\end{split}
\eeqnl
\end{lem}
\proof Fixing a parameterization on $M$, we can assume without loss
of generality  that $M=\mathbb R$. We will prove \eqref{nonlminl} by
showing that ${\bf x}$ is a fixed point of the  contractive function 
$\psi$ defined in \eqref{contraction}.

Notice that
$$c_i=\sum_{j=0}^i \left(f^{i-j}\right)^\prime(\bar{x}_j)\Delta_j$$
If  ${\bf A}=(a_{ij})$ is the lower triangular matrix defined by
\beqn
a_{ij}=\left\{ \begin{array}{cl}
\left(f^{i-j}\right)^\prime(\bar{x}_j)& \text{for}
\quad0\leq j\leq i\leq N\\
0 & \text{otherwise}
\end{array}
\right.
\eeqn
 we  have that
\beqnl{c_A}
 {\bf{c}}={\bf A \Delta} 
\eeqnl
\beqn
\|{\bf{A}}\|=\widehat{\Lambda}(x_0,N)
\eeqn
Let us define a function $\tau:{\mathbb R}^{N+1}\times{\mathbb R}^{N+1}
\longrightarrow{\mathbb R}^{N+1}$
by
$$\tau({\bf z},\alpha)
=(x_0,f(z_0)+\alpha_1,\ldots,f(z_{N-1})+\alpha_{N})^\top$$
Observe that 
\begin{itemize}
\item[\it a)]
$\tau({\bf\bar{x}},0)=\bf{\bar{x}}$
\item[\it b)]
$\tau({\bf z},{\bf\Delta})={\bf z}$ if and only if $ {\bf z}={\bf x}$
\item[\it c)]  For  all 
$({\bf z},\alpha)\in\mathbb R^{N+1}\times\mathbb R^{N+1}$
\beqn
D\tau({\bf z},\alpha)& =& D\tau({\bf z},0)\\
D_{11}\tau({\bf z},\alpha)& =& D_{11}\tau({\bf z},0)
\eeqn
\end{itemize}

The following identities will be useful
\begin{eqnarray}
 D_1\tau({\bf\bar{x}},0){\bf{c}}+D_2\tau({\bf\bar{x}},0)\bf{\Delta}
&=& {\bf c}
\label{iii}\\
(D_1(\tau({\bf{\bar x}},0))-I)^{-1} &=& {\bf A} \label{d1A}
\end{eqnarray}
A direct computation shows that for any point $({\bf z},\alpha)$ 
and vectors $[h,k]$, $[\tilde{h},\tilde{k}]\in{\mathbb R}^{N+1}
\times{\mathbb R}^{N+1}$
\begin{equation*}
\begin{array}{lcr}
D^2\tau({\bf z},\alpha)([h,k],[\tilde{h},\tilde{k}])
&=&(0,\tilde{h}_0f''(z_0)h_0,\ldots,\tilde{h}_{N-1}f''(z_{N-1})h_{N-1})^\top
\end{array}
\end{equation*}
Therefore
\beqnl{boundD2}
 \sup_{({\bf z},\alpha)}\|D^2\tau({\bf z},\alpha)\|=
\sup_{\bf z} \|D_{11}\tau({\bf z},0)\|\leq\|f''\|_{C_0}
\eeqnl
We define an auxiliary function
${\mathcal N}:
\mathbb R^{N+1}\times{\mathbb R}^{N+1}\longrightarrow \mathbb R^{N+1}$ 
by
\beqn
{\mathcal N}({\bf z},\alpha) = -{\bf A}(\tau({\bf z},\alpha)-{\bf z})+
{\bf z}
\eeqn
Using (a)--(c) we have 
\begin{itemize}
\item[\it d)] 
 ${\mathcal N}({\bf z},0)={\bf z}$ if and only if $ {\bf z}={\bf\bar{x}}$
\item[\it e)]
$\mathcal N({\bf z},{\bf \Delta})={\bf z}$ if and only if ${\bf z}={\bf x}$
\item[\it f)] $D_1{\mathcal N}({\bf z},\alpha) =D_1{\mathcal N}({\bf z},0)$
for all $({\bf z},\alpha)$.
\end{itemize}

It follows from \eqref{d1A}
\begin{eqnarray}
{\mathcal N}({\bf z},\alpha)-{\mathcal N}({\bf z},0)&=& 
-{\bf A}[\tau({\bf z},\alpha)-\tau({\bf z},0)]\label{vi}
\end{eqnarray}
\begin{equation}
\begin{array}{lcr}
D_1{\mathcal N}({\bf z},\alpha) &=&
-{\bf A}(D_1\tau({\bf z},0)-I)+I\\
\noalign{\vskip6pt}
&=&{\bf A}[D_1\tau({\bf\bar{x}},0)-D_1\tau({\bf z},0)]
\label{v}
\end{array}
\end{equation}
\begin{equation}
\begin{array}{lcr}
D_2{\mathcal N}({\bf z},\alpha)k &=&-{\bf A}[0,k_1,\ldots,k_N]
\label{v2}
\end{array}
\end{equation}
From \eqref{iii} we have
\beqnl{tautay}
\tau({\bf\bar{x}}+{\bf c},{\bf \Delta})=
{\bf\bar{x}}+{\bf c}+
\frac12 D_{11}\tau({\bf\tilde{\xi}},{\bf\tilde{\Delta}})({\bf c},{\bf c}),
\eeqnl
where  
 $({\bf\tilde{\xi}},{\bf\tilde{\Delta}})$ is a point on 
the  linear segment 
between $(\bar{{\bf x}}+{\bf c},{\bf \Delta})$ and $(\bar{{\bf x}},0)$.
Thus, by  \eqref{c_A}, \eqref{boundD2}, and \eqref{tautay} we have
\beqn
\|{\mathcal N}({\bf\bar{x}}+{\bf c},{\bf \Delta})-({\bf \bar{x}+c})\|
&\leq&
\frac12 \|{\bf A}\|\|{\bf c}\|^2  \|f''\|_{C_0}\\
&\leq&
\frac12 \|{\bf A}\|^3\|{\bf\Delta}\|^2 \|f''\|_{C_0}
\eeqn
  
Let
 $\bar{B}({\bf \bar{x}};r)$ be the closed ball in ${\mathbb R}^{N+1}$ 
centered at 
$\bar{\bf x}$ with radius $r=2\|{\bf \Delta}\|\|{\bf A}\|$. 
Then, by (\ref{v}), we have that  for all
 $(h,\alpha)\in{\mathbb R}^{N+1}\times{\mathbb R}^{N+1}$ and all
${\bf z}\in\bar{B}({\bf \bar{x}};r)$, there is
$\eta$ in the linear segment joining ${\bf z}$ and $\bar{x}$ such that
$$D_1{\mathcal N}({\bf z},\alpha)h=-{\bf A}D_{11}\tau(\eta,0)
(h,{\bf z}-{\bf\bar{x}})$$
Therefore, provided that condition \eqref{cond4cont}  holds, we have that
\beqnl{cont}
\|D_1{\mathcal N}({\bf z},\alpha)\|\leq 2 \|{\bf A}\|^2
\|{\bf\Delta}\|\|f''\|_{C_0}\leq\frac12
\eeqnl
for all ${\bf z}\in\bar{B}({\bf\bar{x}};r)$. 

On the other hand \eqref{vi}  and \eqref{cont} imply that
\beqn
\|{\mathcal N}({\bf z},{\bf\Delta})-{\bf\bar{x}}\|&\leq&
\|{\mathcal N}({\bf z},{\bf\Delta})-{\mathcal N}({\bf\bar{x}},{\bf \Delta})\|+
\|{\mathcal N}({\bf\bar{x}},{\bf\Delta})-{\mathcal N}({\bf\bar{x}},0)\|\\
&\leq&\frac12\|z-{\bf\bar{x}}\|+ \|{\bf A}\|\|{\bf\Delta}\|\leq r
\eeqn
It follows that  the function 
\beqnl{contraction}
\psi({\bf z})={\mathcal N}({\bf z},{\bf\Delta})
\eeqnl
is a contraction  by a factor $1/2$  of the closed ball $\bar{B}({\bf\bar{x}};r)$ into itself. 
Moreover, ${\bf x}$ is the unique fixed point.\\
Inequality \eqref{nonlminl} follows from
\begin{equation*}
\begin{array}{lcl}
\|{\bf\bar{x}}+{\bf c}-{\bf x}\|&\leq&
\|{\bf\bar{x}}+{\bf c}- \mathcal N ({\bf\bar{x}}+{\bf c},{\bf\Delta})\|
+ \|\mathcal N({\bf x},{\bf\Delta}) -
\mathcal N({\bf\bar{x}}+ {\bf c},{\bf\Delta})\|
\\
\noalign{\vskip6pt}
&\leq&
\frac12\|{\bf A}\|\|{\bf c}\|^2 \|f''\|_{C_0} + \frac12\|{\bf\bar{x}}+{\bf c}-{\bf x}\|
\end{array}
\end{equation*}
\endproof
\subsection{Proof of central limit theorem
 (Theorem \ref{main})} \label{proof:clt_gen}
The conclusion of Theorem \ref{main} will follow
by combining  Lemmas \ref{linlem} and  \ref{nonlinlem}. 
The key argument is to define events (outliers)  where
the linear approximation process is very different from 
$x_n(x,\sigma)$, and show that they have small probability.

Notice that if $\|f''\|_{C_0}=0$, then Theorem \ref{main} coincides with 
the Lindeberg--Lyapunov central limit theorem. Therefore, we will 
assume that $\|f''\|_{C_0}>0$.

\subsubsection{Outliers}\label{sec:outliers}
 For each $k$ and $j=1,\ldots, n_k$,  let $\Delta_j=\sigma_k\xi_j$, 
and let  ${\bf \Delta}$,
 $\bar{{\bf x}}$, ${\bf x}$ and ${\bf c}$ be as in
 Lemma \ref{nonlinlem}. Then, by 
 \eqref{lambdarelation}, \eqref{vareqn} and the definition
of the linear approximation process,  we have that
for $j=1,\ldots,n_k$
\beqn
f^j(x)&=&\bar{x}_j\\
y_j(x,\sigma_k)&=& \bar{x}_j + c_j\\
\sigma_k L_j(x)&=& c_j\\
x_j(x,\sigma_k)&=& x_j
\eeqn
 For any sequence
$(\sigma_k)$ of noise levels decreasing to $0$, 
Lemma \ref{nonlinlem} implies that in the event
\beqnl{goodappevent}
\bar{B}_k=\left[ \|f''\|_{C_0}\sigma_k(\widehat{\Lambda}(x,n_k))^2
\max_{1\leq j\leq n_k} |\xi_j|\leq\frac14\right],
\eeqnl
 the linear approximation process
$y_{n_k}(x,\sigma_k)$
\eqref{vareqn} is close to $x_{n_k}(x,\sigma_k)$.
Thus, we will restrict the process $x_{n_k}(x,\sigma_k)$ to events
$B_k\in\mathcal F$ such that
\beqnl{defnoutlier}
B_k=\left[\max_{1\leq j\leq n_k} |\xi_j|\leq a_k\right]\subset
\bar{B}_k
\eeqnl
for some appropriate sequence $a_k$. We will choose $(a_k)$ so that
\beqnl{defnoutlier2}
\lim_{k\rightarrow\infty}
\Pr[\Omega\setminus B_k]\rightarrow0
\eeqnl
For this purpose, it will be enough to define $(a_k)$ by
\beqnl{uppbounda_k}
a_k=\frac{1}{4\|f''\|_{C_0}} (\widehat{\Lambda}(x,n_k))^{-\beta}
\sigma_k^{-\alpha},
\eeqnl
where $\beta$ and $\alpha$ are chosen so that \eqref{defnoutlier}
and \eqref{defnoutlier2}   hold. This means that
\beqnl{condoutlier1}
\sigma_k(\widehat{\Lambda}(x,n_k))^2
\leq \sigma_k^{\alpha}(\widehat{\Lambda}(x,n_k))^{\beta}
\eeqnl
\beqnl{condoutlier2}
\lim_{k\rightarrow\infty}\left
\|\max_{1\leq j\leq n_k}|\xi_j|\right\|_p^{1/\alpha}
(\widehat{\Lambda}(x,n_k))^{\beta/\alpha}\sigma_k = 0
\eeqnl
\begin{remk}
Recall that if $M=I$, we 
define  the process $x_{n_k}(x,\sigma_k)$ conditioned on the event
\beqn
\Omega\setminus C_{n_k}= 
\left\{\sigma_k \max_{1\leq j\leq n_k}|\xi_j|\leq a\right\}
\eeqn
Chebyshev's inequality implies that
\beqn
\Pr[C_{n_k}]\leq
\left(\frac{\sigma_k\left\|\max_{1\leq j\leq n_k}|\xi_j|\right\|_p}{a}
\right)^p
\eeqn
Notice that \eqref{condoutlier2}, and  the fact that 
  $\widehat{\Lambda}(x,n_k)>1$ 
imply that $\Pr[C_{n_k}]\rightarrow0$ as $k\rightarrow\infty$.
\end{remk}
We refer to   $\{\Omega\setminus B_k\}_k$ 
\,($\{\Omega\setminus(B_k\cap C_k\}_k$ if $M=I$) as the 
sequence of outliers.

\subsubsection{Estimates on the effective noise}
In this section, we will show that the outliers have small
probability. This means that linearized effective noise
\eqref{c1} is a good approximation of the effective noise \eqref{efn}
with high probability.
As a consequence,  we will have that the
effective noise scaled by the standard deviation of
the linearized effective noise approaches a Gaussian.

Recall  from  \eqref{tay} that  the effective noise is
decomposed as
\beqn
(x_{n_k}(x,\sigma_k)-f^{n_k}(x))\uno_{B_k}=
\sigma_k L_{n_k}(x)\uno_{B_k}+
 \sigma^2_k\,Q_{n_k}(x,\sigma_k)\uno_{B_k}
\eeqn
To control the effect of nonlinear terms 
in the effective noise,
it will be enough  to require  that the variance of the noise of 
$\sigma_k^2 Q_{n_k}(x,\sigma_k)$
 is  small compared to the variance of the linearized effective
noise $\sigma_k L_{n_k}(x)$. In terms of the 
scaled processes $w_{n_k}(x,\sigma_k)$ and
$l_{n_k}(x)$ defined by \eqref{normprocess}
and \eqref{linclt}, this requirement is equivalent to
\beqnl{goodapprximationlimit}
\lim_{k\rightarrow\infty}
\|(w_{n_k}(x,\sigma_k)-l_{n_k}(x))\uno_{B_k}\|_2 = 0
\eeqnl

By \eqref{defnoutlier} and  Lemma \ref{nonlinlem} we have that
\beqn
|Q_{n_k}(x,\sigma_k)|\uno_{B_k}
&\leq& \|f''\|_{C_0}
(\widehat{\Lambda}(x,n_k))^3\max_{1\leq j\leq n_k}
|\xi_j|^2\uno_{B_k}
\eeqn
Then, by \eqref{uppbounda_k} we get
\beqnl{nonlinearpart}
\|(w_{n_k}(x,\sigma_k)-l_{n_k}(x))\uno_{B_k}\|_2 
\leq C\frac{(\widehat{\Lambda}(x,n_k))^{3-2\beta}
\sigma_k^{1-2\alpha}}{\sqrt{\Lambda_2(x,n_k)}}
\eeqnl
for some $C>0$.
To obtain \eqref{goodapprximationlimit}, it suffices to 
require that 
\beqnl{condoutlier3}
\lim_{k\rightarrow\infty}\frac{(\widehat{\Lambda}(x,n_k))^{3-2\beta}
\sigma_k^{1-2\alpha}}{\sqrt{\Lambda_2(x,n_k)}}=0
\eeqnl
We have the following Theorem
\begin{thm}\label{mainclt1}
Let $f$ and $(\xi_n)$ be as in Theorem \ref{main}.
Assume that the Lyapunov condition \eqref{lyapcon} holds
and that $(\sigma_k)$ satisfies \eqref{weaknoise}.
If $\alpha^{-1}=2=\beta$, then  $w_{n_k}(x,\sigma_k)$ 
converges in distribution to the standard Gaussian.
\end{thm}
\proof
Notice that the Lyapunov condition \eqref{lyapcon} implies that
$$\lim_{k\rightarrow\infty}\Lambda_2(x,n_k)\rightarrow\infty$$
 Hence,
by \eqref{condoutlier1} and \eqref{condoutlier3}, 
it will be enough to consider $0<\alpha\leq 1/2$ and 
$\beta\geq 2$.

In particular, if $\alpha^{-1}=2=\beta$ and $(\sigma_k)$ satisfies
\eqref{weaknoise}, then $\uno_{B_k}$ and 
 $(w_{n_k}(x,\sigma_k)-l_{n_k}(x))\uno_{B_k}$ 
converge in $L_{p/2}(\Pr)$ and hence in probability to $1$ and
$0$ respectively.
By Lemma \ref{linclt}
and the \emph{``converging together''}  Lemma
 \cite[p. 89]{durr:prob}, it follows that $w_{n_k}(x,\sigma_k)$
converges in distribution to the standard Gaussian.
\endproof

The result of Theorem \ref{mainclt1} is improved by 
considering  the weaker condition
\eqref{weaknoise2} for $(\sigma_k)$
 and  smaller outliers $\Omega\setminus\bar{B}_k$.

Notice that by Chebyshev's inequality
\begin{equation}
\begin{array}{rcl}
\Pr[\Omega\setminus \bar{B}_k]&\leq& 
\left(4\|f\|_{C_0}(\widehat{\Lambda}(x,n_k))^2
\left\|\max_{1\leq j\leq n_k}|\xi_n|\right\|_p\sigma_k\right)^p
\label{weakoutliercond}
\end{array}
\end{equation}

On the other hand, we have that
\beqnl{nonlinearpart2}
\left\|\frac{\sigma_kQ_{n_k}(x,\sigma_k)}
{\sqrt{\var[L_{n_k}(x)]}}\uno_{B_k}\right\|_{p/2} 
\leq C\frac{\|\max_{1\leq j\leq n_k}|\xi_j|\|_p^2
(\widehat{\Lambda}(x,n_k))^3
\sigma_k}{\sqrt{\Lambda_2(x,n_k)}}
\eeqnl
The following result is a direct 
consequence of the  ``\emph{convergence together}'' Lemma, 
 \eqref{weakoutliercond}, and 
 \eqref{nonlinearpart2}.
\begin{cor}
Let $f$ and $(\xi_n)$ be as in Theorem \ref{main}. Assume that
 the Lyapunov condition \eqref{lyapcon} holds and that
$(\sigma_k)$ is a sequence of positive numbers that satisfy
\eqref{weaknoise2}. Then, the process 
$w_{n_k}(x,\sigma_k)\uno_{\bar{B}_k}$ and hence $w_{n_k}(x,\sigma_k)$
converge in distribution to the standard Gaussian.
\end{cor}
\subsubsection{Comparison between the variance of the
effective noise and the variance of the linear approximation}
\label{sec:varcompar}
Observed that the  scaling  given
in \eqref{nonlinearpart} uses the variance of a random
variable $L_{n_k}(x)$ which is  defined in the whole space $\Omega$. 
Since  the outliers
occur with very low probability, 
it is  more natural to restrict all the quantities to 
 the complement of the
outliers, that is to the events $B_k$.
 
If $(\sigma_k)$ satisfies \eqref{weaknoise}, 
we will show that the variance of the effective noise
and that if its linear approximation are asymptotically
equal. By the converging together Lemma, we will have as
a consequence a central limit theorem for the 
effective noise scaled  by its variance.

\begin{lem}\label{variancecomparison1}
Let $f$ and $(\xi_n)$ be as in Theorem \ref{main}. 
Assume that the Lyapunov condition \eqref{lyapcon}
holds.
For any sequence $(d_k)$ of positive numbers such that
\beqnl{condond_k}
\lim_{k\rightarrow\infty}
 d^{-2(p-1)}_k
\frac{(\widehat{\Lambda}(x,n_k))^2}{\Lambda_2(x,n_k)}=0,
\eeqnl
define the event
$D_k= [\max_{1\leq j\leq n_k}|\xi_j|< d_k]$.
Then, 
\beqnl{linearoutliervslinear}
\lim_{k\rightarrow\infty}
\frac{\var[L_{n_k}(x)\uno_{D_k}]}
{\var[L_{n_k}(x)]}=1
\eeqnl
\end{lem}
\proof
Notice that
\beqn
\var[L_{n_k}(x)\uno_{D_k}]=\sum_{j,i=1}^{n_k}
\left(f^{n_k-j}\right)^\prime(f^j(x))
\left(f^{n_k-i}\right)^\prime(f^i(x))
\cov[\xi_j\uno_{D_k},\xi_i\uno_{D_k}]
\eeqn
where
$$\cov[\xi_j\uno_{D_k},\xi_i\uno_{D_k}]=
\Ep[\xi_j\xi_i\uno_{D_k}]-\Ep[\xi_j\uno_{D_k}]\Ep[\xi_i\uno_{D_k}]
$$
First, we estimate $\cov[\xi_j\uno_{D_k},\xi_i\uno_{D_k}]$ when
$j\neq i$. 

Since $\Ep[\xi_j]=0$ for all $j$, from the independence of $(\xi_n)$ we get
\beqn
|\Ep[\xi_j\uno_{D_k}]|&=&|\Ep[\xi_j\{|\xi|\leq d_k\}]|
\prod_{m\neq j}\Pr[|\xi_m|\leq d_k]\\
&=&|\Ep[\xi_j\{|\xi|>d_k\}]|
\prod_{m\neq j}\Pr[|\xi_m|\leq d_k]
\eeqn
By H\"{o}lder's and Chebyshev's inequality we get
\beqnl{estimmeanout}
|\Ep[\xi_j\uno_{D_k}]|\leq \|\xi_j\|_p^p d_k^{1-p}
\prod_{m\neq j}\Pr[|\xi_m|\leq d_k]
\eeqnl
Notice that \eqref{condond_k} implies that 
 $\lim_k d_k\infty$ and that $\lim_k\Pr[D_k]\rightarrow1$.
Consequently, the product of probabilities 
on the right hand side of
\eqref{estimmeanout} is close to (and smaller than) $1$.
Furthermore,  A($p$) implies that $(\xi_n)$
is bounded in $L_p(\Pr)$. Hence,
\beqnl{firstgoodboundmean}
|\Ep[\xi_j\uno_{D_k}]|\leq C d_k^{1-p}
\eeqnl

Similarly, the  independence of $(\xi_n)$ implies
\beqn
\left|\Ep[\xi_j\xi_i\uno_{D_k}]\right|&=&
\left|\Ep[\xi_j\{|\xi_j|\geq d_k\}]\Ep[\xi_i\{|\xi_i|
\geq d_k\}]\right|
\prod_{m\neq j,i}\Pr[|\xi_m|\leq d_k]\\
&=&
\left|\Ep[\xi_j\{|\xi_j|>d_k\}]\Ep[\xi_i\{|\xi_i|
> d_k\}]\right|
\prod_{m\neq j,i}\Pr[|\xi_m|\leq d_k]
\eeqn
Using H\"{o}lder's and Chebyshev's inequalities we get
\beqnl{estimcrosterm}
\left|\Ep[\xi_j\xi_i\uno_{D_k}]\right|\leq
\|\xi_j\|_p^p\|\xi_i\|_p^p d_k^{-2(p-1)}\prod_{m\neq j,i}\Pr[|\xi_m|\leq d_k]
\eeqnl
The product of probabilities on the right hand side of
\eqref{estimcrosterm} is close to (and smaller than) $1$. Therefore,
using A3($p$), we have:
\beqnl{secondgoodboundcov}
\left|\Ep[\xi_j\xi_i\uno_{D_k}]\right|\leq C d_k^{-2(p-1)}
\eeqnl

Combining \eqref{firstgoodboundmean} and \eqref{secondgoodboundcov} we
 have, for $j\neq i$, that
\beqnl{covineq}
|\cov[\xi_j\uno_{D_k},\xi_i\uno_{D_k}]|\leq C d_k^{-2(p-1)}
\eeqnl

When $j=i$, we obtain
\beqnl{varestim}
|\var[\xi_j\uno_{D_k}]-\var[\xi]|&=& \Ep[\xi^2_j\uno_{B^c_k}]
+ |\Ep[\xi_j\uno_{D_k}]|^2\nonumber\\
&\leq& C(\Pr[D_k^c]^{(p-2)/p}+ d_k^{-2(p-1)})
\eeqnl

Recall that $\Lambda_2(x,n_k)$ and $\var[L_{n_k}(x)]$ 
are of the same order (see \eqref{order}). Since $p>2$,
by  combining 
\eqref{covineq} and \eqref{varestim}, we get 
\beqn
\left|1-\frac{\var[L_{n_k}(x)\uno_{D_k}]}
{\var[L_{n_k}(x)]}\right|
&\leq& C\left( \frac{d_k^{-2(p-1)}(\Lambda_1(x,n_k))^2}
{\Lambda_2(x,n_k)}
+\Pr[D_k^c]^{(p-2)/p}\right)
\eeqn
Since $\Lambda_1(x,n_k)\leq \widehat{\Lambda}(x,n_k)$,
the conclusion of Lemma \ref{variancecomparison1}
follows from \eqref{condond_k} by passing to the limit.
\endproof
The following result shows that the variance of effective noise and that
of the its linear approximation are asymptotically equal
\begin{lem}\label{varasimptotics}
Let $f$ and $(\xi_n)$ be as in Theorem \ref{main}. Assume that
 the Lyapunov condition \eqref{lyapcon} holds and that
$(\sigma_k)$ is a sequence of positive numbers that satisfy
\eqref{weaknoise}. Then,
\beqnl{varcomparison}
\lim_{k\rightarrow\infty}
\frac{\var[(x_{n_k}(x,\sigma_k)-f^{n_k}(x))\uno_{B_k}]}
{\sigma^2_k\var[L_{n_k}(x)]} = 1,
\eeqnl
where $\{B_k\}$ are the events defined in Section \ref{sec:outliers}.
\end{lem}
\proof
If $(\sigma_k)$ satisfies \eqref{weaknoise}, we know from
Section \ref{sec:outliers} that 
\beqn
\lim_{k\rightarrow\infty}
\|(w_{n_k}(x,\sigma_k)-l_{n_k}(x))\uno_{B_k}\|_2=0
\eeqn
Therefore, by Lemma \ref{variancecomparison1} it is enough to
show that \eqref{condond_k} holds for $d_k=a_k$, where
$a_k$ is defined by \eqref{uppbounda_k}.

 Notice that
\beqn
a_k^{-2(p-1)}\frac{(\widehat{\Lambda}(x,n_k))^2}
{\Lambda_2(x,n_k)}&\leq&
C\left(\frac{\|\max_{1\leq j\leq n_k}|\xi_j|\|^2_p
(\widehat{\Lambda}(x,n_k))^6 \sigma_k}{\sqrt{\Lambda_2(x,n_k)}}
\right)^{4(p-1)}
\eeqn
Thus, the conclusion of Lemma \ref{varasimptotics} follows 
by passing to the limit.
\endproof
\begin{thm}\label{mainclt2}
Let $f$ and $(\xi_n)$ be as in Theorem \ref{main}. Assume that
 the Lyapunov condition \eqref{lyapcon} holds and that
$(\sigma_k)$ is a sequence of positive numbers that satisfy
\eqref{weaknoise}. Then, the process $\widetilde{w}_{n_k}(x,\sigma_k)$
defined by \eqref{normprocess2}
converges in distribution to the standard Gaussian.
\end{thm}
\proof
We know from Theorem \ref{mainclt1} that $w_{n_k}(x,\sigma_k)$
converges in distribution to the standard Gaussian. Notice
that
\beqn
\widetilde{w}_{n_k}(x,\sigma_k)=\sqrt{\frac{\var[L_{n_k}(x)]}
{\var[(x_{n_k}(x,\sigma_k)-f^{n_k}(x))\uno_{B_k}]}}
 w_{n_k}(x,\sigma_k)\uno_{B_k}
\eeqn
The conclusion of Theorem \ref{mainclt2} follows from 
Lemma \ref{varasimptotics} and the ``\emph{converging together}''
 Lemma.
\endproof
\subsubsection{Noise with finite moments of order $p\geq4$}
Now, we will prove Theorem~\ref{main} under hypothesis H2.

We will show that if the sequence $(\xi_n)$ has moments of order
$p\geq4$, then the choice of noise level $(\sigma_k)$ (and hence
the choice of outliers) is improved. 

\begin{thm}\label{mainclt3}
Let $f$ be as in Theorem \ref{main},
assume that the Lyapunov condition \eqref{lyapcon} holds.
Let  $(\xi_n)$ be a sequence of independent random variables
with $p\geq4$ finite moments, and that it satisfies
A3($p$).

For any sequence  $(\sigma_k)$ that satisfies \eqref{weaknoise2}, 
let $\{B_k\}$ be the events defined by
\eqref{goodappevent}. Then, the process
\beqnl{normprocess3}
\widehat{w}_{n_k}=
\frac{(x_n(x,\sigma_k)-f^{n_k}(x))\uno_{\bar{B}_k}]}
{\sigma^2_k\var[L_{n_k}(x)]}
\eeqnl
converges to the standard Gaussian.
\end{thm}
\proof
It suffices to verify that 
\beqn
d_k=\frac{1}{4\|f''\|_{C_0}}(\widehat{\Lambda}(x,n_k))^{-2}\sigma_k^{-1}
\eeqn
satisfy condition \eqref{condond_k} in Lemma \ref{variancecomparison1}.

Since $\sqrt{\Lambda_2(x,n_k)}\leq \Lambda_1(x,n_k)\leq 
\widehat{\Lambda}(x,n_k)$, we have that
\beqn
d_k^{-2(p-1)}\frac{(\widehat{\Lambda}(x,n_k))^2}{
\Lambda_2(x,n_k)} \leq
 C\left(\frac{\|\max_{1\leq j\leq n_k}|\xi_j|\|_p^2
(\widehat{\Lambda}(x,n_k))^3
\sigma_k}{\sqrt{\Lambda_2(x,n_k)}}\right)^{2(p-1)}
\eeqn
The conclusion of Theorem \ref{mainclt3}
 follows from \eqref{weaknoise2}
\endproof
\begin{remk}
 Notice  from the definition of the outliers, 
that if  the sequence of noises $(\xi_n)$ is supported on an 
interval $[-b,b]$, then for $k$ large enough $B_k=\Omega$, that is, the
outliers are empty sets.
\end{remk}
\subsection{Proof of Berry--Esseen estimates 
(Theorem \ref{berry-esseen})}\label{proof:BE_gen}
In this section, we present a proof of 
Theorem \ref{berry-esseen}. First, we will prove two
results (Lemmas \ref{lem_ineq} and \ref{lem_esseen})
that show  how much the distribution of a random variable
changes when we add a small perturbation.

\begin{lem}\label{lem_ineq}
Let $\xi$ be a real random variable defined on some probability 
space $(\Omega,\mathcal F, \Pr)$, let $\mu$ be a probability measure
on $(\mathbb R, \mathcal B(\mathbb R))$. Let $h$ be 
a bounded Borel  measurable function in $\mathbb R$,
 and $B\in\mathcal F$. Then
\beqnl{gen_ineq}
\left|\Ep\left[h\left(\xi\uno_B\right)\right]- \mu h\right|
\leq
2\|h\|_{\infty}\Pr[B^c]+\left|\Ep\left[h(\xi)\right]-\mu h\right|
\eeqnl
where $\mu h= \int h(x) \mu(dx)$.
\end{lem}
\proof
From the identity
$$ h\left(\xi\uno_{B}\right)=h(\xi)\uno_B + h(0)\uno_{B^c}$$
we have
$$
\left|\Ep[h(\xi)] - \Ep\left[h\left(\xi\uno_B\right)\right]
\right|
 \leq 2\|h\|_\infty\Pr[B^c]
$$
We obtain \eqref{gen_ineq} by the triangle inequality
\endproof
\begin{lem}\label{lem_esseen}
Let $\xi$ and $\eta$
be real random variables
defined in some probability space $(\Omega, \mathcal F, \Pr)$.
Assume that for some  Lipschitz function $G$ on $\mathbb R$ and 
some constant $0<a < 1$ we have
\begin{itemize}
\item[1)] $
\sup_{z\in\mathbb R}\left|\Pr[\eta\leq z]-G(z)\right|
\leq a $
\item[2)] $
\Ep[|\xi-\eta|] \leq  a^2
$
\end{itemize}
Then
\beqn
\sup_{z\in\mathbb R}\left|\Pr[\xi\leq z]-G(z)\right|
\leq (2+ L)a 
\eeqn
where $L=L(G)$ is the Lipschitz constant of $G$.
\end{lem}
\proof
Define the event $A=[|\xi-\eta|< a]$.  Then by Chebyshev's inequality,
we have $\Pr[A^c]\leq a$.

Partitioning the probability  space $\Omega$, we have
\beqn
\Pr[\xi\leq z]&=&\Pr[\{\xi\leq z\}\cap A]+ \Pr[\{\xi\leq z\}\cap A^c]\\
&\leq& \Pr[\eta\leq z+a]+ a
\eeqn
Adding and subtracting $G(z+a)-G(z)$ we obtain
\beqnl{wsup}
\sup\left\{ \Pr[\xi \leq z]- G(z)\right\} \leq (2+L)a
\eeqnl

Similarly, by partitioning the space we have
\beqn
\Pr[\eta\leq z-a] &=& \Pr[\{\eta\leq z-a\}\cap A] + 
\Pr[\{\eta\leq z-a\}\cap A^c]\\
&\leq& \Pr[\xi\leq z] + a
\eeqn
Adding and subtracting $G(z)-G(z-a)$ we obtain
\beqnl{winf}
-a(2+L)\leq \inf_{z\in\mathbb R}\left\{\Pr[\xi\leq z]-G(z)\right\}
\eeqnl
We finish the proof by combining \eqref{wsup} and \eqref{winf} 
\endproof
\noindent

Notice that from H\"{o}lder's inequality,
if the Lyapunov condition \eqref{lyapcon} holds at a point
$x\in M$ for some number  $2<p$, then  it also holds for
$s=\min(3,p)$. By the same token, if $(\xi_n)$ is a sequence
 of independent random variables with  
$p$ finite moments that satisfies A3($p$), then it satisfies
A3($s$).

If we denote by   
$$a_k= \frac{\Lambda_s(x,n_k)}{(\Lambda_2(x,n_k))^{s/2}},$$
 and let $b_k$ be the left--hand side of  \eqref{weakernoise}.
The assumption \eqref{weakernoise} is  expressed as
 $$b_k\leq a^2_k$$

Recall the processes $w_{n_k}$ and $l_{n_k}$ defined in
 \eqref{normprocess}  and \eqref{linclt} respectively.
By Theorem \ref{main}, we know that there are events $B_k$ such that
$$\max\left\{\Pr[B^c_k],\Ep[(w_{n_k}-l_{n_k})\uno_{B_k}]\right\}
\leq C a_k
$$

Using  Lemma \ref{lem_ineq} with $h=\uno_{(-\infty,z]}$ and 
$\mu(dt)=\Phi'(t)dt$, we get
\beqn
\sup_{z\in\mathbb R}\left|\Pr\left[l_k\uno_{B_k}\leq z\right]
-\Phi(z)\right|
\leq 2b_k + a_k<\tfrac32 a_k
\eeqn
for $k$ large enough. Then, the conclusion of
Theorem \ref{berry-esseen} follows from Lemma \ref{lem_esseen}
by letting $\xi=w_k\uno_{B_k}$, $\epsilon=l_{n_k}\uno_{B_k}$, 
$G(z)=\Phi(z)$ and $a=3a_k/2$.
$\hfill{\Box}$
\begin{remk}
Notice that the sequence of sizes of noises $\{\sigma_k\}_k$ and the
sequence of  outliers $\{B_k\}_k$ depend on the initial
condition $x_0$ of the process $x_n$. Therefore, the central limit
theorem in
Theorem \ref{main} and the Berry--Esseen estimates in
 Theorem \ref{berry-esseen} are  point wise 
results. In section \ref{examples},  we will give  an example where the 
convergence to central limit theorem holds uniformly 
with respect all initial conditions.
\end{remk}
\subsection{Cumulants}\label{secgencum}
The Lyapunov functions $\Lambda_s$ used in  Theorem \ref{main} have
  a simple interpretation  in terms  of  rather well 
known statistical quantities 
called cumulants in the statistical literature 
or Wick-ordered moments in statistical mechanics.

The cumulants, for a given
random variable $\xi$,
are formally defined as the coefficients $K_p$ in the power
series
\beqnl{cumulantes}
\log \Ep\left[e^{it\xi}\right]= \sum^\infty_{p=0} \frac{K_p}{p!}(it)^p
\eeqnl
Notice that $K_2[\xi]=\var[\xi]$. 
It can be seen \cite[Ch. IV]{petrov:sumran} that
 $K_p[\xi]$  are  homogeneous polynomials 
in moments, of degree $p$. Furthermore, for any $a\in\mathbb R$, 
and any pair of independent random variables $\xi$, $\eta$
we have that
\begin{equation}
\begin{array}{rcl}
K_p [a \xi] & = & a^p K_p [\xi]\\
\noalign{\vskip6pt}
K_p [\xi+\eta] & = & K_p [\xi] + K_p [\eta]
\end{array}
\label{properties}
\end{equation} 
The cumulants  $K_3[\xi]/(K_2[\xi])^{3/2}$ 
and \ $K_4[\xi]/(K_2[\xi])^2$ \ of $\bar{\xi}=\xi/\sqrt{K_2[\xi]}$, 
called {\em skewness} and {\em kurtosis} respectively,
 are empirical measures of  resemblance to Gaussian: the closer
 to zero the closer $\bar{\xi}$ is to a Gaussian \cite[p. 78]{mood}.
Observe that these measures are scale invariant.

\begin{defn}
For each positive integer $p$, we consider the 
cumulant functional
 $\widetilde{\Lambda}_p$  defined by
\beqnl{cumsuma}
\widetilde{\Lambda}_p(x,n)=\sum_{j=1}^n\left(\left(f^{n-j}\right)^\prime
\circ f^j(x)\right)^p
\eeqnl
\end{defn}

Observe that  $|\widetilde{\Lambda}_p(x,n)|\leq \Lambda_p(x,n)$ and that
if  $p$  even, $\widetilde{\Lambda}_p=\Lambda_p$. Furthermore, 
the cumulant functionals $\widetilde{\Lambda}_p(x,n)$ satisfy the analog 
to \eqref{lambdarelation} for the Lyapunov functionals, namely
\beqnl{wlambdarelation}
\widetilde{\Lambda}_p(x,n+m)
=\left((f^m)^\prime\circ f^n(x)\right)^p\widetilde{\Lambda}_p(x,n) 
+ \widetilde{\Lambda}_p(f^n(x),m)
\eeqnl 
In section \ref{cumsec} we will use \eqref{wlambdarelation} to
study the asymptotic behavior of the cumulants in the case of
systems near the accumulation of period doubling.

Notice that when $(\xi_n)_n$ is an i.i.d sequence with
$p\in\mathbb N$ finite moments, we have that
\beqnl{cumulantoflinnoise}
 K_s[L_n(x)] = K_s[\xi_1]\widetilde{\Lambda}_s(x,n)
\eeqnl
for all $0\leq s\leq p$.

If $(\xi_n)_n$ is a sequence of random variables that 
satisfies A3($p$)  for an integer  $p>2$, then from 
 \eqref{c1} and \eqref{properties} 
\beqnl{cumulantoflinearpart}
\left|K_p\left[\frac{L_n(x)}{\sqrt{K_2[L_n(x)]}}
\right]\right|
\leq
c_p \left|
\frac{\widetilde{\Lambda}_p(x,n)}{
\left(\Lambda_2(x,n)\right)^{p/2}}\right|\leq
c_p \frac{\Lambda_p(x,n)}{\left(\Lambda_2(x,n)\right)^{p/2}}
\eeqnl
for some constant $c_p>0$. Therefore,  condition \eqref{lyapcon}
 in Theorem \ref{main}
can be interpreted as a measure of closeness to Gaussian.  

Using renormalization theory, we will obtain sharp asymptotics  for
\eqref{cumulantoflinearpart} for maps near the 
accumulation of period doubling and critical circle maps with 
golden mean rotation number.
\subsection{Examples} \label{examples}
In this section, we consider a nontrivial example of a system to
which Theorem \ref{main} applies. We also show that
systems with enough hyperbolicity do not have a central limit theorem
in the direction of Theorem \ref{main}.
\subsubsection{Smooth diffeomorphisms of the circle}
We start this section by showing a central limit theorem for conjugate maps.
We will assume that $M=I$ or $\mathbb T^1$. 
\begin{lem}\label{conjugate}
Let $f,g\in C^2(M)$, and let $(\xi_n)$ be as in \ref{main}.
Assume $$f=h^{-1}\circ g\circ h$$
 for some map $h$ such that $h, h^{-1}\in C^1(M)$. 
The Lyapunov condition \eqref{lyapcon} holds for  $f$ 
at some point $x$ if and only if it  holds for  $g$ at $h(x)$.
\end{lem}
 \proof
Notice that 
$f^j=h^{-1}\circ g^j\circ h$. Since $h,\, h^{-1}\in C^1(M)$
and $M$ compact, we have that 
$$c \left|(g^j)^\prime\circ h(x)\right|\leq 
\left|(f^j)^\prime(x)\right|\leq C
 \left|(g^j)^\prime\circ h(x)\right|$$
for some constants $C,c >0$. If we denote $y=f(x)$, we have
\beqn
c^s \Lambda^g_s(y,n)\leq \Lambda^f_s(x,n)\leq 
C^s \Lambda^g_s(y,n)
\eeqn
for all $0<s\leq q$. Now,  Lemma \ref{conjugate} follows immediately.
\endproof

Recall that every orientation--preserving homeomorphism $f$
of the circle $\mathbb T^1$ is given by $f(x)=F(x) \mod1$, where
$F$ is a strictly increasing continuous function such that
$F(x)=f(x)+1$ for all $x\in\mathbb R$.  
For any  $\alpha\in\mathbb R$, we denote by $R_\alpha(x)=x+\alpha\mod1$
the rotation of the circle with rotation number $\alpha$. 

We will  use Lemma \ref{conjugate} to show that Theorem \ref{main}
applies for smooth diffeomorphisms of the
circle with Diophantine rotation number.
 This will be a consequence of the following Herman type theorem:
\begin{thm}\label{herman}
Let $k>2$ and assume that  $f\in C^k(\mathbb T^1)$ is an
orientation preserving diffeomorphism, whose
 rotation number $\alpha$ is a Diophantine number
with exponent $\beta$. That is
 \beqn
\left|\alpha -\tfrac{p}{q}\right| > \tfrac{C}{q^{2+\beta}}
\eeqn
for all $\tfrac{p}{q}\in \mathbb Q$, where $C>0$ is a fixed
constant. If $\beta<k-2$, then there is a map
 $h\in C^{k-1-\beta}$\, such that
$$f=h^{-1}\circ R_\alpha\circ h$$
\end{thm}
Theorem \ref{herman}, due to \cite{sinkhan1}, is the most comprehensive
 global result on the problem
of smoothness of conjugacies of diffeomorphisms of
 the circle with rotations. 
The first global result on this problem 
is due to \cite{herman} for $k\geq3$,
and generalized later by \cite{yocc:diocirc}.
 The first mayor improvement ($k>2$)
in the solution of this problem was
obtained by \cite{sinkhan2,katzor,sinkhan1}.
 A brief historical account
of the problem and a detail proof of this and related results 
can be found in \cite{sinkhan1}. 
The technique developed in \cite{sinkhan1}  constitutes 
 one of the most important applications of the  
renormalization group method in dynamical systems.

\bigskip
\noindent
{\em Example 1}.
Notice that Theorem \ref{main} applied to $R_\alpha$ is the
classical central limit theorem.  Thus, we can  conclude that
any Diophantine diffeomorphism of the circle has
a  Gaussian scaling limit for weak noises along the whole
sequence of integers.
Moreover, if the sequence of random variables $\xi_n$ 
have finite moments of order $p=3$ and we let
$\sigma_n=n^{-3/2}$, we can define the outliers by
$$B^c_n=\left[\|f''\|_{C_0}
\max_{i\leq l\leq n}|\xi_l|\geq \sqrt{n} \right]
$$  
Notice that the outliers are independent of the initial
point of the orbit. Furthermore, 
by Theorem \ref{berry-esseen}  we have
\beqn
\sup_{(x,z)\in\mathbb T^1\times \mathbb R}\left|
\Pr[w_n(x)\leq z]-
\Phi(z) \right|\leq \frac{A}{\sqrt{n}}
\eeqn
for some constant $A>0$.
$\hfill{\Box}$

\subsubsection{Systems with positive Lyapunov exponents.}
In this section, we give an example of a system that does not satisfy the
Lyapunov condition \eqref{lyapcon}, and for which the conclusion
of Theorem \ref{main} fails.

\bigskip
\noindent
{\em Example 2.} Consider the map 
on ${\mathbb R}$ (or ${\mathbb T}^1$) given by 
$$f(x)=2 x$$
Condition \eqref{lyapcon}  is not satisfied,
and as we will see the conclusion of Theorem \ref{main} fails.
We will show that
there is a limit for the scaled noise, not necessarily Gaussian, 
 which depends strongly on the distribution of the original noise
$\xi_n$.  Notice that
\beqn
w_n &=&\frac{x_n(x_0,\sigma_n)-2^nx_0}
{\sigma_n\var\left(\sum^n_{j=1}2^{n-j}\xi_j\right)}\\
&=&\frac{3\sqrt{2}}{2\sqrt{1-4^{-n}}}\sum^n_{j=1}2^{-j}\xi_j
\eeqn
If $(\xi_n)$ is an i.i.d. sequence 
with uniform distribution $U[-1,1]$ then, $w_n$ will
converge in law to a compactly supported random variable $\xi$
with characteristic function
$$\phi(z)=\prod\limits^\infty_{k=2}
\frac{2\sqrt{2}\sin(2^{-k}\sqrt{2} 3z)}{2^{-k}3z}$$
On the other hand,  if $(\xi_n)$ is an i.i.d sequence with
 standard normal
distribution, then $w_n$ has standard normal distribution for all
times.

Similar results hold for hyperbolic orbits and the reason is that
derivatives at hyperbolic points grow (or decay) exponentially, that
is $\left|(f^{n-j})^\prime\circ f^j(x)\right|\approx a^{n-j}$
 for some number $a>0$. If  $(\xi_n)$ is an  i.i.d sequence 
for instance,  then  the  cumulants $K_p$ 
of order $p> 2$ of the normalized variable $L_n$
do not decay to zero. In fact we have
\beqn
K_p\left[
\frac{L_n(x)}{\sqrt{\var[L_n(x)]}}\right]
\approx K_p[\xi_1]\sqrt{|a^2-1|^{p-2}}
\eeqn
where $K_p$ is the cumulant of order $p$ of $\xi_1$. 
Therefore, for hyperbolic orbits
the scaling limit depends strongly on the distribution of
the sequence $\xi_n$.
$\hfill{\Box}$

In view of the example above, we can see that the scaling version of
the central limit theorem we have proved does not apply for 
maps which are hyperbolic. In that respect, it is curious 
to mention that a celebrated result \cite{GraczykS96,Lyubich97,GraczykS97}
shows that the  maps in the quadratic family are hyperbolic 
for a dense set of parameters.

\section{Central limit theorem for maps the accumulation of period
 doubling}\label{perioddoubling}
In this section,  we will show that the orbits \eqref{ransys} at
 the accumulation of period--doubling satisfy 
a central limit theorem (Theorem \ref{fclt}). We consider some 
statistical characteristics of the effective noise,
i.e. Lindeberg--Lyapunov sums,
 and show that they satisfy some scaling relations
(Corollary \ref{kurto2n}). We use this scaling relations to show that
 the Lyapunov condition \ref{lyapcon} holds for certain
class of initial conditions (Theorem \ref{kurtnx}).
\subsection{The period doubling renormalization
group operator}\label{symmunimodalclt}
We will consider analytic even maps $f$ of
the interval $I=[-1,1]$ into itself such that
\begin{itemize}
\item[P1. ] $f(0)=1$
\item[P2. ]  $xf'(x)<0$ for $x\neq0$
\item[P3. ] $Sf(x)<0$ for all $x\neq0$ where $S$ is the Schwartzian 
derivative ($Sf(x)=(f'''(x)/f'(x))-(3/2)(f''(x)/f'(x))^2$)
\end{itemize}
The set of functions that satisfy conditions P1--P3  is called the set of 
unimodal maps.\\
Let us denote by
$\lambda_f=f(1)$ and $b_f=f(\lambda_f)$. We will further assume that
\begin{itemize}
\item[P4.] $0<|\lambda_f|<b_f$ 
\item[P5.]  $0<f(b_f)<|\lambda_f|$
\end{itemize}
The set $\mathcal P$ of functions that satisfy 
P1--P5 is called the space of period--doubling renormalizable
functions.

It follows from P4 that  the intervals 
$I_0=[-|\lambda_f|,|\lambda_f|]$ and $I_1=[b_f,1]$ do not overlap. 
If $f\in{\mathcal P}$, then $f\circ f|_{I_0}$ has a single critical point, 
which is a minimum. The change of variables $x\mapsto \lambda_fx$ 
replaces $I_0$  by $I$ and the minimum by a maximum.
\begin{defn}
 The period doubling 
renormalization operator $T$  on $\mathcal P$ is defined by
$$Tf(x)=\frac{1}{\lambda_f}f\circ f(\lambda_f x)$$
If $f,\, Tf,\ldots,T^{(n-1)}f\in\mathcal P$, we 
 denote  by $\Gamma_n=f^{2^n}(0)$ (or  
 $\Gamma^f_n$ if we need to emphasize the dependence on $f$).
We have that
\beqn
T^nf(x)&=&\frac{f^{2^n}(\Gamma_nx)}{\Gamma_n}\\
\lambda_{T^nf}&=&\frac{\Gamma_{n+1}}{\Gamma_n}
\eeqn
\end{defn}
The renormalization group operator $T$ has been well studied by many 
authors, \cite{feig:univ,Coullet-Tresser,coeck:itemap,
coecklan:univ,lan:casproof,vul:feiguniv,eckwitt,
eps,colles:renno,sull,mart,mel,JacobsonS02}.
 The following results, whose proofs can be found in the reference
above, will be very useful for our purposes.
 \begin{itemize}
\item[R1.] For each $k\in\mathbb N$, there is an analytic 
 function $h_k$ such that $g_k(x)=h_k(x^{2k})$
is analytic on some open bounded domain $D_k\subset\mathbb C$\,
containing  the interval $[-1,1]$, $g_k$ restricted to $I$ is
concave, and for all $x\in I$
\beqnl{feigfixpointeq}
Tg_k(x)=g_k(x)
\eeqnl
Furthermore, since $Sf(x)<0$ for all $x\neq0$, 
it follows that  $g''(x)<0$ for all $x\neq0$.
\end{itemize}
Notice that \eqref{feigfixpointeq} implies that
$$g_k'(x)=g'_k(g_k(\lambda_k\,x))g'_k(\lambda_k x)$$
Then, by taking the limit $x\rightarrow0$ we get
$$
g'_k(1)=\lim_{x\rightarrow0}\frac{g_k'(x)}{g'_k(\lambda_k x)}=\lambda_k^{1-2k}
$$
\begin{itemize}
\item[R2.] The domain $D_k$ can be chosen so that $\partial D_k$ is
smooth, $\overline{\lambda_k\,D_k}\subset D_k$ and
$\overline{g_k(\lambda_k\,D_k)}\subset D_k$.
\end{itemize}
This condition is closely related to the fact that $g_k$ is in the domain 
of the period doubling renormalization operator.

For each $k$, denote by $H(D_k)$ the space of analytic functions on
$D_k$, and let  $\mathcal G_k\subset H(D_k)$ be the Banach space
 of analytic functions
\beqn
\mathcal G_k=\left\{f:\, f(z)=z^{2k}\hat{f}(z),\, \hat{f}(z)\in H(D_k),\,
\hat{f}\in C_0(\overline{D_k})\right\}
\eeqn
endowed with the sup norm. 
Notice that if  $f\in \mathcal G_k$, then
 $f^{(j)}(0)=0$ for all derivatives of order $0\leq j<2k$. 
\begin{itemize}
\item[R3.] There
 is a  neighborhood $\mathcal V$ of $g_k$ in $(\mathcal G_k+1) \cap\mathcal P$
   where  $T$ is differentiable.
\item[R4.] A consequence of R2 is that if $f\in\mathcal V$,
 the derivative $DT(f)$   is a compact operator on $\mathcal G^k$.
\item[R5.] $T$ admits   an invariant stable manifold $\mathcal W_s(g_k)$, 
such that if $f\in\mathcal W_s(g_k)$ then, 
\beqnl{rgconvergence}
\lim_{n\rightarrow\infty}\|T^n f - g_k\|_{C(\bar{D}_k)}=O(\omega_k^n)
\eeqnl
where $0<\omega_k<1$ is a universal constant.
\end{itemize}
 The set $\mathcal W_s(g_k)$
is called the universality class of $g_k$.

If we denote by $\lambda_{g_k}=\lambda_k$ then,  for each 
$n\in\mathbb N$ we have that
$\Gamma^{g_k}_n=(\lambda_k)^n$. Moreover,
the exponential convergence (R4) implies  that
for any $f\in\mathcal W_s(g_k)$,\, $\lambda_{T^n}$ converges to
$\lambda_k$ exponentially fast.

\subsection{Renormalization theory of the noise}\label{cumulantoper}
In this section, we develop a 
theory of renormalization for the noise.
Since the renormalization group operator $T$  gives information for times
 $2^n$ at small scales,  we will study first the effective noise 
at times $2^n$ for  orbits $x_n$  starting at zero.

Let $f\in\mathcal W_s(g)$ be fixed, and for each  $p\geq0$ denote by 
$k^{(n)}_p(x)=\Lambda^f_p(x,2^n)$. From  \eqref{lambdarelation} we have 
\beqnl{cumopfeig2n}
k^{(n)}_p(x)= \left|
\left(f^{2^{n-1}}\right)^\prime\circ f^{2^{n-1}}(x)\right|^p
k^{(n-1)}_p(x) +
k^{(n-1)}_p\circ f^{2^{n-1}}(x) 
\eeqnl
Since  $f\in\mathcal P$, then  $-f^\prime(x) >1>0$ for all
 $x\in[f(\lambda_f),1]$. Therefore, by taking a
determination of   logarithm
defined on
$$\mathbb C\setminus \{x+iy: x\leq0,\, y=0\}$$
we can analytically define the map
 $z\mapsto (-f^\prime\circ f(\lambda_f z))^p$ 
 on $ D_k$.  In particular, notice that if $z\in D_k\cap\mathbb R$, then 
\beqnl{observation}
(-f^\prime\circ f(\lambda_f z))^p=
\left|f^\prime\circ f(\lambda_f z)\right|^p
\eeqnl
\subsubsection{Lindeberg--Lyapunov 
operators}\label{linlyapoperatorfeig}
The following family of operators will be very useful in
the study of the propagation of noise.
\begin{defn}
For each $f\in\mathcal W_s(g)$, the family of Lindeberg--Lyapunov operators
$\{\mathcal K_{f,p}: p\geq0\}$  acting  on the 
space $H^r(\bar{D}_k)$ of real and bounded analytic
 functions on $D_k$ which are continuous in $\bar{D}_k$, 
is defined by
\beqnl{cumop}
\mathcal K_{f,p}h(z)=\frac{1}{(-\lambda_f)^p}\left[
\left(-f'\circ f(\lambda_f z)\right)^p h(\lambda_f z)
+ h(f(\lambda_f z))\right]
\eeqnl 
where $h\in H^r.(\bar{D}_k)$.
\end{defn}
The operators defined by \eqref{cumop} belong to a class of operators
 called transfer operators, which have been extensively studied
\cite{baladi,mayer}.

It follows from \eqref{cumopfeig2n} and \eqref{cumop} that if
\beqn
x=\Gamma^f_n z,\qquad
\widetilde{k}^{(n)}_p(z)=
\left|\Gamma^f_n\right|^{-p}k^{(n)}_p(\Gamma^f_n z),
\eeqn
 then,  defining the function $\uno$ by  $\uno(z)=1$
\beqnl{feig_smallscale}
\widetilde{k}^{(n)}_p(z)
=
(\mathcal K_{T^{n-1}f,p}
\cdots\mathcal K_{f,p} \uno)(z)
\eeqnl
Equation \eqref{feig_smallscale} provides a relation between the
 Lindeberg--Lyapunov  operators
defined by \eqref{cumop} and the effective noise. Indeed,
 recall that the sequence of noises
 $(\xi_n)_n$ satisfies 
$c^{-1}<\inf_n\|\xi\|_2\leq \sup_n\|\xi\|_p<c$. Then,
for $s=2,\,p$ we have that
\beqn
c^{-1}\widetilde{k}^{(n)}_s(z)
\leq
|\lambda_k|^{n\,s}\sum_{j=1}^{2^n}\left|\left(f^{2^n-j}\right)^\prime
\circ f^j(\lambda^n_k z)\right|^s\|\xi\|^s_s
 \leq c\,\widetilde{k}^{(n)}_s(z)
\eeqn
This means 
that \eqref{feig_smallscale}
estimates the growth of  the linearized  propagation of noise
at times $2^n$  with initial condition at the origin.

\subsubsection{Exponential convergence of the Lindeberg--Lyapunov
operators}
Recall that for a map $f\in \mathcal W_s(g_k)$, we have that
$T^nf(z)$ converges exponentially fast to $g_k$. Our next result
implies that  the sequence of Lindeberg--Lyapunov 
operators $\mathcal K_{T^nf,p}$
also converges exponentially to $\mathcal K_{g_k,p}$. 
This implication will be
important in Section \ref{sec31}, Proposition 
\ref{productalignment}, where
we analyze the asymptotic behavior of the product of operators
in \eqref{feig_smallscale}.

\begin{lem}\label{Kumconvergence}
Let $f\in\mathcal W_s(g_k)$ close enough to $g_k$. For
all $p\geq0$,
$\mathcal K_{f,p}$ is a compact operator on $H^r(\overline{D}_k)$.

Furthermore, there is a constant $C_p$ such that
\beqnl{expconvergence}
\|\mathcal K_{f,p} -\mathcal K_{g_k,p}\|< C_p\|f- g_k\|_{C(\bar{D}_k)}
\eeqnl
\end{lem}
\proof
For any set $D\subset \mathbb C$, let 
 $D^\epsilon=\{z:d(x,D)\leq \epsilon\}$.
Recall that 
\beqn
\overline{\lambda_k D_k}&\subset& D_k\\
\overline{g_k(\lambda_k D_k)}&\subset& D_k
\eeqn
Since for all $z\in \bar{D}_k$
$$|\lambda_f\,z-\lambda_k\,z|\leq \|f-g_k\|_{\bar{D}_k}\diam(D_k),$$
we can  choose a neighborhood $\mathcal U$ of $g_k$ such that for all
$f\in\mathcal W_s(g_k)\cap\mathcal U$,
\beqnl{setrel1}
\overline{\lambda_f D_k}\subset \overline{\lambda_k D_k}^\epsilon
\subset D_k
\eeqnl
where $\epsilon > 0$ is small enough.
Let us denote by $V_k=\overline{\lambda_k D_k}^\epsilon$.

Notice that for all $z\in D_k$,
\beqn
|g_k(\lambda_k\,z)-g_k(\lambda_f\,z)|&\leq &
\|g'_k\|_{V_k}\diam(D_k)\|f-g_k\|_{D_k}\\
|f(\lambda_f\,z)-g_k(\lambda_f\,z)|&\leq&
\|f-g_k\|_{V_k}
\eeqn
Then, for all $z\in D_k$ we have that
\beqnl{firstgoodbound}
|f(\lambda_f\, z)- g_k(\lambda_k\,z)|\leq C\|f-g_k\|_{D_k}
\eeqnl
for some constant $C$. Therefore, shrinking $\mathcal U$ if necessary,
 we can assume that
\beqnl{setrel2}
\overline{f(\lambda_f D_k)}\subset\overline{g_k(\lambda_k D_k)}^\epsilon
\subset D_k
\eeqnl
 Let us denote by 
$W_k=\overline{g_k(\lambda_k D_k)}^\epsilon$. 

For all $z\in D_k$ we have that
\begin{multline}
|f'(f(\lambda_f\,z))-g'_k(g_k(\lambda_k\,z))\|\leq 
|g'_k(g_k(\lambda_k\,z))-g'_k(f(\lambda_k\,z))|\nonumber\\+
|g'_k(f(\lambda_f\,z))-f'(f(\lambda_f\,z))|
\end{multline}
Using Cauchy estimates and the following inequalities
\beqn
|g'_k(g_k(\lambda_k\,z))-g'_k(f(\lambda_k\,z))|
&\leq&\|g''_k\|_{W_k}\|g_k-f\|_{V_k}\\
|g'_k(f(\lambda_f\,z))-f'(f(\lambda_f\,z))|&\leq&
\|g'_k-f'\|_{W_k},
\eeqn
we obtain that
\beqnl{compositionbound}
\|f'\circ f\circ\lambda_f-g'_k\circ g_k\circ\lambda_k\|_{\bar{D}_k}
\leq \tilde{C}\|g_k-f\|_{\bar{D}_k}
\eeqnl
Let $p>0$ be fixed. Notice that the Lindeberg--Lyapunov
 operator $\mathcal K_{f,p}$ has the form
$$\mathcal K_{f,p} h = U_f h + R_f h$$
where 
$U:h\mapsto (\psi\circ f'\circ f\circ\lambda_f)(h\circ\lambda_f)$
and $R:h\mapsto h\circ f\circ\lambda_f$, where $\psi(z)=z^p$ is 
defined on $O=\mathbb C\setminus\{x+i\,y: x<0,\,y=0\}$. By
shrinking $\mathcal U$ if necessary,
we can assume by \eqref{compositionbound} that 
if $f\in\mathcal W_s(g_k)\cap \mathcal U$ then,
$$
\overline{f'(f(\lambda_f D_k))}\subset
 \overline{g'_k(g_k(\lambda_k D_k))}^\epsilon\subset O$$
Let us denote by $Y_k= \overline{g'_k(g_k(\lambda_k D_k))}^\epsilon$
From \eqref{firstgoodbound}, we have that
\begin{equation}
\begin{array}{lcr}
\|(R_f-R_{g_k})h\|_{D_k}&\leq& C\|h'\|_{W_k}\|f-g_k\|_{\bar{D}_k}\\
\noalign{\vskip6pt}
&\leq& \hat{C}\|h\|_{\bar{D}_k}\|f-g_k\|_{\bar{D}_k}
\end{array}
\label{StepI}
\end{equation}
From \eqref{compositionbound}, we have that
\begin{multline}
\|(U_f-U_{g_k})h\|\leq 
\hat{C}\|\psi\|_{Y_k}\|h\|_{\bar{D}_k}\|f-g_k\|_{\bar{D}_k}\label{StepII}\\
+ \tilde{C}\|\psi'\|_{Y_k}\|h\|_{V_k}\|f-g_k\|_{\bar{D}_k}
\end{multline}
The exponential convergence follows
from \eqref{StepI} and \eqref{StepII}.

\bigskip
Notice that the operators $\mathcal K_{f,p}$ send bounded
sets in $H^r(\overline{D}_k)$ into bounded sets in 
$H^r(\overline{D}_k)$. 

Furthermore, for all $z_1,\,z_2 \in D_k$ we have
from \eqref{setrel1} and \eqref{setrel2} that
\beqn
| R_fh(z_1)- R_fh(z_2)|&\leq& 
\|(h\circ f)^\prime\|_{V_k}|z_1-z_2|\\
|U_fh(z_1)-U_fh(z_2)| &\leq& 
\left|\left((\psi\circ f)h\right)^\prime\right|_{V_k}|z_1-z_2|
\eeqn
The compactness of the operators $\mathcal K_{f,p}$ follows
from Cauchy estimates and the Arzela-Ascoli theorem.
\endproof

We will study  the spectral properties of the Lindeberg--Lyapunov
 operators  to show that  \eqref{lyapcon} holds, and then  
by  Theorem \ref{main}  we will obtain the Gaussian scaling 
limit for systems at the accumulation of period--doubling.
\subsection{Spectral analysis of the Lindeberg--Lyapunov 
operators}\label{sec31} 
Let $\mathcal C$ be the cone of functions
 on $H^r(\bar{D}_k)$ which are non-negative when restricted
to $\bar{D}_k \cap \real$.  Note that 
$\text{int} \mathcal C$, the interior of $\mathcal C$ consists of 
the functions in  $H^r(\bar{D}_k)$
which are strictly possible when restricted to 
$\bar{D}_k \cap \real$.

Let $f\in\mathcal W_s(g)$
and $p\geq0$ fixed, and denote $\mathcal K_p=\mathcal K_{f,p}$.

We have that  
\begin{itemize}
\item[(i)] $\mathcal K_p(\mathcal C\setminus\{0\})
 \subset \mathcal C\setminus\{0\}$
\item[(ii)] $\mathcal K_p(\text{int}(\mathcal C))\subset 
\text{int}(\mathcal C)$, and 
\item[(iii)] for each  $h\in \mathcal C\setminus\{0\}$, there is an 
integer $n=n(h)$ such that
 $\mathcal K^n_p h \in \text{int}(\mathcal C)$.
\end{itemize}

Properties (i), (ii) are quite obvious. Property (iii) follows from 
the observation that if a function is non-negative and 
strictly positive in an interval $I^*$ -- the set of zeros 
has to be isolated --  then,  
by the definition of ${\mathcal K}_p$ \eqref{cumop}
it is strictly positive in 
the image of $I^*$ under the inverse of the maps $\lambda_f \cdot $,
$f(\lambda_f \cdot))$. Since the functions $\lambda_f, f(\lambda_f \cdot ) $
are have derivatives strictly smaller than $1$, it follows
that a finite number of iterations of 
the inverse functions covers the whole 
interval. 

By  Kre\v{\i}n--Rutman's Theorem  \cite[p. 265]{schaefer},
\cite{tack:kreinrut}, we have that
\begin{prop}\label{KR}
Let $f\in\mathcal W_s(g)$ and $p\geq0$, and denote
 by $\mathcal K_p=\mathcal K_{f,p}$. Then,
\begin{itemize}
\item[K1.] The spectral radius $\rho_p$ of $\mathcal K_p$ is a positive
 simple  eigenvalue of $\mathcal K_p$.
\item[K2.] An eigenvector $\psi_p\in H^r(\bar{D}_k)\setminus\{0\}$ 
associated with $\rho_p$ can be chosen in $\text{int}(\mathcal C)$.
\item[K3.] If $\mu$ is in the spectrum of $\mathcal K_p$,
 $0\neq \mu\neq \rho_p$, 
then $\mu$ is an eigenvalue of $\mathcal K_{f,p}$ and $|\mu|<\rho_p$.
\item[K4.] If $h\in \mathcal C\setminus\{0\}$ is an eigenvector of 
$\mathcal K_p$, then the corresponding eigenvalue  is $\rho_p$.
\end{itemize}
\end{prop}

\begin{remk}
We emphasize that Krein-Rutman's theorem only needs to assume 
that the operators are compact and that preserve a cone. Hence, 
the only properties of the operators ${\mathcal K}_p$ that 
are used are precisely $i),ii), iii)$ above. Later in 
Section~\ref{critcircmap} we will see another application of 
the abstract set up. In the case of Section~\ref{critcircmap}
the space will be a space whose elements are
pairs of functions. It is important to note that 
Corollary~\ref{cor:consequenceKR}, Proposition~\ref{uniqueness}
and  Proposition~\ref{KRbound} remain valid in the case of
several components. 
\end{remk}

We will say that two sequences of functions 
$\{a_n(z)\}_n$ and $\{b_n(z)\}_n$
of function in $H^r( \bar{D}_k) $
are asymptotically
 similar (denoted by  $a_n\asymp b_n$) if 
$$\lim\limits_{n\rightarrow\infty}\left\|a_n/b_n\right\|_{H^r(\bar{D}_k)} 
= 1$$

As standard, we use the notation $\psi > \phi$ 
to denote $\psi - \phi \in \text{int} \mathcal{C} $.  
In particular $\phi  > 0$ means $\phi \in \text{int} \mathcal{C} $.  
For function in $H^r( {\bar D}_k)$, $\phi > \psi$ is just means 
that $\phi(z) > \psi(z)$ for all $z \in {\bar D}_k \cap \real$. 

One property of the cones in the 
spaces we are considering in this Section and in Section
\ref{critcircmap} is that:

\begin{property}\label{ordercone}
Let $\phi, \psi \in \text{int} \mathcal{C}$ 
Then, there exists $\delta > 0 $ such that 
$\psi > \delta \phi$. 
\end{property}

A  consequence of Proposition \ref{KR}
is:
\begin{cor}\label{cor:consequenceKR}
 $h\in {\mathcal C}\setminus\{0\}$
 then,  there is a constant $c_p=c_p(h) > 0$ such that
\beqnl{consequenceKR1}
||{\mathcal K}_p^n h(z) - c_p \rho_p^n \psi_p(z)|| 
\le C (\rho_p - \delta)^n
\eeqnl 
for some $C >0, \delta > 0$. 

In the one-dimensional case, this implies. 
\beqnl{consequenceKR}
{\mathcal K}_p^n h(z) \asymp c_p \rho_p^n \psi_p(z)
\eeqnl
\end{cor}

To prove \eqref{consequenceKR1} we 
note that we can define spectral projections 
corresponding to the spectral radius and to the
complement of the spectrum. 
If $h = c_p  \psi_p + h^{<}$ with $h^{<}$ in the space
corresponding to the spectrum away from the spectral 
radius. 

By the spectral radius formula, we have
$ || {\mathcal K}_p^n h^{<}|| \le C (\rho_p -\delta)^n $
for some $\delta > 0$. On the other hand, we have 
that for some $\delta > 0 $, $h > \delta \psi_p$ 
over the reals. Hence, 
${\mathcal K}_p^n h > \delta \rho_p^n \psi_p$ over the reals. 
By comparing with 
${\mathcal K}_p^n h =  \delta \rho_p^n  \psi_p
+ {\mathcal K}_p^n h^<$, we obtain that $c_p > 0$. 

The conclusion 
\eqref{consequenceKR} follows from the fact that $\psi_p$ 
is bounded away from zero. 

Another important consequence 
of Proposition~\ref{KR} and the analyticity improving is 
that the spectrum and the eigenvalues 
are  largely independent of the
domain of the functions we are considering. This 
justifies that we have not included it in the notation. 

\begin{prop}\label{uniqueness}
Assume that $D_k$, $\tilde D_k$ are domains 
satisfying the assumptions $R_2$.

Assume that applying a finite number of
times the inverse functions of either 
$\lambda_k \cdot$ or $g( \lambda_k \cdot)$ 
$D_k$ we obtain a domain which contains 
$\tilde D_k$.
Then, the spectrum of
${\mathcal K}_{f,p}$ acting on $H^r(D_k)$ is contained
in the spectrum of ${\mathcal K}_{f,p}$ acting on 
$H^r(\tilde D_k)$. 
\end{prop}

We know that the spectrum is just eigenvalues (and zero). 
If we use the functional equation satisfied by an 
eigenvalue of ${\mathcal K}_{f,p}$, we see that the eigenfunction 
extends to the images of $D_k$ under 
the inverse images of 
$\lambda_k \cdot$ or $g( \lambda_k \cdot)$ 

Hence, by the assumption, we obtain that the 
eigenvalues in $H(D_k)$ are eigenvalues in $H( \tilde D_k)$.
\qed

Of course, if the hypothesis of Proposition~\ref{uniqueness}
is satisfied also when we exchange $D_k$ and $\tilde D_k$, 
we conclude that the spectrum is the same. 

Note that in our case, since the inverse maps are expansive, 
Proposition~\ref{uniqueness} shows that the
spectrum of the operator is independent of the domain, provided
that the domain is not too far away from the unit interval. 

For the case of the quadratic fixed point there is a very 
detailed study of the maximal domains in \cite{EpsteinL81}. 
The results of \cite{EpsteinL81} imply that there is a
natural boundary for the eigenfunctions of the operators. 
Of course, if we consider domains larger than that, the results are
different. Nevertheless, for domains inside this natural boundary, 
the spectrum is independent of the domain.

\subsubsection{Properties of the spectral radius
of the Lindeberg--Lyapunov operators}
For $f\in\mathcal W_s(g)$, we have that the spectral radii of
the operators $\mathcal K_{f,p}$  satisfy some  convexity
properties. 

\begin{thm} \label{Comparison}
Let $f$ be a map in ${\mathcal W}_s(g_k)$ close enough to $g_k$. For
each $p\geq0$, denote by $\mathcal K_p=\mathcal K_{f,p}$, and 
the spectral radius of $\mathcal K_p$ by
 $\rho_p$. Then
\begin{itemize}
\item[S1.] 
\beqnl{compa1}
(\lambda^{-1}_f f^\prime(1))^p < \rho_p < 
(\lambda^{-1}_f f^\prime(1))^p + (-\lambda_f)^{-p}
\eeqnl
In particular,  if $f=g_k$ then
\beqnl{radspecbound}
 1<\lambda_k^{2kp}\rho_p<1+ (-\lambda_k)^{(2k-1)p}
\eeqnl
\item[S2.]
For $p \ge 0$.
The map $p\mapsto \rho_p$ 
is strictly increasing and strictly log--convex.

\item[S3.]  The map  $p \mapsto\log[\rho_p]/p$ 
is strictly decreasing.

\end{itemize}
\end{thm}

For future purposes in this paper, [S3.] is the most important of
the consequences claimed in Theorem~\ref{Comparison}.

\proof
S1: For  fixed $p\geq0$, let  $\psi_p$ be a positive 
eigenfunction of $\mathcal K_p$.
Denote by  $b=(-\lambda_f)^{-p}$, and  notice that $b>1$.
Let $U_p$, $R$  and $S_p$ be the nonnegative compact operators defined by
\beqn
U_ph(z)&=&\left(-f^\prime\circ f(\lambda_f z)\right)^{p}h(\lambda_f z)\\
R h(z)&=&h\circ f(\lambda_f z)\\
S_p h(z)&=& (-f'(1))^p h(\lambda_f z)
\eeqn
Since $\mathcal K_p\,h= b\left(U_p + R\right)h$,
for
any $h\in \mathcal C\setminus\{0\}$ and any $n\in{\mathbb N}$
we have that, by the binomial 
theorem
\beqnl{binomial}
 \mathcal K^n_ph >  b^n(U_p^n + R^n)h > b^n\,U_p^n h.
\eeqnl

Notice that
\beqn
U_p^nh(z)=\left\{
\prod^n_{j=1}\left(- f^\prime\circ f(\lambda^j_f z)\right)\right\}^p
h(\lambda^n_f z)
\eeqn
If we  take  $h(z)=\psi_p(z)$  and then let $z=0$, 
the left--hand side of \eqref{compa1} follows.

Since  $1<-f'(z)<-f'(1)$ for all $z\in[f(\lambda_f),1]$,  we have that
$$\mathcal K^n_ph < b^n( S_p + R)^n h$$
positive function $h$. Notice that the spectral radius of the positive
 operator $S_p+R$ is  $(-f'(1))^p+ 1$. Then, 
the  right--hand side of \eqref{compa1} follows. 

The claim in [S1.] for  the  special case  $f=g_k$ follows  
from  
$$g_k^\prime(1)=\lambda^{1-2k}_{g_k}=\lambda_k^{1-2k}$$

\noindent

To prove [S2.], [S3.] we will use the following Proposition

\begin{prop}\label{KRbound} 
In the conditions of Proposition~\ref{KR}.
Assume that $\phi$ is a positive function 
$\lambda$ is a positive number such that 
\[
{\mathcal K}_p \phi > \lambda \phi
\]
Then, 
\[
\rho_p > \lambda
\]
\end{prop}

Note that because the functions are strictly positive, 
we have ${\mathcal K}_p \phi \ge  (\lambda + \delta)  \phi$ for 
some $\delta > 0$. Then, applying this repeatedly, we
obtain ${\mathcal K}_p^n \phi \ge  (\lambda + \delta)^n  \phi$. 
By Corollary~\ref{cor:consequenceKR}, we have
${\mathcal K}_p^n \phi = c_p \rho_p^n \psi_p + o((\rho- \delta)^n)$. 
Since the function $\psi_p$ is strictly positive, we have 
the desired result. 
\qed

Now, we continue with the proof of [S2.], [S3.]

We recall that the operator ${\mathcal K}_p$ has the structure:
\beqnl{structure} 
{\mathcal K}_p h = A^p h \circ f_1 + h \circ f_2
\eeqnl
where $A(z) = -f' \circ f(\lambda_f z)$, $f_1(z) =  
\lambda_f z$, $f_2(z) = f(\lambda_f z)$.  In particular, $A > 1$. 

The strict monotonicity of $\rho_p$ 
follows from the fact that if $q > p$
\[
\begin{split}
\rho_q \psi_q  &= A^q \psi_q \circ f_1 + \psi_q \circ f_2 \\
&>  A^{p} \psi_q \circ f_1 + \psi_q \circ f_2 \\
&= {\mathcal K}_p \psi_q
\end{split}
\]

To prove that $\log( \rho_p)/p$ is strictly decreasing,
we will show that  for any $0 < p$ and any $\alpha > 1$ we have
$\rho_{p \alpha} < \rho_p^\alpha$, which is equivalent to 
the desired result.  

If we raise the eigenvalue equation for ${\mathcal K}_p$ to the 
$\alpha$ power, we obtain. 
\begin{equation} 
\begin{split} 
\rho_p^\alpha \psi_p^\alpha & = 
( A^p \psi_p  \circ f_1 + \psi_p \circ f_2)^\alpha \\
& > A^{p \alpha} \Psi_p^\alpha + \psi_p^\alpha \circ f_2 \\
& = {\mathcal K}_{p \alpha} \psi_p^\alpha
\end{split} 
\end{equation}
The inequality above is a consequence of the 
binomial theorem for fractional powers.
Applying  Proposition~\ref{KRbound} we obtain the desired result
$\rho_{p \alpha} < \rho_p^\alpha$.

To prove the strict log-convexity of $\rho_p$, we 
argue similarly.
We will show that for $0< p < q$ and for 
$\rho_{(p +  q)/2} < \rho_p^{1/2} \rho_q^{1/2}$, which implies
strict convexity. 

Multiplying the eigenvalue equations 
satisfied by the spectral radius
for ${\mathcal K}_p$ and for  ${\mathcal K}_q$ 
and raising them to the $1/2$, we have:
\[
\begin{split}
(\rho_p \rho_q)^{1/2} &\psi_p^{1/2}\psi_q^{1/2}
 = 
\big[ ( A^p \psi_p  \circ f_1 + \psi_p \circ f_2) \cdot 
( A^q \psi_q  \circ f_1 + \psi_q \circ f_2) \big]^{1/2} \\
&  >
\left( A^{p +q} (\psi_p \psi_q) \circ f_1 
+  (\psi_p \psi_q) \circ f_1 \right)^{1/2} \\
&  >  A^{(p +q)/2} (\psi_p \psi_q)^{1/2} \circ f_1 
  + (\psi_p \psi_q)^{1/2} \circ f_1 )^{1/2} \\
&= \rho_{(p+q)/2}(\psi_p \psi_q)^{1/2}
\end{split}
\]
The first inequality comes by expanding the product and 
ignoring the cross terms, and the second inequality is
a consequence of $(1 +x)^{1/2} < 1 + x^{1/2}$ for $x > 0$.

Applying Proposition~\ref{KRbound}, we obtain the desired 
result.
\qed

\endproof
\subsection{Asymptotic properties of the renormalization}
From Lemma \eqref{Kumconvergence}, we know that
$\|T^nf-g_k\|_{D_k}=O(\omega^n_k)$
implies 
\beqnl{Kfsconvergence}
\|\mathcal K_{T^nf,p} - \mathcal K_{g_k,p}\|
 = O(\omega^n_k)
\eeqnl
for all $p>0$. Therefore,  the asymptotic properties of
 $\mathcal K_{T^nf,p}$ \, become similar to those of $\mathcal K_{g_k,p}$
for large $n$. In Proposition \ref{productalignment} below, we 
show that the dominant eigenspaces of the Lindeberg--Lyapunov operators, 
after a large number of renormalizations,  align to  
the dominant eigenspace of the cumulant operator at the Feigenbaum
fixed point $g_k$.
\begin{prop}\label{productalignment}
Let $\rho_p$  be  the spectral radius of 
$\mathcal K_p=\mathcal K_{g_k,p}$.  There exists a 
neighborhood  $\mathcal U$ of $g_k$ such that if
 $f\in\mathcal W_s(g_k)\cap \mathcal U$ and $h$ is a
positive function
 $h\in H^r$ then,
\beqnl{alignment}
c^{-1}\rho_p^n \leq \mathcal K_{T^{n-1}f,p}
\cdots\mathcal K_{f,p}h(z)
\leq c \rho_p^n
\eeqnl
for some constant $c$.
\end{prop}
The proof of this Proposition will be done after stating a
consequences.
A  Corollary  of \eqref{alignment} is the following
\begin{cor}\label{kurto2n}
For any $p\geq0$, let   $\rho_p$ be the spectral radius
 of $\mathcal K_p=\mathcal K_{g_k,p}$, then:
\begin{itemize}
\item[1.] There are constants $0<c<C$ such that
 \beqnl{lambdaasympat0}
c\leq \inf_{z\in I}\frac{\Lambda^f_p(\Gamma^f_n z,2^n)}
{\{|\lambda_k|^p\rho_p\}^n} \leq \sup_{z\in I}
\frac{\Lambda^f_p(\Gamma^f_n z,2^n)}
{\{|\lambda_k|^p\rho_p\}^n}\leq C
\eeqnl
\item[2.] If $f\in\mathcal W_s(g_k)$ is close to $g_k$ then,
 for any $p>2$
\beqn
\lim_{n\rightarrow\infty}\sup_{z\in I} 
\frac{\Lambda^f_p(\Gamma^f_n z,2^n)}
{\{\Lambda^f_2(\Gamma^f_n z, 2^n)\}^{p/2}}= 0
\eeqn
\end{itemize}
\end{cor}
\proof (1.) From \eqref{feig_smallscale} and Proposition
\ref{productalignment} we have that
\beqn
c \rho^n_p |\Gamma^f_n|^p\leq \Lambda_p(\Gamma^f_n z,2^n)\leq
C \rho^n_p |\Gamma^f_n|^p
\eeqn
for constants $0<c<C$. Part (1) follows from the 
exponential convergence since \eqref{rgconvergence}
implies that 
$$\hat{c} |\lambda_k|^n\leq |\Gamma^f_n|\leq \hat{C} |\lambda_k|^n$$
for some positive constants $\hat{c}$ and $\hat{C}$. 

\bigskip
\noindent
Parts (2) follows from part (1) and Proposition 
 \ref{Comparison} since for all $p>2$, 
we have $\rho_p<\sqrt{\rho^p_2}$.
\endproof
The proof of Proposition \ref{productalignment} is based on the 
following general result (Proposition \eqref{al})
on convergent sequences of  operators on a Banach space.
The method that we will use  is a standard technique in hyperbolic
theory \cite{hirshpugh}, where it is used to 
establish the stability of hyperbolic splittings

\begin{lem}\label{al}
Let $X$ be a Banach space, $\{ \mathcal K_n\}_{n =1}^\infty$
be a convergent sequence
 of
operators
 on $X$ such that 
$\lim_n \|\mathcal K-\mathcal K_n\|=0$.
Assume that there is a
 decomposition of $X=E^<\oplus
E^>$, and constants $0<\lambda_-<\lambda_+$ 
such that $E^{<,>}$ are closed subspaces and
\beqn
\mathcal K\left( E^> \right) &=& E^>\\
\mathcal K\left(E^<\right) &\subset& E^<\\
\spec\left(\mathcal K|_{E^<}\right)&\subset&
\{z\in\mathbb C: |z|< \lambda_-\}\\
\spec\left(\mathcal K^{-1}|_{E^>}\right)&\subset&
\left\{z\in\mathbb C:|z|<\frac{1}{\lambda_+}\right\}.
\eeqn
There exist $\epsilon>0$ such that if
$\sup_n\|\mathcal K_n-\mathcal K\| <\epsilon$, then
there are sequences of linear  operators $\{A_n\}_{n=0}^\infty
\subset \mathcal L(E^<,E^>)$
and  $\{B_n\}_{n=0}^\infty \subset \mathcal L(E^>,E^<)$
with $\lim_n\|A_n\|=0$  and $\lim_n\|B_n\|=0$,  such that
the spaces
\beqnl{defspaces}
E^>_n &=&\left\{x \oplus A_n x:\quad x\in E^<\right\}\\
E^<_n &=&\left\{x \oplus B_n x:\quad x\in E^<\right\}
\eeqnl
satisfy $X = E^>_n \oplus E^<_n$
and 
\beqnl{defdecomp1}
E^>_n =  \mathcal K_n\cdots\mathcal K_1 (E^>_0)
\eeqnl
\beqnl{defdecomp2}
E^<_n =  \mathcal K_n\cdots\mathcal K_1 (E^<_0)
\eeqnl
\end{lem}

\proof
 Since $X$ is the direct sum of invariant closed subspaces
 with respect to $\mathcal K$, using adapted norms 
\cite[Appen. A]{cabrefontichdelallave} we can assume without loss of
generality that
\begin{eqnarray}
\left\|\mathcal K|_{E^<}\right\|&\leq& \lambda_- \label{adaptednorm1}\\
\left\|\mathcal K^{-1}|_{E^>}\right\|&\leq& \frac{1}{\lambda_+}
\label{adaptednorm2}
\end{eqnarray}
and that if $x=x^<+ x^>$ with $x^<\in E^<$ and
$x^>\in E^>$,  then $\|x\|=\|x^<\| + \|x^>\|$.

The operators $\mathcal K$ and $\mathcal K_n$ are represented by
the following  matrices of operators with respect to the fixed 
decomposition $X = E^< \oplus E^>$
\begin{eqnarray}
\mathcal K &=& \begin{pmatrix}\mathcal K^< & 0 \cr
         0 & \mathcal K^>\cr
\end{pmatrix}\label{Koperator}\\
\mathcal K_n
 &=&\begin{pmatrix} \mathcal K^<+\delta^{<<}_n & \delta^{<>}_n \cr
                    \delta^{><}_n     &  \mathcal K^>+\delta^{>>}_n \cr
\end{pmatrix}\label{Knoperator}
\end{eqnarray}
where $\|\delta^{<<}_n\|+\|\delta^{<>}_n\|+\|\delta^{><}_n\|+
\|\delta^{>>}_n\|\rightarrow0$. 

\noindent
I. We show first the existence of the sequence of 
transformations $\{A_n\}$ in \eqref{defdecomp1}. 

If use \eqref{Koperator} to compute the image of a 
point $x + A_n x$ we obtain that it has its 
$E^>, E^<$ components are, respectively:
\beqnl{iteratedecomposition}
& & (K^> +  \delta^{>>}) A_n x  + \delta^<> x \\
& & (K^< + \delta^{<<}) x  + A_n \delta^>< x 
\eeqnl

The invariance condition that 
$K_n( x + A_n x)$ is in the graph of $A_{n+1}$
is equivalent to saying that  applying 
$A_{n+1}$ to the second row of \eqref{iteratedecomposition},
we obtain the first. 

Therefore, the invariance condition is equivalent to
\[
\begin{split}
& A_{n+1}[ (K^< + \delta^{<<}) + A_n \delta^>< ]  \\
& = (K^> +  \delta^{>>}) A_n   + \delta^<> 
\end{split}
\]

Isolating the first $A_n$ in the RHS, we obtain that 
the invariance  \eqref{defdecomp1} for the spaces  is equivalent to
the following invariance equation for $\{A_n\}$
\beqnl{inveq1}
A_n=\left(\mathcal K^>+\delta^{>>}_n\right)^{-1}\left\{A_{n+1}
\left(\mathcal K^< + \delta^{<<}_n\right) +
 A_{n+1}\delta^{<>}_nA_n -\delta^{><}_n\right\}
\eeqnl
Notice that  if $\{A_n\}\in\mathcal A$ satisfies \eqref{inveq1},
then \eqref{defdecomp1} holds.

We will show that there is a unique solution $\{A_n\}$ of 
 the invariance equation \eqref{inveq1} in
 an appropriate space 
of sequences. To be more precise, let
 ${\mathcal A}$ be the space of sequences of  operators in 
${\mathcal L}(E^<,E^>)$  that converge to zero. Observe 
that $\|A\|=\sup_n\|A_n\|$
defines a  norm on $\mathcal A$, which makes it into a 
Banach space.

Let $\tau$ be the operator on ${\mathcal A}$
 defined by the right--hand side of equation \eqref{inveq1}.
 For $A$ and $\widetilde{A}$ in ${\mathcal A}$ we have that
\begin{multline}
|\tau(A)_n - \tau(\widetilde{A})_n| \leq 
|(\mathcal K^>+\delta^{>>}_n)^{-1}(A_{n+1}-\widetilde{A}_{n+1})
(\mathcal K^<+\delta^{<<}_n)| +\nonumber\\
|(\mathcal K^>+\delta^{>>}_n)^{-1}||(A_{n+1}-
\widetilde{A}_{n+1})\delta^{<>}_nA_n|+ \nonumber\\
|(\mathcal K^>+\delta^{>>}_n)^{-1}||\widetilde{A}_{n+1}\delta^{<>}_n
(A_n-\widetilde{A}_n)|
\end{multline}
Let $R\geq\max(\|A\|,\|\widetilde{A}\|)$. 
Since $\delta_n\rightarrow0$, we can assume with no loss of 
generality that  for a given $\epsilon>0$,
$\sup_n\|\delta_n\|<\epsilon$.
\beqn
\|\tau(A)-\tau(\widetilde{A})\|\leq 
(\lambda^{-1}_+\lambda_- +\epsilon)\|A-\widetilde{A}\|+
2(\lambda_+^{-1}+\epsilon)\epsilon\|A - \widetilde{A}\|R
\eeqn
Since $\lambda_+^{-1}\lambda_-<1$, we can  choose
  $\epsilon>0$ small enough so that
\begin{equation*}
\begin{array}{lcl}
\epsilon+2R(\lambda^{-1}_+\epsilon)\epsilon &<&
 \frac{1-\lambda^{-1}_+\lambda_-}{2}\\
\noalign{\vskip6pt}
(\lambda^{-1}_+\lambda_-+\epsilon)R+
(\lambda^{-1}_+\epsilon)(\epsilon R^2+\epsilon)&<& R
\end{array}
\end{equation*}
we get that  $\tau$ is a contraction on 
the ball $\bar{B}(0;R)\subset{\mathcal A}$.
Therefore,  there exists a unique solution, $ A$,\  to the invariance 
equation \eqref{inveq1}.

\bigskip
\noindent
II. An  argument very similar to the one above shows
 that \eqref{defdecomp2} holds if and only 
if $\{B_n\}\subset\mathcal L(E^>,E^<)$ satisfies the
following invariance equation
\beqnl{inveq2}
B_{n+1}=\left\{(\mathcal K^<+\delta^{<<}_n)B_n +
\delta^{<<}_n-B_{n+1}\delta^{><}_n B_n \right\}
(\mathcal K^>+\delta^{>>}_n)^{-1}
\eeqnl
Let $\mathcal B$ be the space of sequences in 
$\mathcal L(E^>, E^<)$ that converge
 to zero with norm $\|B\|=\sup_n\|B_n\|$, and let $\tilde{\tau}$ be the
operator on $\mathcal B$ defined by the right hand side of \eqref{inveq2}.
An idea similar to the one  used in part (I),
shows that $\tilde{\tau}$ is a contraction on an appropriate closed ball.
Therefore, there is a unique solution to the invariance 
equation \eqref{inveq2}.
\endproof
\begin{remk}
$\|A_n\|$ and $\|B_n\|$ are measures of distance between the spaces
$E^{<,>}_n$ and $E^{<,>}$. Lemma \ref{al}  states
 that these respective spaces  align in the limit. 
\end{remk}
The next Corollary gives control on the growth of the products 
of operators $\mathcal K_n$.
\begin{cor}\label{cumopproduct}
Under the hypothesis of Proposition~\ref{al}, let us assume that dim$(E^>)=1$.
Then:
\begin{itemize}
\item[(1)]  For any $v\in X$
\beqnl{alignment2}
\prod_{j=1}^n(\lambda_+  -\epsilon_j)\,\|v\|\leq \|\mathcal K_n\cdots 
\mathcal K_0 v\|\leq
\prod_{j=1}^n(\lambda_+ + \epsilon_j)\, \|v\|
\eeqnl
where $\epsilon_n=\epsilon_n(\|\mathcal K_j -\mathcal K\|)$,
is such that $\lim_{t\rightarrow0}\epsilon_n(t)=0$.
\item[(2)] If   $\|\mathcal K_n-\mathcal K\|=O(\omega^n)$ 
with $0<\omega<1$ , then we have that 
\beqnl{alignmentexpo}
c^{-1}\lambda_+^n\|v\|\leq \|\mathcal K_n\cdots 
\mathcal K_0 v\|\leq
c \lambda_+^n\, \|v\|
\eeqnl
for some  constant $c>0$.
\end{itemize}
\end{cor} 
\proof
Using adapted norms, we can assume that \eqref{adaptednorm1}, 
\eqref{adaptednorm2} hold, and that
 if $x=x^<+ x^>$ with $x^<\in E^<$
and $ x^>\in E^>$, then  $\|x\|=\|x^<\|+\|x^>\|$.
For $s=$``$<$'', `` $>$'', let $E^s_0=E^s$ and 
define $E^s_{n+1}= \mathcal K_n(E^s_n)$.

\noindent
(1) Let $v \in E^<_n$, then by Lemma \ref{al} there exists
a unique $w\in E^<$ such that
$$v=w+A_nw$$
Since $\|A_n\|,\,\|B_n\|\rightarrow0$, we may further assume that 
\beqn
\sup_n\|A_n\|&\leq&\frac{\delta}{2\|\mathcal K\|}\\
\sup_n\|B_n\|&\leq&\frac{\delta}{2\|\mathcal K\|}\\
\|\mathcal K_n-\mathcal K\|&\leq&\frac{\delta}{2}
\eeqn
for $\delta>0$ small. 
Since $\|\mathcal K w\|< \lambda_-\|w\|\leq \lambda_+\|v\|$, 
from triangle inequality we have
\beqn
\|\mathcal K_n v\|&\leq& \|(\mathcal K_n- \mathcal K)v\| +
\|\mathcal K v\|\\
&\leq& \frac{\delta}{2}\|v\| + \|\mathcal K (w + A_nw)\|\\
&\leq& (\lambda_+ + \delta)\|v\|.
\eeqn

Let  $v\in E^>_n$ then, by Lemma \ref{al} there is unique
 $z\in E^>$ such that
$$ v = B_n z + z.$$
Therefore, $\|v\|\leq \frac{\delta}{2\lambda_+}\|z\| + \|z\|$, which
is equivalent to
$$\|z\|\geq \frac{\|v\|}{1+\delta}$$
From the triangle inequality, we have that
\beqn
\|\mathcal K_n v\|&\geq& \|\mathcal K v\|-
\|(\mathcal K_n - \mathcal K)v\|\\
&>& \|\mathcal K (B_n z + z)\| - \frac{\delta}{2}\|v\|\\
&>& \|\mathcal K z\| - \frac{\delta}{2}\|v\|\\
&>& (\lambda_+ -\delta) \|v\|
\eeqn
Combining these results, we obtain inequality \eqref{alignment2}.

(2) In the case of exponential convergence, \eqref{alignmentexpo} follows
from an elementary result on sums and products, namely: if
$\lim_n \sup_j{|c_{nj}|}=0$, and $\sup_n\sum_{j=1}^n|c_{n,j}|<\infty$, 
 then
$$\prod_j(1+c_{nj})\rightarrow e^{c}
\quad\text{if}\quad\sum_jc_{nj}\rightarrow c$$
See \cite[p. 78]{durr:prob}.
 \endproof

\bigskip
\noindent
{\em Proof of Proposition \ref{productalignment}.}
By Proposition \ref{KR} we have that
 that the spectral radius $\rho_p$ of $\mathcal K_p$ is a simple
 eigenvalue. Then, \eqref{alignment} follows from Corollary
\ref{cumopproduct}  and the exponential convergence, see
\eqref{Kfsconvergence}, of the Lindeberg--Lyapunov operators.
$\hfill{\Box}$
\begin{remk}
The estimates \eqref{lambdaasympat0} for $p=2$
are obtained in \cite{vul:feiguniv} using the 
Thermodynamic formalism.
\end{remk}
\subsection{Cumulant operators}\label{cumsec}
In Section \ref{secgencum} we introduced the concept of
 cumulants, see \eqref{cumulantes}.
One important fact is the cumulants of random variables, normalized
to have variance one, provide an  empirical measures of resemblance to
Gaussian. In this section, we make a connection between the spectral
properties of the Lindeberg--Lyapunov operators studied in
Section \ref{sec31} and the cumulants of the propagation
of noise.

Recall that the functions $f\in\mathcal W_s(g_k)$ are defined in 
a domain $D_k$, see R2 in section \ref{symmunimodalclt}.
\begin{defn}\label{cumoperatorfeig}
For each $f\in\mathcal W_s(g)$, the family of cumulant operators
$\{\tilde{\mathcal K}_{f,p}: p\in\mathbb N\}$  acting  on 
the space  $H^r(\bar{D}_k)$ are defined by
\beqnl{cumop1}
\tilde{\mathcal K}_{f,p}h(z)=\frac{1}{(\lambda_f)^p}\left[
\left(f'\circ f(\lambda_f z)\right)^p h(\lambda_f z)
+ h(f(\lambda_f z))\right]
\eeqnl 
\end{defn}

By an argument identical to that given in
 Lemma \ref{Kumconvergence},
 it follows from  the exponential convergence of 
$T^n f$  that 
\begin{lem}\label{cumlemma}
Let $f\in\mathcal W_s(g_k)$ close to $g_k$. For each nonnegative
integer $p$,  $\widetilde{\mathcal K}_{f,p}$ is a
compact operator on $H^r(\overline{D}_k)$. Furthermore, 
and $\widetilde{\mathcal K}_{T^{n-1}f,p}$ converges
to  $\widetilde{\mathcal K}_{g_k,p}$ exponentially fast.
\end{lem}

Notice  from Definition \ref{cumoperatorfeig} that
\begin{equation}
\begin{array}{lllr}
\widetilde{\mathcal K}_{f,p}&=&
\mathcal K_{f,p}& p\, \mathop{even}\\ 
\noalign{\vskip6pt}
\|\widetilde{\mathcal K}_{f,p}\|&\leq& 
\|\mathcal K_{f,p}\|&p\, \mathop{odd} 
\end{array}
\label{oprelation}
\end{equation}
Denote the spectral radius
 of $\widetilde{\mathcal K}_{f,p}$ and $\mathcal K_{f,p}$ 
by $\tilde{\rho}_{f,p}$ and $\rho_{f,p}$ respectively. Then,
from \eqref{oprelation} we have that 
$\tilde{\rho}_{f,p}=\rho_{f,p}$ for all $p$ even, 
and for all integers $p>2$
\beqnl{plargerthan2}
\tilde{\rho}_{f,p}\leq
\rho_{f,p}\leq \left(\rho_{f,2}\right)^{p/2}
\eeqnl

Let $f\in\mathcal W_s(g_k)$ be fixed, 
and recall from Section \ref{secgencum} 
the cumulant functional $\widetilde{\Lambda}$. 
For $x\in I$,  let
$$\tilde{k}_{n,p}(x)=
 \Gamma_n^{-p}\widetilde{\Lambda}_p(\Gamma_n x,2^n)$$
Then, from \eqref{wlambdarelation}, we have that
\beqnl{cumulantproduct}
\tilde{k}_{n,p}(s)(x)=\widetilde{\mathcal K}_{T^{n-1}f,p}\cdots
\widetilde{\mathcal K}_{f,p}\uno(x)
\eeqnl

We have that from Corollary \ref{cumopproduct}
that for very $h\in H^r(\overline{D}_k)$ there is a constant $c$
such that
\beqn
c^{-1}\tilde{\rho}^n_p\leq \|\widetilde{\mathcal K}_{T^{n-1}f,p}\cdots
\widetilde{\mathcal K}_{f,p}h\|\leq
c \tilde{\rho}_p^n
\eeqn
Recalling the notation $K_p[\xi]$ for 
the $p$--th order cumulant of a random
variable $\xi$, we have 
from \eqref{plargerthan2} that 
\beqnl{kurtosisfeig}
K_p\left[\frac{L_{2^n}(\Gamma_n^f\,x)}
{\sqrt{\var[L_{2^n}(\Gamma_n^f\,x)]}}\right]\asymp
\left(\frac{\tilde{\rho}_p}{{\tilde{\rho}_2}^{p/2}}\right)^n
\frac{h_{g_k,p}(x)}{(h_{g_k,2}(x))^{p/2}}
\rightarrow0
\eeqnl
as $n\rightarrow\infty$, where $h_{g_k,p}$ is
a positive eigenfunction of $\mathcal K_{g_k,p}$.
From the  compactness of the cumulant and Lindeberg--Lyapunov
operators, the asymptotic expansions \eqref{kurtosisfeig}
are  improved by adding more terms that involve eigenfunctions
of smaller eigenvalues \cite{fontichdelallave}. That is, 
we start from the standard asymptotic expansion of powers of 
a linear compact operator:
\beqnl{asympexpancum}
\mathcal K^n_{g_k,p} u(z)=\rho_p^n  h_{g_k,p} c(u) +
\sum_{j=1}^N \bar{\mu}_j^n c_j(u) \Psi_j (z)+ O(\mu^n_{N+1})
\eeqnl
where $\bar{\mu}_j$ are the Jordan blocks 
associated to the eigenvalues  $\mu_j$ (in decreasing order of 
size), $c_j(u)$ are the projections on this space and $\Psi_j(z)$
is a basis for the eigenspace. 

For the nonlinear renormalization operator, the theory of 
\cite{fontichdelallave} shows that there are very similar 
expansions for the asymptotic noise. The expansion involves 
not only powers of the eigenvalues of $\mathcal K$ but also powers of 
products of eigenvalues.  

These asymptotic approximations are very reminiscent 
to the 
higher order asymptotic expansions in the central limit theorem,
 such as Edgeworth expansions \cite[p. 535]{feller:prob}. 
(The asymptotic expansions \eqref{asympexpancum}
corresponds to Edgeworth expansions of $2^n$ terms). 
It seems possible that one could also prove Edgeworth expansions 
for the asymptotics of the scaling limits but it looks like 
the powers appearing in the expansions will be very different form 
the standard semi-integer powers and will be related to 
the spectrum of some of the operators giving the renormalization of 
the cumulants.

\subsubsection{Numerical conjectures on  spectral properties of
the cumulant operators and the Lindeberg--Lyapunov operators}
If $p$ is even, we know from Proposition \ref{KR} that
$\rho_{f,p}$ is a that single eigenvalue of 
$\widetilde{\mathop K}_{f,p}=\mathcal K_{f,p}$.

If $p$ is odd, we have  by compactness that 
 $\widetilde{\mathcal K}_{f,p}$ has an 
eigenvalue
 $\tilde{\mu}_{f,p}\in\{z\in\mathbb C: |z|=\tilde{\rho}_{f,p}\}$.

In the quadratic case ($k=1$), numerical computations
 \cite{DiazL} suggest that 
\begin{itemize}
\item[C1.] If  $\spec(\mathcal K_{g_1,p})\setminus\{0\}=\{\mu_{n,p}\}_n$ with
$|\mu_{n,p}|\leq |\mu_{n-1,p}|$, then $\mu_{n,p}\sim\lambda_1^n$
(Here, $u_n\sim v_n$ means that  $\lim_n u_n/v_n =1$).
\item[C2.] $\spec({\mathcal K}_{g_1,p})\setminus\{0\}$
 is asymptotically real,
that is: $\mu_n\in\mathbb R$ for all $n$ large enough. 
\end{itemize}
For all $p$ odd 
\begin{itemize}
\item[C3.] $\tilde{\rho}_{g_1,p}$ is an eigenvalue of
 $\widetilde{\mathcal K}_{g_1,p}$. The rest of 
$\spec(\widetilde{\mathcal K}_{g_1,p})$
is contained in the interior 
of the ball of radius $\tilde{\rho}_{g_1,p}$.
\item[C4.] $\lambda_1^{-2p}\sim{\rho}_{g_1,p}\sim \rho_{g_1,p}$ 
\, as  $p\rightarrow\infty$.
\item[C5.] $\spec(\widetilde{\mathcal K}_{g_1,p})\setminus\{0\}
=\{\tilde{\mu}_{n,p}\}$, where $|\tilde{\mu}_{n,p}|\leq |\nu_{n-1,p}|$ 
is asymptotically real.
\item[C6.] 
$\tilde{\mu}_{n,p}\sim\mu_{n,p}\sim \lambda_1^n$.
\end{itemize}
Conjectures C1--C6 are related to conjectures 
on  behavior of 
the Perron-Frobenius operators of real analytic expanding maps,
see  \cite{baladi,mayer,rugh}. The case $p=1$ corresponds to the
problem of reality of the spectrum of the linearized
period doubling operator. It is observed in  \cite{cvirugh} that
the spectrum in this special case appears to be  real, and that all the 
eigenvalues behave as $\lambda_1^n$.

\subsection{Proof of the central limit theorem
 for systems in a domain of universality (Theorem \ref{fclt})}
\label{sec:feig_clt}
\subsubsection{Introduction}
In this section, we prove the central limit theorem
(Theorem \ref{fclt}) for  maps $f$ in the 
domain of universality $\mathcal W_s(g_k)$. 

The first step is to prove a central limit theorem and Berry--Esseen
estimates for  $x_n(x,n)$ with initial condition $x$ in the 
orbit of $0$. 

We will verify the Lyapunov condition \eqref{lyapcon}  along
the whole sequence of positive integers. Then, we apply the result of 
Theorems \ref{main} and \ref{berry-esseen}.

The main observation is that the renormalization  theory developed
in Section~\ref{cumulantoper}
gives control on $\Lambda_s(x,2^n)$ for $x$ sufficiently small
(roughly $|x|\leq \left|\Gamma_{n+1}\right|$), see Corollary \ref{kurto2n}.
Then, the main tool to obtain control of the effect of the noise
on a segment of orbit, will be to decompose the segment into 
pieces whose lengths are powers of $2$ and chosen so that the 
renormalization theory applies. 

This decomposition will be accomplished in 
Section~\ref{binarydecomposition}.
The main conclusion of this section will be that this decomposition is 
possible and that, the final effect over the whole interval is 
the sum of the binary blocks -- which are
understandable by renormalization affected by weights. See
\eqref{decompositionoflambda}. The weights measure how the effect of 
the dyadic block propagates to the end of the interval.

The effect of the dyadic blocks in the final result depends on 
somewhat delicate arguments. This follows from the 
the following heuristic 
argument.

Note that when an orbit passes though zero, 
the derivatives become very small (hence the effect of the noise 
up to this time is erased). Hence, when an orbit includes extremely close 
passages through zero, the next iterations tend to behave as if
they were short orbits (and, hence the noise is far from Gaussian). 
On the other hand, the close passage to zero are unavoidable
(and they are the basis of the renormalization). 

Hence, our final result will be some quantitative 
estimates that measure the effect of the returns. The main 
ingredients is how close the returns are and how long do we
need to wait for the following close return (this is closely 
related to the gaps in the binary decomposition). 
The combination
of these properties depends on quantitative properties of 
the fixed point. See Theorem~\ref{kurtnx}. 

The fact that the gaps in the binary decomposition play a role 
is clearly illustrated in the numerical calculations in 
\cite{DiazL}. 
\subsubsection{The binary decomposition}
\label{binarydecomposition}
Given a number $n \in \mathbb N$, we will use the binary expansion
\beqnl{nbin}
n=2^{m_0}+\cdots+2^{m_{r_n}}
\eeqnl
where  $m_0> m_1>\ldots >m_{r_n}\geq0$.

If  necessary,   we will use $m_j(n)$ to emphasize the
dependence on  $n$ of the
power $m_j$ in \eqref{nbin}.

We will denote  the integer part of any number 
$u$ by  $[u]$.  Notice that
 that $m_0(n)= [\log_2(n)]$, and 
$r_n \leq m_0$.

Using \eqref{lambdarelation} several times we have 
for all $p\geq0$ and for any $x \in [-1,1]$ that 

\beqnl{decompositionoflambda}
\Lambda_p(x,n)=\sum^{r_n}_{j=0}\left|\Psi_{j,n}(x)\right|^p\Lambda_p\left(
f^{n-2^{m_j}-\cdots -2^{m_r}}(x), 2^{m_j}\right)
\eeqnl
\beqnl{psiderivatives}
\Psi_{j,n}(x)&=&
\left(f^{n-2^{m_0}-\cdots-2^{m_j}}\right)^\prime\circ
\left(f^{2^{m_j}+\cdots 2^{m_0}}\right)(x)
\eeqnl

When the  initial value $x_0=x$ of  $x_n=f(x_{n-1})$ 
satisfies  $|x|\leq \left|\Gamma_{m_0+1}\right|$, the points 
 \beqnl{return}
\upsilon_j=f^{2^{m_j}}(\upsilon_{j-1}),
\quad\upsilon_{-1}=x
\eeqnl
where $m_0>\ldots m_{r_n}\geq0$, are return points to a
small neighborhood of $0$ in the domain of renormalization.

The decomposition in dyadic scales \eqref{decompositionoflambda}
will be our main tool for the analysis of the growth of 
\eqref{decompositionoflambda}.

 We will show in Lemma \ref{growth0} that the size
 of the weights $\Psi_{j,n}(x)$ 
 are  controlled using the renormalization.
Then, using  Corollary \ref{kurto2n} we will
get control over the size of
the $\Lambda_p(\upsilon_{j-1},2^{m_j})$. Combining  these
results, the Lyapunov condition \ref{lyapcon} will follow for the orbit
of zero, thus establishing the central limit theorem.
\subsubsection{Estimating the weights}
The main tasks of this section will be to estimate the 
weights $\Psi_{n,j}$ in \eqref{decompositionoflambda}
and to show that they do not affect too much 
the conclusions of renormalization.

For $f\in\mathcal W_s(g_k)$ fixed, recall the notation $\lambda_f=f(1)$,
and for each nonnegative integer $m$, $\Gamma_m= f^{2^m}(0)$. Then  
$$\lambda_{T^nf}=\frac{\Gamma_{n+1}}{\Gamma_n}$$
If necessary, we will use $\Gamma^f_n$ to emphasize the 
role of $f$.  For the Feigenbaum fixed point $g_k$ of
order $2k$,  we use  $\lambda_k=\lambda_{g_k}$.  We have that
\beqn
\Gamma^{g_k}_n= \lambda_k^n
\eeqn

Recall   that the map $g_k$ has a series expansion of the form
$$g_k(x)=1+\sum_{l=1}^\infty b_l x^{2kl}$$
and that   $g'_k(1)=\lambda_k^{1-2k}$.
\begin{lem}\label{growthgk}
Let $h_k$ be the function  over the interval $[-1,1]$  defined by
\beqnl{gderivoveroderminusone}
h_k(x)=\frac{g'_k(x)}{x^{2k-1}}
\eeqnl
for $x\neq0$ and  $h_k(0)=2k\,b_1$. Then $h$ is a negative 
even function, increasing on $[0,1]$.
\end{lem}
\proof  From the fixed point equation
\beqn
g_k(x)=\frac{1}{\lambda_k}g_k(g_k(\lambda_k x))
\eeqn
we have that
\beqnl{formula}
h_k(x)=\frac{g'_k(g_k(\lambda_k x))}{g'_k(1)}h(\lambda_k x)
\eeqnl
Repeated applications on \eqref{formula} gives
\beqnl{formulaproduct}
h_k(x)=
h_k(\lambda_k^{n+1}x)\prod_{j=1}^n
 \frac{g'_k(g_k(\lambda_k^jx))}{g'_k(1)}
\eeqnl
Notice that $h_k(\lambda_k^{n+1}x)\rightarrow h_k(0)$ uniformly on 
$[-1,1]$ as $n\rightarrow\infty$.
Furthermore,  the factors  in the product \eqref{formulaproduct}
converge exponentially to $1$. Therefore, we can pass to the limit
as $n\rightarrow\infty$ in \eqref{formulaproduct}
and  get
$$h_k(x)=
h_k(0)\prod_{n=1}^\infty \frac{g'_k(g_k(\lambda_k^nx))}{g'_k(1)}$$
Each term in the product is a nonnegative decreasing function on $[0,1]$
therefore, $h_k$ is increasing on $[0,1]$.
\endproof

In the rest of this section, we assume that the 
 domain of universality  $\mathcal W_s(g_k)$  is fixed.
From the  exponential convergence  of $T^mf$ to $g_k$, 
we can choose a small neighborhood $\mathcal U$ of
$g_k$ so that of 
\beqnl{closeTmftog}
T(\mathcal W_s(g_k)\cap\mathcal U)\subset 
\mathcal W_s(g_k)\cap\mathcal U
\eeqnl
The following result gives control of derivatives for points
 close to the origin.
Some of the bounds for the derivatives in Lemma \ref{growth0}
are done in \cite{vul:feiguniv} using the thermodynamic formalism.

\begin{lem}\label{growth0}
Assume that $f\in\mathcal W_s(g_k)\cap \mathcal U$ and as before,
let $h_k(x)=g'_k(x)x^{1-2k}$. For each finite decreasing 
sequence of nonnegative
 integers  $\{m_j\}^r_{j=0}$, and $|x|\leq|\Gamma_{m_0+1}|$, let 
 $\{\upsilon_j\}^r_{j=-1}$ with $\upsilon_{-1}=x$ be as 
 in \ref{return}. There exists $0\leq\epsilon\ll 1$,
 depending only on $\mathcal U$, such that if 
\begin{eqnarray}
G&=&g_k(\lambda_k)-\epsilon\label{constant1}\\
c&=&|h(\lambda_k)|-\epsilon\label{constant2}\\
d&=&|h_k(0)|+\epsilon\label{constant3}
\end{eqnarray}
 then $1<c\,G^{2k-1}<d$, and
\begin{itemize}
\item[(1)] For each  $0\leq j\leq r$,
\beqnl{rel1}
G|\Gamma_{m_j}|\leq|\upsilon_j|\leq|\Gamma_{m_j}|
\eeqnl
\item[(2)] For $1\leq j\leq r$ 
\beqnl{rel2}
c\,G^{2k-1}\left|\frac{\Gamma_{m_{j-1}}}{\Gamma_{m_j}}\right|^{2k-1}
\leq
\left|\left(f^{2^{m_j}}\right)^\prime\left(\upsilon_{j-1}\right)\right|
\leq d \left|\frac{\Gamma_{m_{j-1}}}{\Gamma_{m_j}}\right|^{2k-1}
\eeqnl
\end{itemize}
\end{lem}

\bigskip
\noindent
\proof
By \eqref{closeTmftog}, we can assume without loss of
generality that 
\begin{eqnarray}
\sup_{m\in\mathbb N}\left|\lambda_{T^mf}- \lambda_k\right|\leq \epsilon
\label{closenessoflambda}\\
\left|T^mf(\lambda_k)-g_k(\lambda_k)\right|\leq \epsilon
\label{closenessofTfatlambda}
\end{eqnarray}
for some $0\leq \epsilon\ll 1$.
Since $g_k(\lambda_k)> |\lambda_k|$ and 
$$|g_k'(\lambda_k)|=
\left|\frac{g'_k(1)}{g'_k(g_k(\lambda_k))}\right|>1,
$$
by  choosing  $0\leq\epsilon\ll 1$ smaller if necessary, we get 
 $1<c\,G^{2k-1}< d$.

The exponential convergence 
implies that there is $0\leq\delta\ll 1$ such that
\beqnl{ratioofgamma}
\sup_{m\in\mathbb N}
\left|\Gamma_{m+1}\Gamma^{-1}_{m}\right|\leq \left|\lambda_{T^mf}\right|
<|\lambda_k|+\delta
\eeqnl
We know that or all $m\in\mathbb N$,   $T^mf$ is  even,
 decreasing in $[0,1]$,  and 
$\sup_{v\in I}|T^mf(v)|\leq1$. Shrinking 
$\mathcal U$ if necessary,   from \eqref{closenessofTfatlambda},
we can assume without loss of generality
that 
\beqnl{generaltmgcomparison}
\{g_k(\lambda_k)-\epsilon\}
\leq |T^m f (x)|\leq 1
\eeqnl
for all $x\in[\lambda_k-\delta,|\lambda_k|+\delta]$.

\noindent
 (1) For $j=0$ we have that
\beqn
\upsilon_0=f^{2^{m_0}}(\upsilon_{-1})=
\Gamma_{m_0} T^{m_0}f(\Gamma^{-1}_{m_0}\upsilon_{-1})
\eeqn
It follows from \eqref{closenessoflambda} that
 $$|\Gamma^{-1}_{m_0}\upsilon_{-1}|<|\lambda_k|+\delta$$
Then, by  \eqref{generaltmgcomparison} we get
\beqn
\{g_k(\lambda_k)-\epsilon\}|\Gamma_{m_0}|
\leq |\upsilon_0|\leq |\Gamma_{m_0}|
\eeqn

 We continue  by induction on $j$.  
Assume \eqref{rel1} holds for all numbers
$0\leq i\leq j$, where $j<r$.   For $j+1$ we have that
$$\upsilon_{j+1}=f^{2^{m_{j+1}}}(\upsilon_j)=\Gamma_{m_{j+1}}
T^{m_{j+1}}f\left(\Gamma^{-1}_{m_{j+1}}\upsilon_j\right)
$$
From
 \eqref{closenessoflambda} and \eqref{ratioofgamma}, 
it follows that
\beqn
\{g_k(\lambda_k)-\epsilon\}|\Gamma_{m_{j+1}}|
\leq |\upsilon_{j+1}|\leq |\Gamma_{m_{j+1}}|
\eeqn

\bigskip
(2) Shrinking $\mathcal U$ if necessary, we can  assume that
\beqnl{fclosetogk}
\sup_{x\in I}\left|\frac{f^\prime(x)}{x^{2k-1}}
-h_k(x)\right|<\epsilon
\eeqnl
for all $f\in\mathcal W_s(g_k)\cap\mathcal U$.

We know from  Lemma \ref{growthgk} that 
$$|h_k(\lambda_k)||x|^{2k-1}
\leq |g'_k(x)| \leq |h_k(0)|$$
 for all $|x|\leq |\lambda_k|$. Hence, for each $m\in\mathbb N$, 
 \eqref{fclosetogk} and the exponential convergence of
$T^mf$  imply that 
\beqnl{dtmfdgk}
 \left(|h(\lambda_k)|-\epsilon\right)|x|^{2k-1}
 \leq \left|\left(T^mf\right)^\prime(x)\right|
\leq \left(|h_k(0)|+\epsilon\right)|x|^{2k-1}
\eeqnl
for  $x\in[f(\lambda_f)|\Gamma_m|,\left|\Gamma_m\right|]$.

Therefore, from part (1) of the present Lemma, 
\eqref{dtmfdgk},  and  the identity
$$\left(f^{2^{m_{j+1}}}\right)^\prime(\upsilon_j)=
\left(T^{m_{j+1}}f\right)'\left(\Gamma^{-1}_{m_{j+1}}\upsilon_j\right)$$
we have that
\beqn
\left|\left(f^{2^{m_{j+1}}}\right)^\prime(\upsilon_j)\right|\geq&
c|\Gamma^{-1}_{m_{j+1}} \upsilon_j|^{2k-1}\geq& 
c\,G^{2k-1} \left|\frac{\Gamma_{m_j}}{\Gamma_{m_{j+1}}}\right|^{2k-1}\\
\left|\left(f^{2^{m_{j+1}}}\right)^\prime(\upsilon_j)\right|\leq&
 d |\Gamma^{-1}_{m_{j+1}}\upsilon_j|^{2k-1}
\leq& d \left|\frac{\Gamma_{m_j}}{\Gamma_{m_{j+1}}}\right|^{2k-1}
\eeqn
for all $0\leq j\leq r-1$. 
\endproof
\subsubsection{Lyapunov condition for all times for functions
 in the domain of universality}
In this section, we show now that the Lyapunov condition \ref{lyapcon}
 holds along the whole
sequence of positive integers for
orbits starting at any point of the form $x=f^l(0)$, $l\in\mathbb N$, with
 $f\in\mathcal W_s(g_k)\cap \mathcal U$.

\begin{prop}\label{kurton0}
Let $f\in\mathcal W_s(g_k)\cap\mathcal U$ be fixed, and for each
$n\in\mathbb N$, let  $m_0(n)=[ \log_2(n)]$, $r_n+1$ 
be the number of ones in the binary expansion of $n$, 
and $I_{m_0(n)}$ be the interval
 $\left[-|\Gamma_{m_0(n)+1}|, |\Gamma_{m_0(n)+1}|\right]$. 
\begin{itemize}
\item[(1)]  Let $c$, $G$, and $d$ be as in Lemma \ref{growth0}.
For any real number $p>0$, there are constants  $a_p$, $b_p$
 such that
\beqnl{growthofLambdap}
a_p\left(\frac{(c\, G^{2k-1})^{r_n}}{\lambda_k^{(2k-1)m_{r_n}}}\right)^p 
\leq 
\frac{\Lambda_p(x ,n)}
{(\lambda_k^{2kp}\rho_p)^{m_0(n)}}
 \leq  b_p\left(\frac{d^{r_n}}{\lambda_k^{(2k-1)m_{r_n}}}\right)^p
\eeqnl
for all $x\in I_{m_0(n)}$.
\item[(2)] For any real number $p>2$, we have that
\beqnl{ratio0}
\lim\limits_{n\rightarrow\infty}\sup_{x\in I_{m_0(n)}}
\frac{\Lambda_p(x,n)}
{\left\{\Lambda_2(x,n)\right\}^{p/2}}=0
\eeqnl
\end{itemize}
\end{prop}
\proof
For each $m\in\mathbb N$, let 
$I_m=\left[-|\Gamma_{m+1}|,|\Gamma_{m+1}|\right]$.
Let $n$ be fixed for the moment, and consider the sequence of
powers,  $\{m_j\}_{j=0}^{r_n}$, in the binary expansion
of $n$, see \eqref{nbin}.

Let  $x\in I_{m_0}$ be fixed, and 
consider the sequence of returns $\{\upsilon_j\}_{j=-1}^{r_n}$
with  $\upsilon_{-1}=x$ as in \eqref{return}.

From \eqref{decompositionoflambda}, we have 
for  each $p\geq0$ that
\beqnl{decomplambdanearzero}
\Lambda_p(x,n)= 
\sum^{r_n}_{j=0}|\Psi_{j,n}(x)|^p\Lambda_p(\upsilon_{j-1},2^{m_j}),
\eeqnl

\bigskip
\noindent
(1) {From}  Lemma \ref{growth0}, we have that
$$\left(c\, G^{2k-1}\right)^{r_n-j+1}
\leq
|\Psi_{j,n}(x)|\left|\frac{\Gamma_{m_{r_n}}}{\Gamma_{m_j}}\right|^{2k-1}
\leq d^{r_n-j+1}
$$
where   $1<c\,G^{2k-1}<d$. Since
$$|\Gamma_{m+1}\Gamma^{-1}_m-\lambda_k|\leq C\omega^m_k$$
with $0<\omega_k<1$, there
is a constant $a$ such that
\beqnl{weightfactor}
a^{-1}(c \, G^{2k-1})^{r_n-j+1}
\leq|\Psi_{j,n}(x)||\lambda_k|^{(m_{r_n}-m_j)(2k-1)}
\leq a d^{r_n-j+1}
\eeqnl
for all $x\in I_{m_0(n)}$.
By Corollary \ref{kurto2n}, there are 
positive constants $\hat{c}_p$, $\hat{C}_p$  such that
\beqnl{blocfactor}
\hat{c}_p\leq \Lambda_p(\nu_{j-1},2^{m_j})(|\lambda_k|^p\rho_p)^{-m_j}
\leq \hat{C}_p
\eeqnl
Using \eqref{weightfactor} and \eqref{blocfactor}, we obtain
estimates for the size of each term in \eqref{decomplambdanearzero}.
Since $\lambda_k^{2kp}\rho_p>1$, the estimate
 \eqref{decomplambdanearzero} follows by
noticing that  $m_0-m_j\geq j$.

\bigskip
\noindent
(2)  Since $0\leq r_n\leq m_0=[\log_2(n)]$, 
\eqref{growthofLambdap} implies that
there are  constants $0<C_p<D_p$ such that for all $n$ large enough
\beqnl{lambda_sgrowth}
C_p\{\lambda_k^{2k p}\rho_p\}^{m_0} \leq
\Lambda_p(x,n)\leq D_p \{d^p\rho_p\}^{m_0}
\eeqnl
From  Proposition \ref{Comparison}, we know that
 $\lambda_k^{2kp}\rho_p>1$ for all $p>0$, therefore
\beqnl{bigL2}
\lim_{n\rightarrow\infty}\inf_{x\in I_{m_0}(n)}\Lambda_p(x,n)=\infty
\eeqnl

Assume that $p>2$. By Corollary \ref{kurto2n}, we have that
for any given $\epsilon>0$, 
there is  an integer $M_\epsilon$ such that if $m\geq M_\epsilon$, 
then
 $$ \Lambda_p(x, 2^m)<\epsilon 
(\Lambda_2(x, 2^m))^{p/2}$$
for all $|x|\leq |\Gamma_m|$.

Denote by   $\mathcal S_{r_n}=\{0,\ldots,r_n\}$  and let 
$$
\mathcal A^1_n=\left\{j\in \mathcal S_{r_n}:
\sup_{u\in I_{m_j}}\frac{ \Lambda_p\left(u, 2^{m_j}\right)}
{(\Lambda_2\left(u, 2^{m_j}\right))^{p/2}}
\leq\epsilon\right\}$$
and $\mathcal A^2_n=\mathcal S_{r_n}\setminus\mathcal A^1_n$.
For each $x\in I_{m_0}$ and $s=2,p$, we split $\Lambda_s(x,n)$ as
$$\Lambda_s(x,n)=H_s(x,\mathcal A^1_n)+
H_s(x, \mathcal A^2_n),$$
where  $H_s(x,\mathcal A^i_n)$, $i=1,2$ is defined by
$$H_s(x,\mathcal A^i_n)=
\sum_{j\in \mathcal A^i_n}
|\Psi_{j,n}(x)|^s\Lambda_s(\upsilon_{j-1},2^{m_j})$$

Notice that 
 $H_p(x, \mathcal A^2_n)$ 
is the sum of bounded number of terms (at most $M_\epsilon$)
all of which are bounded. Hence,
\beqnl{firstpiece}
H_p(x,\mathcal A^2_n)\leq C_\epsilon
\eeqnl
for some constant $C_\epsilon>0$.

In addition,  notice that
\beqnl{secondpiece}
H_p(x, \mathcal A^1_n)\leq \epsilon (H_2(x,\mathcal A^1_n))^{p/2}
\leq \epsilon (\Lambda_2(x,n))^{p/2}
\eeqnl

Combining  \eqref{firstpiece}, \eqref{secondpiece},  we get
\beqn
\frac{\Lambda_p(x,n)}{(\Lambda_2(x,n))^{p/2}}&\leq&
\frac{H_p(x, \mathcal A^1_n)+H_p(x, \mathcal A^2_n)}
{(\Lambda_2(x,n))^{p/2}}\\
&\leq&\epsilon + \frac{C_\epsilon}{(\Lambda_2(x,n))^{p/2}}
\eeqn
Passing to the limit as $n\rightarrow\infty$,  we obtain from
 \eqref{bigL2}, that
$$\limsup\limits_{n\rightarrow\infty}
\sup_{x\in I_{m_0(n)}}\frac{\Lambda_p(x,n)}
{\{\Lambda_2(x,n)\}^{p/2}}
 \leq \epsilon.$$
By letting  $\epsilon\rightarrow0$, we get \eqref{ratio0}.
\endproof

For any  $l\in\mathbb N$  fixed and
for all $m>[\log_2(l)]$,  denote by
$I_{m}=\left[-|\Gamma_{m+1}|,|\Gamma_{m+1}|\right]$.
We  will show  that the Lyapunov condition \ref{lyapcon} holds
for orbits starting at any  point on $\{f^l(0)\}_{l\in\mathbb N}$.
\begin{thm}\label{kurtnx}
Let $f\in\mathcal W_s(g_k)\cap \mathcal U$ be fixed, 
and  $l\in\mathbb N$. For all $p>2$, 
\beqnl{ratio}
\lim\limits_{n\rightarrow\infty}\sup_{x\in I_{m_0(n+l)}}
\frac{\Lambda_p(f^l(x),n)}{\left(\Lambda_2(f^l(x),n)\right)^{p/2}}=0
\eeqnl
\end{thm}

The limit \eqref{ratio} is not uniform in $l$.

\proof
The idea of the proof is to use the dyadic decomposition 
\eqref{decompositionoflambda} and  show that for any fixed
integer $l\geq1$, the Lyapunov functionals
$\Lambda_p(f^l(x),n)$ and $\Lambda_p(x,n+l)$ for 
$x\in I_{m_0(n+l)}$ are very similar for large values of
$n$ (see \eqref{blockasymp}).

Let $n\in\mathbb N$ be a large integer and decompose $n+l$ in 
its binary decomposition:
$$n+l= 2^{m_0}+\cdots 2^{m_r}$$
where $m_0=m_0(n+l)=[\log_2(n+l)]$. Consider  $x\in I_{m_0(n+l)+1}$,
and observe that
$$\Lambda_p(f^l(x),n)=\sum_{j=l+1}^{n+l}
\left|\left(f^{n+l-j}\right)\circ f^j(x)\right|^p$$
Using $\{m_j\}$, consider the sequence of returns $\{\upsilon\}_j$
with $\upsilon_{-1}=x$ (see \eqref{return}). It follows that 
\beqnl{almostdyadic}
\begin{split}
\Lambda_p(f^l(x),n)= & |\Psi_{n+l,0}(x)|^p A_p(x,2^{m_0})
+ \sum_{i=1}^{r}|\Psi_{n+l,i}(x)\|^p\Lambda_p(\upsilon_{i-1},2^{m_i})
\end{split}
\eeqnl
where 
\beqnl{difference}
\Lambda_p(x,2^{m_0})-A_p(x,2^{m_0})=
\sum_{j=1}^{l}\left|\left(f^{2^{m_0}-j}\right)\circ f^j(x)\right|^p
\eeqnl
We will show that the difference \eqref{difference} is small
compared to $\Lambda_p(x,2^{m_0})$. Indeed, since
\beqn
\left(f^{2^{m_0}-j}\right)^\prime\circ f^j(x)=
\frac{\left(f^{2^{m_0}}\right)^\prime(x)}{\left(f^j\right)^\prime(x)},
\eeqn
it follows from  Lemma \ref{growth0} that
\beqnl{goodbound}
\frac{c_j}{|\lambda_k|^{(2k-1)m_0}}
\leq \left|\left(f^{2^{m_0}-j}\right)^\prime\circ f^j(x)\right|
\leq \frac{d_j}{|\lambda_k|^{(2k-1)m_0}}
\eeqnl
where  $c_j$ and $d_j$ depend only  on $j$. 
Since\eqref{difference} contains $l$ terms, 
 \eqref{radspecbound} and Corollary \ref{kurto2n} imply that
\beqnl{goodbound2}
\frac{\Lambda_p(x,2^{m_0})-A_p(x,2^{m_0})}
{\Lambda_p(x,2^{m_0})}\leq C_{l,p}\frac{1}{(\lambda_k^{2kp}\rho_p)^{m_0}}
\rightarrow0
\eeqnl
as $m_0\rightarrow\infty$. Therefore
\beqnl{blockasymp}
\lim_{n\rightarrow}\sup_{x\in I_{m_0(n+l)}}
\frac{\Lambda_p(f^l(x),n)}{\Lambda_p(x,n+l)}=1
\eeqnl

The limit \eqref{ratio}  follows from Proposition \ref{kurton0} 
and \eqref{blockasymp}.
\endproof

We will use the result of Theorem \ref{kurtnx} together with the results
of Theorems \ref{main} and \ref{berry-esseen} to prove the central 
limit theorem
for systems in the domain of universality $\mathcal W_s(g_k)$.

\subsubsection{Proof of Theorem \ref{fclt}}
By Theorem \ref{main} and Theorem \ref{kurtnx}, it suffices to
choose $\sigma_n$ that satisfies \eqref{weaknoise}.

Recall that for any integer $n$, $m_0(n)=[\log_2(n)]$ and 
$r_n+1$ equals the number of ones in the binary expansion
of $n$.

For  $l\in\mathbb N$ fixed, denote by $y_l=f^l(0)$. 
Let $c$, $G$, and $d$ be as in Lemma \ref{growth0}. 
From \eqref{blockasymp} and Proposition
\ref{kurton0}, for each
$p>0$ there  there are constants $C_{p,l}$ and $D_{p,l}$
such that for all for all $x\in I_{m_0(n+l)}$,
\beqn
C_{p,l}\left(\frac{(c\, G^{2k-1})^{r_{n+l}}}
{\lambda_k^{(2k-1)m_{r_{n+l}}}}\right)^p 
\leq 
\frac{\Lambda_p(f^l(x),n)}
{(\lambda_k^{2kp}\rho_p)^{m_0(n+l)}}
 \leq  D_{p,l}\left(\frac{d^{r_{n+l}}}
{\lambda_k^{(2k-1)m_{r_{n+l}}}}\right)^p
\eeqn

 Therefore, we have that
\beqnl{sigoptfeig}
1\leq \frac{(\widehat{\Lambda}(f^l(x),n))^3}{\sqrt{\Lambda_2(f^l(x),n)}}
\leq C\,\left(\frac{d^3\lambda_k^{4k}}{c\, G^{2k-1}}\right)^{r_{n+l}}
\left(\frac{\rho_1^3}{\sqrt{\rho_2}}\right)^{m_0(n+l)}
\eeqnl
for some constant $C>0$.

Hence, if  
\beqnl{caseobs}
\left|\frac{d^3\lambda_k^{4k}}{c\, G^{2k-1}}\right|\leq 1
\eeqnl
then,  it suffices to consider 
$$\sigma_m=\frac{1}{n^{\gamma+1}}$$
with $\gamma> \log_2(\rho_1^3/\sqrt{\rho_2})$.

In the case 
\beqnl{casevac}
\left|\frac{d^3\lambda_k^{4k}}{c\, G^{2k-1}}\right|> 1,
\eeqnl
it is enough to consider 
$$\sigma_n=\frac{1}{n^{\gamma^*+1}}$$
with 
$\gamma^*> \log_2(d^3\lambda_k^{4k}\rho_1^3)+
\log_2(c\,G^{2k-1}/\sqrt{\rho_2})$.

To obtain the Berry--Esseen estimates of 
Theorem \ref{berry-esseen}, we choose 
$\sigma_n$ that satisfies \eqref{weakernoise}

If \eqref{caseobs} holds, then it 
suffices to consider
$$\sigma_n=\frac{1}{n^\theta}$$
with
$\theta> \log_2(\rho_1^3/\sqrt{\rho_2})+\log_2(\rho_2^3/\rho_3)+1$.

In the case \eqref{casevac}, it suffices to consider
$$\sigma_n=\frac{1}{n^{\theta^*}}$$
with 
$\theta^*>  \log_2(d^3\lambda_k^{4k}\rho_1^3)+
\log_2(c\,G^{2k-1}/\sqrt{\rho_2})
+ \log_2(\rho_2^3/\rho_3)$ +1
$\hfill{\Box}$

\begin{remk}
If we consider an increasing sequence of iterations $\{n_k\}_k$ which
have  lacunar binary expansions, that is 
$$ \frac{r_{n_k}}{\log_2(n_k)}\rightarrow0$$
as $k\rightarrow\infty$ then,   
from \eqref{sigoptfeig}, it suffices to consider
$\sigma_{n_k}=n_k^{-(\gamma+1)}$
to have a  central limit holds along $\{n_k\}$. 
\end{remk}

\begin{remk}
In the quadratic ($k=1$) case, the power
  which arises from the argument
in the central limit theorem  presented here
 computed in \cite{DiazL}. The result is 
$\gamma=3.8836\ldots$
The method used in \cite{DiazL} follows  a numerical scheme similar
to \cite{lan:casproof} to compute $g_1$ and the
spectral  radii of the operators $\mathcal K_{g_1,1}$
and $\mathcal K_{g_1,2}$. 

The numerics in \cite{DiazL} give a confirmation that
the dependence of the optimal  $\sigma_n$ on the number 
of iterations is a power of exponent $\gamma_*$ which is
not too far from the  the $\gamma$ coming from the arguments in this 
paper.
\end{remk}

\subsubsection{Conjectures for orbits of points in the basin
of attraction}
The topological theory of interval maps
\cite{coeck:itemap, coecklan:univ} implies that for
 each $f\in\mathcal V\cap\mathcal W_s(g_k)$, the set
 set $\mathcal C_f=\overline{\{f^n(0)\}}_{n\in\mathbb N}$
 attracts  all points in $[-1,1]$,
 except for periodic points of period $2^n$ -- which are unstable -
and  their preimages.
More formally,  for any point $x\in I$, either
\begin{itemize}
\item[(a)] $\lim_{n\rightarrow\infty}d(f^{n}(x),\mathcal C_f)=0$, or
\item[(b)]$f^{2^n}(f^l(x))=f^l(x)$  for some
 nonnegative integers $l$ and $n$,
\end{itemize}
 In the second alternative, we have 
$\left|(f^{2^n})^\prime(f^l(x))\right|>1$.

 The set $\mathcal C_f$ is commonly referred to as the Feigenbaum attractor
since it attracts a set of full measure (not open, however). 
 It is also known that $f|_{\mathcal C_f}$ is uniquely
ergodic since it is a substitution system \cite{coecklan:univ}.
The set of points $\mathcal B_f$  that are attracted to $\mathcal C_f$ is
called the basin of attraction of $f$.

Numerical simulations \cite{DiazL} for the quadratic Feigenbaum
fixed point and the quartic map suggest that the effective noise
of orbits starting
 in the basin 
of attraction $\mathcal B_f$ and affected by weak noise
approaches a Gaussian.
We conjecture that
the is indeed the case.

 Of course, the  convergence 
to Gaussian
is not expected to be  uniform in the point in the basin. 
We have already shown that the 
convergence to Gaussian is false for the preimages of 
unstable periodic orbits 
which are accumulation points for the basin.

\section{Central limit theorem
 for critical circle maps}\label{critcircmap}

In this section, we consider another example of dynamical
systems with a non trivial renormalization theory,
 namely circle maps with a critical
point and golden mean rotation number. The theory has been 
developed both
heuristically and rigorously in 
\cite{feigkadashenk,ostrandseth,mestel,Lanf,sinaikhanTheRenorm,deFaAsymp} 

We will adapt   the argument developed in section 
\ref{symmunimodalclt} to  the case of critical circle maps with
golden mean rotation number.
\subsection{Critical circle maps
with golden mean rotation number}
We will consider  the following class of maps of the circle.
\begin{defn}\label{critcircmapdefn}
The space of critical circle maps is defined as the set of
analytic functions $f$,  that are strictly increasing in 
 in $\mathbb R$ and satisfy
\begin{itemize}
\item[G1.] $f(x+1)=f(x)+1$
\item[G2.] $f$ has rotation number
 $\beta=\frac{\sqrt{5}-1}{2}$.
Recall that for a circle map $f$ the rotation number $r(f)$
is defined by
$$r(f)=\lim_{n\rightarrow\infty}\frac{f^n(x)-x}{n}$$
\item[G3.]  $f^{(j)}(0)=0$ for all $0\leq j\leq 2k$,
and $f^{(2k+1)}(0)\neq0$.
\end{itemize}
\end{defn}
From the well known relation between the golden mean and the Fibonacci 
numbers, $Q_{0}=0$, $Q_1=1$ and $Q_{n+1}=Q_{n-1}+ Q_n$, it follows that
$Q_n \beta - Q_{n-1}=(-1)^{n-1}\beta^n$ is the rotation number of the map
\beqnl{newrotnumbermap}
f_{(n)}(x)= f^{Q_n}(x)-Q_{n-1}
\eeqnl

The sequence of maps $f_{(n)}$ will be useful in section \ref{cumopcirc}
where we analyze the cumulants of the noise.
\subsection{Renormalization theory of circle maps}
There are different rigorous renormalization formalisms  for 
circle maps. For our purposes,  we  will need very little about the 
renormalization group, so that we will  use the very basic 
 formalism of scaling limits.
\begin{defn}
Given a map $f$ as in Definition \ref{critcircmap},  
let $f_{(n)}$ be as in \eqref{newrotnumbermap} and
denote  $\lambda_{(n)}=f_{(n)}(0)$.
The $n$--th renormalization of $f$ is defined by
\beqnl{renormcirc}
f_n(x)=\frac{1}{\lambda_{(n-1)}}f_{(n)}(\lambda_{(n-1)} x)
\eeqnl
\end{defn}
\begin{remk}
Since the rotation number of $f_{(n)}(x)$ is
$(-1)^{n-1}\beta^n$, we have that 
$$(-1)^{n-1}\left(f_{(n)}(x)-x\right)>0$$
 for all $n\in\mathbb N$ and $x\in \mathbb R$. In particular, for
$x=0$ we obtain that $(-1)^{n-1}\lambda_{(n)}>0$.
Therefore, it follows that each function
$f_n(x)$, defined by \eqref{renormcirc},  is increasing in $x$ and 
 satisfies $f_n(x)<x$
\end{remk}
It is known \cite{deFariaI,Yampolsky01,Yampolsky02,Yampolsky03} 
that for every $k\in\mathbb N$  there is a 
universal constant  $-1<\lambda_k<0$
 and universal function $\eta_k$  such that:
\begin{itemize}
\item[RC1.] $\eta_k$ is increasing, $\eta_k(x)\leq x$, and
$\eta_k(x)=H_k(x^{2k+1})$ for some analytic function $H_k$.
\item[RC2.] The sequence of ratios 
$\alpha_n=\lambda_{(n)}/\lambda_{(n-1)}$
 converges to a limit $-1<\lambda_k<0$
\item[RC3.] For some $0<\delta_k<1$, 
$ \|f_n(x)-\eta_k(x)\|_{C_0}=O(\delta_k^n)$
where the  norm $\|\cdot\|_{C_0}$ is taken over an appropriate 
complex open domain $D_k$ that contains the real line $\mathbb R$.
\item[RC4.] The function  $\eta_k$ is a solution of the functional  
equations
\begin{equation}
\begin{array}{rcl}
\eta_k(x)&=&\frac1\lambda_k
 \eta_k\left(\frac1\lambda_k\eta_k(\lambda_k^2 x)\right)
\label{fixedeqn1}
\end{array}
\end{equation}
\begin{equation}
\begin{array}{rcl}
\eta_k(x)&=&
\frac{1}{\lambda_k^2} \eta_k\left(\lambda_k\eta_k(\lambda_k x)\right)
\label{fixedeqn2}
\end{array}
\end{equation}
\item[RC5.] The domain $D_k$ can be taken so
 that $\partial D_k$ is smooth and
\beqnl{setcirc1}
\overline{\lambda_k D_k}\subset D_k
\eeqnl
\beqnl{setcirc2}
\lim_{n\rightarrow\infty}\sup\{|\Imag(z)|: z\in D_k,\, |z|>n\}=0
\eeqnl
\begin{eqnarray}
\overline{\lambda^{-1}_k\eta_k(\lambda^2_k D_k)}&\subset& D_k
\label{setcirc3}\\
\overline{\lambda_k\eta_k(\lambda_k D_k)}&\subset& D_k
\label{setcirc4}
\end{eqnarray}
\end{itemize}

A few  useful relations can be obtain directly from
 \eqref{fixedeqn1} and \eqref{fixedeqn2}. For instance, by
letting $x=0$ we get
\beqn
\eta_k(1)=\lambda^2_k\qquad \eta_k(\lambda_k^2)=\lambda_k^3
\eeqn
Taking derivatives on \eqref{fixedeqn1} and \eqref{fixedeqn2}
 we have that
\beqn
\frac{\eta'_k(x)}{\eta'_k(\lambda_k^2 x)}&=&
\eta'_k\left(\frac{1}{\lambda_k}\eta_k(\lambda_k^2 x)\right)\\
\frac{\eta'_k(x)}{\eta'_k(\lambda_k x)}&=&
\eta'_k\left(\lambda_k\eta_k(\lambda_k x)\right)
\eeqn
Then letting $x\rightarrow0$ we get
\beqn
\eta'_k(1)&=&\frac{1}{\lambda^{4k}_k}\\
\eta'_k(\lambda^2_k)&=&\frac{1}{\lambda^{2k}_k}
\eeqn

Solutions of the equations \eqref{fixedeqn1}, \eqref{fixedeqn2} 
are constructed in \cite{EpsFixed} for all orders
of tangency at the critical point.

In the case of cubic circle maps, \cite{mestel,LanfordL} using
computer assisted proofs  constructed a fixed point. 
The exponential convergence
follows from the compactness of the derivative of the renormalization 
transformation at $\eta_1$, which is a consequence of 
the analyticity improving of the auxiliary functions. 
This is part of the conclusions of the computer assisted 
proofs.

\subsection{Renormalization theory of the noise}\label{cumopcirc}
In this section, we will develop in parallel two renormalization
theories for the noise. For the purposes of these paper, either
one will be enough so we will not
show to which extent they are equivalent.
 We will assume that the order of tangency
$k$ is fixed.

The renormalization group scheme for critical maps 
with rotation number $\beta$ gives information at small
scales and at Fibonacci times. Observe  that if $f$ is a circle
maps, then 
 Lyapunov functions $\widehat{\Lambda}^f$ 
and $\Lambda^f_p$ defined by \eqref{lambdahat} and \eqref{lambdap}
are periodic with period $1$.

 Denote by 
 $$k_n(x)=\Lambda^f_p(x, Q_n)$$
Since 
$f$ is increasing, equation \eqref{lambdarelation} 
for the Lyapunov functions $\Lambda_p$  implies  that
\begin{eqnarray}
k_n(x)&=&\left(f^\prime_{(n-1)}\circ f_{(n-2)}(x)\right)^p k_{n-2}(x)
+ k_{n-1}(f_{(n-2)}(x))\label{cumcirc1}\\
k_n(x)&=&\left(f^\prime_{(n-2)}\circ f_{(n-1)}(x)\right)^p k_{n-1}(x)
+ k_{n-2}(f_{(n-1)}(x))\label{cumcirc2}
\end{eqnarray}
\subsubsection{Lindeberg--Lyapunov operators}
To study the propagation of noise at small scales, we will use the 
following definition
\begin{defn}\label{deffcomopcirc}
Let $f$ be a critical circle map, and for  each $p\geq0$ and 
$n\in \mathbb N$, define the operators
\begin{equation*}
\begin{array}{rcl}
U_{n,p}h(z) &=& \left[f_{n-2}^\prime(\alpha_{n-2}
f_{n-1}(\alpha_{n-1} z))\right]^p h(\alpha_{n-1} z)\\
\noalign{\vskip6pt}
T_{n,p} q(z) &=&\left[f_{n-1}^\prime\left(\alpha^{-1}_{n-2}
f_{n-2}(\alpha_{n-1}\alpha_{n-2} z)\right)\right]^p
q(\alpha_{n-1}\alpha_{n-2}z)\\
\noalign{\vskip6pt}
R_{n} h(z) &=& h(\alpha^{-1}_{n-2}f_{n-2}(\alpha_{n-1}\alpha_{n-2} z))\\
\noalign{\vskip6pt}
P_n q(z) &=& q\left(\alpha_{n-2}f_{n-1}(\alpha_{n-1} z)\right)
 \end{array}
\end{equation*}
acting on the space or real analytic functions in the domain $D_k$.
For $p\geq0$, and $n\in\mathbb N$, 
 the \emph{Lindeberg--Lyapunov}  operators, 
$\mathcal K_{n,p}$
and $\widehat{\mathcal K}_{n,p}$,
acting on pairs of real analytic functions in $D_k$ are defined by 
the matrices of operators
\begin{equation}
\mathcal K_{n,p} =\left(
\begin{array}{ll} R_n & T_{n,p}\\
                  I   & 0 
\end{array}\right) \quad\widehat{\mathcal K}_{n,p} =\left(
\begin{array}{ll} U_{n,p} & P_n\\
                  I   & 0
\end{array}\right)\label{cumopcirc1}
\end{equation}
where $I$ is the identity map and $0$ is the zero operator.
Similarly, consider the operators acting on the 
space of real analytic functions on $D_k$ defined by 
\begin{equation*}
\begin{array}{rcl}
U_ph(z) &=& \left[\eta_k^\prime(\lambda_k
\eta_k\lambda_k z))\right]^p h(\lambda_k z)\\
\noalign{\vskip6pt}
T_p q(z) &=&\left[\eta_k^\prime\left(\lambda_k^{-1}
\eta_k(\lambda^2_k z)\right)\right]^p
q(\lambda^2_kz)\\
\noalign{\vskip6pt}
R h(z) &=& h(\lambda_k^{-1}\eta_k(\lambda_k^2 z))\\
\noalign{\vskip6pt}
P q(z) &=& q\left(\lambda_k \eta_k(\lambda_k z)\right)
 \end{array}
\end{equation*}
The Lindeberg--Lyapunov  operators
 $\mathcal K_p$ and $\widehat{\mathcal K}_p$ are defined by
\begin{equation}
\mathcal K_p =\left(
\begin{array}{ll} R & T_p\\
                  I   & 0 
\end{array}\right) \quad\widehat{\mathcal K}_p =\left(
\begin{array}{ll} U_p & P\\
                  I   & 0
\end{array}\right)\label{cumopfixed}
\end{equation}
\end{defn}

\bigskip
Notice from \eqref{cumcirc1} and \eqref{cumcirc2} that 
the change  variables 
$$x=\lambda_{(n-1)} z,\qquad \tilde{k}_n(z)=k_n(\lambda_{(n-1)}z),$$
implies that
\begin{eqnarray}
\left[\begin{array}{l}\tilde{k}_n\\ \tilde{k}_{n-1}
    \end{array}
\right]
&=& \mathcal K_{n,p}\cdots\mathcal K_{3,p}
\left[\begin{array}{l}\tilde{k}_2\\ \tilde{k}_1
    \end{array}\right]\label{circsmallscale}\\
\left[\begin{array}{l}\tilde{k}_n\\ \tilde{k}_{n-1}
    \end{array}\right]&=& \widehat{\mathcal K}_{n,p}\cdots
\widehat{\mathcal K}_{3,p}
\left[\begin{array}{l}\tilde{k}_2\\ \tilde{k}_1
    \end{array}\right]\label{circsmallscale2}
\end{eqnarray}

\begin{remk}
Notice that  the derivatives involved in Definition
\ref{deffcomopcirc} are positive. Therefore,
 the notions of cumulant
operators and Lindeberg--Lyapunov operators coincide
in this case.
\end{remk}

Equations \eqref{circsmallscale} and \eqref{circsmallscale2}
are the analogs to  \eqref{feig_smallscale}
for  period doubling. These equations measure the 
growth of the linearized propagation of noise at
Fibonacci times for an orbit starting at the critical point $0$.

\subsubsection{Exponential convergence of the Lindeberg--Lyapunov
operators}
An important consequence of the exponential convergence 
of $f_n$ to $\eta_k$ is that 
the  Lindeberg--Lyapunov operators $\mathcal K_{n,p}$ and
 $\widehat{\mathcal K}_{n,p}$
will converge exponentially fast to operators $\mathcal K_p$ and 
$\widehat{\mathcal K}_p$ respectively, as $n\rightarrow\infty$.
Indeed, we have that
\begin{lem}\label{circconvexp}
Let $f$ be a circle map of order $k$ satisfying G1--G3.
The operators $\mathcal K_{n,p}$, $\mathcal K_p$,
$\widehat{\mathcal K}_{n,p}$ and $\widehat{\mathcal K}_p$
acting on the space of pairs of bounded  analytic functions defined
on some compact set $B_k\subset D_k$ that  has smooth boundary and 
contains $[-1,1]$ in its interior, are compact.
Furthermore,  for $p$ fixed, 
\begin{eqnarray}
\|\mathcal K_{n,p}-\mathcal K_p\|&\leq& c_p
\left(\|f_n-\eta_k\|_{B_k}
+ \|f_{n-1}-\eta_k\|_{B_k}\right)\label{knp}\\
\|\widehat{\mathcal K}_{n,p}-
\widehat{\mathcal K}_p\|&\leq& c_p
\left(\|f_n-\eta_k\|_{B_k}
+ \|f_{n-1}-\eta_k\|_{B_k}\right)\label{hatknp}
\end{eqnarray}
\end{lem}
\proof
We  only show \eqref{hatknp}, since the proof of \eqref{knp}
is very similar.

For any $A\subset \mathbb C$, recall the notation 
$A^{\epsilon}=\{z\in \mathbb C: d(z,A)\leq \epsilon\}$.
Since $|\alpha_{n-1}-\lambda_k|=O(\delta_k^n)$ for some
$0<\delta_k<1$, we have that for any compact $B_k\subset D_k$
and for all $n$ large enough 
$$\overline{\alpha_{n-1}B_k}\subset 
\overline{\lambda_kB_k}^\epsilon \subset D_k$$
where  $\epsilon>0$ is small.

For $z\in D_k$, we have that
\begin{multline*}
|\alpha_{n-2}f_{n-1}(\alpha_{n-1}z)-\lambda_k\eta_k(\lambda_k z)|
\leq |\alpha_{n-2}||f_{n-1}(\alpha_{n-1}z)- 
\eta_k(\alpha_{n-1} z)|\\
+|\alpha_{n-2}-\lambda_k||\eta_k(\alpha_{n-1}z)|
+ |\lambda_k||\eta_k(\alpha_{n-1}z)-\eta_k(\lambda_kz)|
\end{multline*}
Therefore, using Cauchy estimates, we have that for any compact set 
 $B_k\subset D_k$ containing $0$ 
\begin{multline}
\|\alpha_{n-2}f_{n-1}\circ\alpha_{n-1} - 
\lambda_k\eta_k\circ\lambda_k\|_{B_k}\leq
 C\left(\|f_{n-1}-\eta_k\|_{B_k}\right.\label{bcirc5}\\
+\left. \|f_{n-2}-\eta_k\|_{B_k}\right)
\end{multline}
Therefore, we by taking $\epsilon$ small enough, we have that 
for all $n$ large enough
\beqnl{setcirc5}
\overline{\alpha_{n-2}f_{n-1}(\alpha_{n-1}B_k)}^\epsilon
\subset \overline{\lambda_k\eta_k(\lambda_k B_k)}^\epsilon
\subset D_k
\eeqnl
Consider  $B_k$ with smooth boundary  and let 
$W_k=\overline{\lambda_k\eta_k(\lambda_k B_k)}^\epsilon$. 
For  $z\in B_k$ we have
\begin{multline}
|f'_{n-2}(\alpha_{n-2}f_{n-1}(\alpha_{n-1} z))
-\eta'_k(\lambda_k\eta_k\lambda_k z))|
\leq \left|f'_{n-2}(\alpha_{n-2}f_{n-1}(\alpha_{n-1} z))
\right.\nonumber\\
\left. -\eta'_k(\alpha_{n-2}f_{n-1}(\alpha_{n-1} z))\right|\nonumber\\
+|\eta'_k(\alpha_{n-2}f_{n-1}(\alpha_{n-1} z))-
\eta'_k(\lambda_k\eta_k\lambda_k z))|
\end{multline}
Then, from \eqref{setcirc5} and \eqref{bcirc5} we get
\begin{multline*}
\|f'_{n-2}\circ(\alpha_{n-2}f_{n-1}\circ\alpha_{n-1})
-\eta'_k\circ(\lambda_k\eta\circ\lambda_k )\|
\leq |f'_{n-2} -\eta'_k|_{W_k} \\
+C\|\eta''_k\|_{W_k}\|\left(
\|f_{n-1}-\eta_k\|_{B_k}+\|f_{n-2}-\eta_k|_{B_k}\right)
\end{multline*}
Using Cauchy estimates we obtain
\begin{multline}
\|f'_{n-2}\circ(\alpha_{n-2}f_{n-1}\circ\alpha_{n-1})
-\eta'_k\circ(\lambda_k\eta\circ\lambda_k )\|_{B_k}
\leq c\left( \|f_{n-1}-\eta_k\|_{D_k}\right.\label{bcirc6}\\
\left.+\|f_{n-2}-\eta_k|_{D_k}\right)
\end{multline}
Denote by $\mathcal O=\mathbb C\setminus\{x+iy:\, x\leq0,\, y=0\}$, 
and let $B_k\subset D_k$ be a bounded domain containing
 $[-1,1]$ with smooth boundary
such that 
\beqn
 \overline{\lambda_k\eta_k(\lambda_k B_k)}
\subset
\mathcal O
\eeqn
From \eqref{bcirc6} we can assume, by taking $\epsilon$ smaller if 
necessary, that
\beqn
\overline{\alpha_{n-2}f_{n-1}(\alpha_{n-1}B_k)}
\subset\overline{\lambda_k\eta_k(\lambda_k B_k)}^{\epsilon}
\subset\mathcal O
\eeqn
Recall the operators $U_{n,p}$, $U_p$, $P_{n,p}$ and $P_p$ defined
in Definition \ref{deffcomopcirc}. We restrict these  operators
to the space of bounded  analytic functions on $B_k$.
 Denote by
$Y_k=\overline{\lambda_k\eta_k(\lambda_k B_k)}^{\epsilon}$
and let  $\psi(z)=z^p$ defined on $Y_k$. Hence, using the 
triangle inequality and Cauchy estimates, we get for any bounded
analytic functions $h$, $q$ on $B_k$ that 
\begin{multline}
\|(U_{n,p}-U_p)h\|_{B_k}\leq 
\|\psi\|_{Y_k}\|h\|_{\lambda_kB_k}|\alpha_{n-1}-\lambda_k|\label{bcirc7}\\
+\tilde{C}\|\psi'\|_{Y_k}\|h\|_{\lambda_kB_k}\left(
\|f_{n-1}-\eta_k\|_{B_k}+\|f_{n-2}-\eta_k\|_{B_k}\right)
\end{multline}
\begin{multline}
\|(P_{n,p}-P_p)q\|_{B_k}\leq C\|q'\|_{W_k}\left(
\|f_{n-1}-\eta_k\|_{B_k}\right.\label{bcirc8}\\
+ \|f_{n-2}-\eta_k\|_{B_k}
\end{multline}
Combining \eqref{bcirc7} and \eqref{bcirc8} we obtain
\eqref{hatknp}.

For $z_1,\,z_2\in B_k$ close enough we have that
\beqn
|T_{n,p}h(z_1)-T_{n,p}h(z_2)|&\leq& \left\|\left(\psi\circ
 f'_{n-2}(\alpha_{n-2}f_{n-1}\circ\alpha_{n-1})\right)^\prime
\right\|_{B_k}
|z_1-z_2|\\
|P_{n,p}h(z_1)-P_{n,p}h(z_2)|&\leq&
\left\|\left(q\circ(\alpha_{n-2}f_{n-1}\circ\alpha_{n-1})\right)^\prime
\right\||z_1-z_2|
\eeqn
Therefore, the compactness of the operators 
$\widehat{\mathcal K}_{n,p}$ follows
from Cauchy estimates and the Arzela--Ascoli theorem.
\endproof
As in section \ref{cumulantoper}, the spectral properties of
the operators $\mathcal K_{n,p}$,\, $\mathcal K_p$,
 $\widehat{\mathcal K}_{n,p}$ and $\widehat{\mathcal K}_p$ will
prove important for our analysis of the Gaussian properties of the
scaling limit.

Notice that for $p\in\mathbb N$, since
$f_n$ and $\eta_k$ are increasing,  the cumulant functions
defined by \eqref{cumsuma} and the Lyapunov functionals \eqref{lambdap}
coincide.  Hence, the study of the properties of cumulants
is equivalent to that of  the properties of the Lindeberg--Lyapunov
functionals.

\subsection{Spectral analysis of Lindeberg--Lyapunov operators}
\label{circcumopspect}
{From}  Lemma \ref{circconvexp},
we know that the Lindeberg--Lyapunov  operators
${\mathcal K}_p$, $\widehat{\mathcal K}_p$ when defined on 
in a space space of  real valued analytic functions with domains in 
$D$ that satisfy \eqref{setcirc1}, \eqref{setcirc2}, \eqref{setcirc3},
\eqref{setcirc4} are such that their squares are compact. 
Again, the compactness follows from the fact that 
  the operators 
$U_{n,p}, T_{n,p}, R_n, Q_n$ in Definition \ref{deffcomopcirc}
are analyticity improving.

The compactness of the operators ${\mathcal K}_{n,p}$, 
$\widehat {\mathcal K}_{n,p}$,  is slightly more subtle 
since they have a block (in the lower diagonal) 
that is the identity. Using the fact that this block
appear in the lower diagonal, the operator can be made 
compact by considering that the second component 
is analytic in a slightly smaller domain. In this way, the 
identity becomes an immersion, which is compact. 
If the domain for the second component is chosen 
only slightly smaller, the analyticity improving properties of 
the other operators are still maintained. 

\begin{remk} An slightly different approach is 
to remark that the square of the operators
consist only of compact operators. 
For our purposes, the study of the squares will 
be very similar. 
\end{remk}

Notice that  the operators $\mathcal K_{n,p}, 
\widehat {\mathcal K}_{n,p}$
preserve the cone $\mathcal C$ of pairs of functions, both of whose components 
are strictly positive  when restricted to the reals. 
Note that the interior of the cone $\mathcal C$ are the 
pairs of functions whose components  are strictly positive when 
restricted to the reals. 

It is also easy to verify the properties (i), (ii), (iii) 
in Section~\ref{sec31}.
Again, (i), (ii) are obvious and we note that if a pair is not identically
zero, then, one of the components has to be strictly positive in 
an interval. The structure of the square of the operator implies 
that the component is strictly positive in an interval which is 
larger. This, in turn, implies that the second component 
becomes positive. Both components 
become strictly positive by repeating the operation  a finite number of times.

Therefore, we can apply the  
the  Kre\v{\i}n--Rutman theorem and obtain a   result similar to
Proposition \ref{KR}.
\begin{prop}\label{KR2} 
Let
 $\mathcal K_{n,p}$ and $\widehat{\mathcal K}_{n,p}$ be the
Lindeberg--Lyapunov operators defined by \eqref{cumopcirc1}.
Denote by  $\mathcal K_{\infty,p}=\mathcal K_p$ and
$\widehat{\mathcal K}_{\infty,p}=\widehat{\mathcal K}_p$. 
For all $p\geq0$ and $n\in\mathbb N\cup\{\infty\}$,
let $\rho_{n,p}$ ($\widehat{\rho}_{n,p}$) be the spectral radius of
$\mathcal K_{n,p}$ ($\widehat{\mathcal K}_{n,p}$). Then,
\begin{itemize}
\item[{\it a)}]  $\rho_{n,p}$  ($\widehat{\rho}_{n,p}$) is a positive
eigenvalue of $\mathcal K_{n,p}$ ($\widehat{\mathcal K}_{n,p}$).
\item[{\it b)}]  The rest of 
$\spec(\mathcal K_{n,p})\setminus\{0\}$  
($\spec(\widehat{\mathcal K}_{n,p})\setminus\{0\}$)
consists of eigenvalues $\mu$ with  $|\mu|<\rho_{n,p}$
($|\mu|<\widehat{\rho}_{n,p}$).
\item[{\it c)}] A pair of positive functions
$(\psi_{n,p},\phi_{n,p})$ is an eigenvector of $\mathcal K_{n,p}$
($\widehat{\mathcal K}_{n,p}$) if and only if the corresponding 
eigenvalue is $\rho_{n,p}$ ($\widehat{\rho}_{n,p}$).
\end{itemize}
\end{prop}
\subsubsection{Properties of the spectral radius of the 
Lindeberg--Lyapunov operators}
The convexity properties of the Lyapunov functions $\Lambda_p$ will
also imply similar properties for the cumulants $\mathcal K_{n,p}$
and $\mathcal K_{n,p}$. In the next Theorem, we will only consider the
the operators $\mathcal K_p$ and $\widehat{\mathcal K}_p$.
We also remark that the argument in Proposition~\ref{uniqueness}
applies in our case also, so that the spectrum is 
independent of the domain considered. Hence, we 
will not include the domain in the notation.

\begin{thm}\label{comparison2}
Let $\rho_p$  be the spectral radius of
$\mathcal K_p$. ( Similarly for $\widehat{\rho}_p$ and
  $\widehat{\mathcal K}_p$.) 
\begin{itemize}
\item[K1.] For all $p>0$
$$ \lambda_k^{2kp} \rho_p > 1$$
\item[K2.]The map  $p\mapsto \rho_p$ is strictly log--convex
\item[K3.] The map $p\mapsto \log(\rho_p)/p$ is strictly decreasing.
\end{itemize}
A Similar result holds for $\widehat{\rho}_p$.
\end{thm}

\proof (K1) 
For any pair of positive functions $[h,\,\ell]$, we have that
\beqn
\mathcal K^{2m}_p[h,\,\ell] > \left[ T^m_p h, T^m_p \ell\right]
\eeqn
If $[h,\,\ell]$ is the dominant eigenvector of $\mathcal K_p$, 
then we have that
\beqn
\rho^{2m}_p[h(z),\,\ell(z)]& >& \left[ T^m_p h, T^m_p \ell\right](z)\\
&=& 
\prod_{j=1}^m\left(\eta'\left(
\lambda_k^{-1}\eta_k\lambda_k^{2j}z)\right)\right)^p
[h(\lambda_k^{2mz}), \ell(\lambda_k^{2m}z)]
\eeqn
Since   $\eta_k'(1)=\lambda_k^{-4k}$, the 
conclusion of (M1) follows by letting  $z=0$.

\noindent
To prove [K2.] [K3.] we proceed
in a similar way as in the proof of Theorem~\ref{Comparison}.
We study the eigenvalue equation for the leading eigenvalue
and use it to construct positive functions that satisfy
eigenvalue equations with a positive remainder and then, we 
apply Proposition~\ref{KRbound}. We suppress the subindex $n$
for typographical clarity.  We will use the notations for 
operators introduced in Definition \ref{deffcomopcirc}.

The eigenvalue equation for the leading 
eigenvalue of ${\mathcal K}_p$ is
\beqnl{coordinates} 
\rho_p \ell_p &=& R h_p  + T_p  \ell_p  \\
\rho_p \ell_p &=& h_p
\eeqnl

We raise the equations \eqref{coordinates} to the power $\alpha > 1$
and use the inequality for the binomial 
theorem and  the elementary identities
 $(T_{p} \ell_p)^\alpha = T_{p \alpha} \ell_p^\alpha$, 
$(R h_p)^\alpha = R h_p^\alpha$. Then, we have:
\beqnl{coordinates2} 
\rho_p^\alpha  \ell_p^\alpha &=& (R h_p  + T_p  \ell_p)^\alpha
\\
& > & (R h_p)^\alpha  + (T_p  \ell_p)^\alpha\\
& = & R h_p^\alpha + T_{p \alpha} \ell_p^\alpha \\
\rho_p^\alpha \ell_p^\alpha &=& h_p^\alpha
\eeqnl

Hence, we obtain 
\[
{\mathcal K}_{\alpha p}[h^\alpha, \ell^\alpha] \ge \rho_p^\alpha 
[h^\alpha, \ell^\alpha]
\]
and we have that the inequality is strict in the first component. 
Therefore
\[
{\mathcal K}_{\alpha p}^2 [h^\alpha, \ell^\alpha] > 
 \rho_p^{ 2 \alpha} [h^\alpha, \ell^\alpha]
\]
Using that the leading eigenvalue of ${\mathcal K}_p^2$ is 
the square of the leading eigenvalue of ${\mathcal K}_p$, we
obtain, applying  Proposition~\ref{KRbound} that 
\[
\rho_{p \alpha} < \rho_p^\alpha
\]
which is property [K3.]. 

To prove property [K2], we multiply the eigenvalue
equation for leading eigenvectors of 
${\mathcal K}_p$, $\widehat {\mathcal K}_p$ and raise 
to the $1/2$ power.   We use the identities 
$(R h_p) (R h_q )   = R (h_p h_q)$, 
$[ R (h_p h_q)]^{1/2} = R \left[ (h_p h_q)^{1/2}\right] $
as well as similar identities for the operator $T$.

\beqn
\begin{split}
(\rho_p \rho_q)^{1/2} (h_p h_q)^{1/2} &= 
\left[ (R h_p + T_p \ell_p)\cdot
(R h_q + T_q \ell_q ) \right] \\
& >( R h_p R h_p + T_p \ell_p T_q \ell_q)^{1/2} \\
&> R( h_p h_q)^{1/2} + T_{(p +q)/2} (\ell_p \ell_q)^{1/2} 
\end{split}
\eeqn

Multiplying the second components we have 
\[
(h_p h_q)^{1/2}  = (\rho_p \rho_q)^{1/2}
(\ell_p \ell_q)^{1/2}
\]

In other words, we have 
\[
{\mathcal K}_{(p+q)/2} [ (h_p h_q)^{1/2},  (\ell_p \ell_q)^{1/2} ]
\le (\rho_p \rho_q)^{1/2}  [ (h_p h_q)^{1/2},  (\ell_p \ell_q)^{1/2}  ]
\]
with the inequality being strict in the first component. 
Proceeding as before, we obtain 
\[
\rho_{(p +q)/2 } <  (\rho_p \rho_q)^{1/2}
\]
which is equivalent to strict log-convexity.

\endproof

\subsection{Asymptotic properties of the renormalization}
A consequence of Lemma \ref{circconvexp} is that
the asymptotic properties of the operators $\mathcal K_{n,p}$
and $\widehat{\mathcal K}_{n,p}$ become similar to those of
$\mathcal K_p$ and $\widehat{\mathcal K}_p$ respectively.

The following result, Corollary \ref{kurtoQnresult},  is 
obtained as an application of  Proposition  \ref{productalignment}.
Since the proof is very similar to that of
Corollary \ref{kurto2n} for period doubling,
 we will omit the details.

\begin{cor}\label{kurtoQnresult}
Let $\mathcal K_{n,p}$,  $\mathcal K_p$, and 
$\rho_p$ be  as in Proposition \ref{KR}. Then
\begin{itemize}
\item[1.] There is a constant $c_p>0$  such that for all
 positive analytic pairs of  functions $[h,q]$ in
the domain $\mathcal D$  of the Lindeberg--Lyapunov operators
\beqnl{alignmentcirc}
c^{-1}_p\rho_p^n \leq \mathcal K_{n,p}
\cdots\mathcal K_{1,p}[h(z),q(z)]
\leq c_p\rho_p^n
\eeqnl
\item[2.] For $p>2$, we have  $\rho_p< (\rho_2)^{p/2}$
\item[3.] For any $p>2$
\beqn
\lim_{n\rightarrow\infty}\sup_{x\in I_{n}} 
\frac{\Lambda_p^f(x,Q_n)}
{\{\Lambda_2^f(x, Q_n)\}^{p/2}}= 0
\eeqn
where $I_m=[-|\lambda_{(n)}|, |\lambda_{(n)}|]$.
\end{itemize}
\end{cor}
\subsection{Proof of the central limit theorem for circle maps
(Theorems \ref{circclt} and \ref{circclt2}).}
\label{sec:circ_clt}
In this section, we will prove a 
 central limit for critical circle maps  (Theorem
\ref{circclt}).  The method of the proof is  similar 
to the method developed in 
 Section \ref{sec:feig_clt} for period doubling.
That is, we will verify the Lyapunov condition \eqref{lyapcon}
along the whole sequence of integers and then use Theorems
\ref{main} and \ref{berry-esseen}.

The renormalization theory developed in Section \ref{circcumopspect}
gives control on $\Lambda_s(x,Q_m)$ for all
$x$ sufficiently small ($|x|\leq |\lambda_{(m+2)}|$),
see Corollary \ref{kurtoQnresult}. The main tool to obtain control
of the noise on a segment of the orbit will be to
decompose the segment into pieces where renormalization applies, that 
is segments of Fibonacci length.

This decomposition is develop in Section \ref{decompfiboblock}. 
The main conclusion of this section is that the effect of the 
noise over a period of time $n$ equals the sum of
Fibonacci blocks -- where renormalization applies -- with
some weights. See \eqref{decompositionofcirclambda}. This weights
measure how the effect of Fibonacci blocks propagates to
the end of the interval.

The effect of weights on the Fibonacci blocks is measure in 
Section \ref{circweigths}.

\subsubsection{Fibonacci decomposition}\label{decompfiboblock}
Given $n\in \mathbb N$, it admits a unique decomposition
\beqnl{fibodecomp}
n=Q_{m_0} + \cdots + Q_{m_{r_n}}
\eeqnl
where $m_0> \ldots > m_{r_n}>0$  are  non--consecutive integers
and $Q_{m_j}$ is the $m_j$--th Fibonacci number. If necessary,
we will use $m_j(n)$ to emphasize the dependence on $n$.
Notice that $r_n\leq m_0(n)\leq[\log_{\beta^{-1}}n+1]$.

We have the following Fibonacci decomposition for  $\Lambda_p(x,n)$
\beqnl{decompositionofcirclambda}
\Lambda_p(x,n)=\sum^{r_n}_{j=0}\left|\Psi_{j,n}(x)\right|^p\Lambda_p\left(
f^{n-Q_{m_j}-\cdots -Q_{m_r}}(x), Q_{m_j}\right)
\eeqnl
where 
\beqnl{psicircderivatives}
\Psi_{j,n}(x)&=&
\left(f^{n-Q_{m_0}-\cdots-Q_{m_j}}\right)^\prime\circ
\left(f^{Q_{m_j}+\cdots Q_{m_0}}\right)(x)
\eeqnl

Since  $(f^m)^\prime(x)$ is periodic of 
period $1$ for any $m\in\mathbb N$,   $\Lambda_p(x,n)$ can be
expressed in terms of the functions $f_{(n)}$ 
defined by \eqref{newrotnumbermap}. Indeed, let
\beqnl{retcirck}
\upsilon_j=f_{(m_j)}(\upsilon_{j-1})\qquad\upsilon_{-1}=x
\eeqnl
then,
\beqn
\Psi_{n,j}=\prod^{r_n}_{l=j+1} (f_{(m_l)})^\prime(\upsilon_{l-1})
\eeqn
and
\beqnl{lyapdecompcirc}
\Lambda_p(x,n)=\sum^{r_n}_{j=0}
 (\Psi_{n,j}(x))^p\Lambda_p(\upsilon_{j-1},Q_{m_j})
\eeqnl

Let us fix a critical map of the circle of order $k$, and let 
$0<\epsilon\ll 1$. For all   $n$ large enough we have that
\beqnl{closeness2}
\|f_n-\eta_k\|_{C(D_k)}<\epsilon
\eeqnl
Using renormalization, we will control the
size of the size of the weights 
 $\Psi_{j,n}(x)$ for 
$|x|\leq |\lambda_{(m_0+2)}|$.
\subsubsection{Estimation of the weights} \label{circweigths}
In this section we estimate the weights $\Psi_{n,j}$ defined in
\eqref{lyapdecompcirc}. The method used is very similar to
to one develop for period doubling.

Recall that 
 $\eta_k(x)=\lambda_k + \sum_{j=1}^\infty b_j x^{(2k-1)j}$,
with $b_1>0$. 
Consider the function $H_k$ defined by
\beqn
H_k(x)=\frac{\eta'_k(x)}{x^{2k}}
\eeqn
and define
\beqnl{infsup}
s_k=\inf_{\{x:|x|\leq\lambda^2_k\}} H_k(x) \qquad
u_k=\sup_{\{x:|x|\leq\lambda^2_k\}} H_k(x)
\eeqnl
It follows that 
\beqnl{etabound}
c x^{2k}\leq \eta'_k(x)\leq d x^{2k}
\eeqnl

In the rest of this section, we will assume with no loss of generality that
all $n$ are large enough so that
\beqnl{closecirc}
\sup_{\{x:|x|\leq 1\}}|f_n(x)-\eta_k(x)|<\epsilon
\eeqnl
with $0<\epsilon\ll 1$.

\begin{lem}\label{growthcircle}
Let $\{m_j\}^r_{j=0}$ be a decreasing sequence of non--consecutive 
positive integers. For $|x|\leq \lambda_{(m_0+2)}$, let
  $\{\upsilon_{j}\}_{j=-1}^r$ be defined by \eqref{retcirck}.
\begin{itemize}
\item[(1)] For all $0\leq j\leq r$ for which $m_j$ is large enough
\beqnl{control1}
|\lambda^3_k||\lambda_{(m_j-1)}|\leq |\upsilon_j| \leq 
|\eta_k(-\lambda^2_k)| |\lambda_{(m_j-1)}|
\eeqnl
\item[(2)]  Let $c= s_k-\epsilon$ and 
$d=u_k+\epsilon$. Then,
\beqnl{control2}
c \lambda^{6k}
\left(\frac{\lambda_{(m_{j-1}-1)}}{\lambda^{{(m_j-1)}}}\right)^{2k}
\leq (f_{(m_j)})^\prime(\upsilon_{j-1})
\leq d \left(\frac{\lambda_{(m_{j-1}-1)}}{\lambda_{(m_j-1)}}\right)^{2k}
\eeqnl
\end{itemize}
\end{lem}
\proof
(1) For $j=0$, we have  that
$$ \upsilon_0= \lambda_{(m_0-1)}f_{m_0}\left(
\frac{\upsilon_{-1}}{\lambda_{(m_0-1)}}\right)$$
The exponential convergence of $f_m$ to $\eta_k$ implies that
$$\left|\frac{\upsilon_{-1}}{\lambda_{(m_0-1)}}\right|\leq 
\left|\frac{\lambda_{(m_0+2)}}{\lambda_{(m_0-1)}}\right|\leq 
|\lambda_k|^3+\epsilon$$
for some $0<\epsilon\ll 1$. Hence, 
$$(|\lambda_k|^3-\epsilon)|\lambda_{(m_0-1)}|
\leq |\upsilon_0|\leq |\eta_k(-\lambda_k^2)||\lambda_{(m_0-1)}|
$$

By induction, assume that
$$|\upsilon_{j-1}|< |\eta_k(-\lambda_k^2)|\left|\lambda_{(m_{j-1}-1)}\right|$$
Notice that
$$\upsilon_j=f_{(m_j)}(\upsilon_{j-1})=
\lambda_{(m_j-1)}f_{m_j}\left(
\frac{\upsilon_{j-1}}{\lambda_{(m_j-1)}}\right)
$$
Since $m_{j-1}-m_j\geq 2$, we have that
\beqnl{lesslambda2}
\left|\frac{\upsilon_{j-1}}{\lambda_{(m_j-1)}}\right|\leq
|\eta_k(-\lambda_k^2)| (\lambda_k^2 +\epsilon)\leq \lambda_k^2
\eeqnl
since $\epsilon\ll 1$. Therefore, we have that
$$|\lambda_k|^3\left|\lambda_{(m_j-1)}\right|
\leq |\upsilon_j|\leq |\eta_k(-\lambda_k^2)|
\left|\lambda_{(m_j-1)}\right|$$

\bigskip
\noindent
(2) Taking the neighborhood of the fixed point 
smaller if necessary, assume that
for all $m$ large enough
$$\sup_{|x|\leq\lambda_k^2}\left|\frac{(f_m)^\prime(x)}{x^{2k}} - 
\frac{\eta_k^\prime(x)}{x^{2k}}\right|<\epsilon
$$
Notice that
$$(f_{(m_j)})^\prime(\upsilon_{j-1})=
 (f_{m_j})^\prime\left(\frac{\upsilon_{j-1}}{\lambda_{(m_j-1)}}\right)
$$
If  $m_j$ is large enough, by the exponential convergence 
of the renormalized functions to the fixed point, we have that
$$\left|\frac{\upsilon_{j-1}}{\lambda_{(m_j-1)}}\right|\leq 
\lambda_k^2+\epsilon$$
Therefore,  from \eqref{etabound} we get
\beqn
 c\lambda_k^{6k}
\left(\frac{\lambda_{(m_{j-1}-1)}}
{\lambda_{(m_j-1)}}\right)^{2k}
\leq \left|(f_{(m_j)})^\prime(\upsilon_{j-1})\right|\leq d
\left(\frac{\lambda_{(m_{j-1}-1)}}
{\lambda_{(m_j-1)}}\right)^{2k}
\eeqn
\endproof
\subsubsection{Lyapunov condition for critical circle maps}\label{techcircsec}
In this section we show that the Lyapunov condition \ref{lyapcon}
holds for subsequences of the orbits 
starting at points
of the form $x=f^l(0)$, $l\in\mathbb N$.
The main result of this section is Theorem
\ref{lyapconorbit0circ} which proves 
that the Lyapunov condition holds always if 
we choose appropriately subsequences (they are 
subsequences of numbers with very few terms in the 
Fibonnaci expansion. 
The precise conditions depend on numerical properties 
of the fixed point. We show that if the fixed point satisfies
some conditions, we can obtain the limit along the full sequence. 

First, we  use Lemma \ref{growthcircle} to
estimate the growth  of the Lyapunov functions $\Lambda_p(0,n)$
 at zero at all times.

\begin{prop}\label{kurto0ncirc}
Let $c$ and $d$ as in Proposition \ref{growthcircle}. For each
integer $n$, let $Q_{m_0(n)}$  be the largest Fibonacci number 
in the Fibonacci expansion of $n$
and  let $r_n$ be the number of terms in 
expansion \eqref{fibodecomp} of $n$. Denote
by $I_{m_0}=[-|\lambda_{(m_0+2)}|,|\lambda_{(m_0+2)}|]$.
\begin{itemize}
\item[1.] For each  $p>0$, there are constants $a_p$ and $b_p$ such
that
\beqnl{lambdacircgrowthat0}
a_p\left(\frac{c\,\lambda^{6k}_k}{\lambda_k^{2km_{r_n}}}\right)^p
\leq
\frac{\Lambda_p(x,n)}{(\lambda_k^{2kp}\rho_p)^{m_0(n)}}
\leq
b_p\left(\frac{d^{r_n}}{\lambda_k^{2k m_{r_n}}}\right)^p
\eeqnl 
for all $x\in I_{m_0(n)}$.
\item[2.] If $\{n_i\}$ is an increasing sequence of
integers such that
\beqnl{laccirc}
\lim_{i\rightarrow\infty}\frac{r_{n_i}}{m_0(n_i)}=0
\eeqnl
then
\beqnl{circcurtoratiolac}
\lim_{j\rightarrow\infty}
\sup_{x\in I_{m_0(n_j)}}\frac{\Lambda_p(x,n_j)}{(\Lambda_2(x,n_j))^{p/2}}
= 0
\eeqnl
\item[3.] If the following condition
\beqnl{hypothesiscirc}
\left(s_k \lambda_k^{6k}\right)^p \lambda_k^{2kp}\rho_p >1,
\eeqnl
 holds ($s_k$ defined as in \eqref{infsup}) then,
 for each $p>2$
\beqnl{circcurtoratio}
\lim_{n\rightarrow\infty}
\sup_{x\in I_{m_0(n)}}\frac{\Lambda_p(x,n)}{(\Lambda_2(x,n))^{p/2}}
= 0
\eeqnl
\end{itemize}
\end{prop}
\proof
For a given $n\in\mathbb N$, 
 let $\{m_j\}_{j=0}^{r_n}$ be the decreasing sequence of numbers
  that appear in the 
Fibonacci expansion  \eqref{fibodecomp} of $n$.

\noindent
(1) Recall the weights  $\Psi_{n,i}$ defined in \eqref{lyapdecompcirc}.
By  \eqref{control2} in Lemma \ref{growthcircle} and 
\eqref{closecirc}, we can assume without  loss of generality that
all $m_j$ are large enough, so that
\beqn
(c\lambda^{6k}_k)^{r-j}
\leq \Psi_{n,j}(x)
\left(\frac{\lambda_{(m_{r_n}-1)}}{\lambda_{(m_j-1)}}\right)^{2k}
\leq d^{r-i}
\eeqn
for all $x\in I_{m_0(n)}$. 

The exponential convergence of $f_n$ to $\eta_k$, 
implies that there is a constant $C>0$ such that
\beqnl{boundlambdacirc1}
C^{-1}(c\lambda^{6k}_k)^{r_n-j+1} 
\leq \lambda_k^{2k(m_{r_n}-m_j)}\Psi_{n,j}(x)\leq 
C\,d^{r_n-j+1}
\eeqnl
By Corollary \ref{kurtoQnresult}, there is a constant $D>0$
such that 
\beqnl{boundlambdacirc2}
D^{-1}\leq \rho^{-m_j}_p\Lambda_p(\upsilon_{j-1},Q_{m_j})< D
\eeqnl
for all $j=0,\ldots r_n$. 

Since each term in the decomposition \eqref{lyapdecompcirc} 
of $\Lambda_p(0,n)$ is positive, for some constants $a_p$ and
$b_p$ we have that
\beqn
a_p (c\lambda^{6k}_k)^{p\,r_n}\left(\frac{\lambda_k^{m_0}}
{\lambda_k^{m_{r_n}}}\right)^{2kp}\rho^{m_0}_p
 \leq \Lambda_p(x,n)
\leq b_p
\sum_{j=0}^{r_n} d^{p(r_n-j)} \left(\frac{\lambda_k^{m_j}}
{\lambda_k^{m_{r_n}}}\right)^{2kp}\rho^{m_j}_p
\eeqn
for all $x\in I_{m_0(n)}$.
We obtain \eqref{lambdacircgrowthat0} from
\eqref{boundlambdacirc1}, \eqref{boundlambdacirc2} 
by noticing that $m_0-m_j\geq 2j$.

Notice that if $\{n_i\}$ is a sequence of
 integers that satisfy \eqref{laccirc}
then
\beqnl{biglambdacirclac}
\lim_{i\rightarrow\infty}\inf_{x\in I_{m_0(n_i)}}\Lambda_p(x,n_i)
=\infty
\eeqnl

On the other hand, if the condition \eqref{hypothesiscirc} holds,
then $(c\lambda^{6k}_k)^p\lambda_k^{2kp}\rho_p>1$. Therefore,
\beqnl{biglambdacirc}
\lim_{n\rightarrow\infty}\inf_{x\in I_{m_0(n)}}\Lambda_p(x,n)
=\infty
\eeqnl

\bigskip
\noindent
(2) Corollary \ref{kurtoQnresult} implies that  for any $\epsilon>0$, 
there is  an integer $M_\epsilon$ such that if $m\geq M_\epsilon$, 
then
 $$ \Lambda_p(x, Q_m)<\epsilon 
\left(\Lambda_2(x, Q_m)\right)^{p/2}$$
for all $x\in [-|\lambda_{(m)}|, |\lambda_{(m)}|]$.
 Let  $\mathcal S_{r_n}=\{0,\ldots,r_n\}$ and define
$$
\mathcal A^1_n=\left\{j\in S_{r_n}:
\sup_{x\in I_{m_j}}\frac{ \Lambda_p\left(x, Q_{m_j}\right)}
{\left(\Lambda_2\left(x, Q_{m_j}\right)\right)^{p/2}}
\leq\epsilon\right\}$$
and $\mathcal A^2_n=\mathcal S_{r_n}\setminus\mathcal A^1_n$.
For each $x\in I_{m_0}$ and $s=2,p$, we split $\Lambda_s(x,n)$ as
$$\Lambda_s(x,n)=H_s(x,\mathcal A^1_n)+
H_s(x, \mathcal A^2_n),$$
where  $H_s(x,\mathcal A^i_n)$, $i=1,2$ is defined by
$$H_s(x,\mathcal A^i_n)=
\sum_{j\in \mathcal A^i_n}
|\Psi_{j,n}(x)|^s\Lambda_s(\upsilon_{j-1},Q_{m_j})$$

Notice that
 $H_p(x, \mathcal A^2_n)$ 
is the sum of at most $M_\epsilon$ bounded terms. Therefore, for some
constant $C_\epsilon>0$
\beqnl{firstpiececirc}
H_p(x,\mathcal A^2_n)\leq C_\epsilon
\eeqnl

For the term $H_p(x, \mathcal A^1_n)$ we have
\beqnl{secondpiececirc}
H_p(x, \mathcal A^1_n)\leq \epsilon (H_2(x,\mathcal A^1_n))^{p/2}
\leq \epsilon (\Lambda_2(x,n))^{p/2}
\eeqnl
Therefore, combining  \eqref{firstpiececirc}, \eqref{secondpiececirc}
 we get
\beqn
\frac{\Lambda_p(x,n)}{(\Lambda_2(x,n))^{p/2}}
\leq \frac{C_\epsilon}{(\Lambda_2(x,n))^{p/2}} +\epsilon
\eeqn

Limit \eqref{circcurtoratiolac} follows from \eqref{biglambdacirclac}.

\bigskip
\noindent
(3) If condition \eqref{hypothesiscirc} holds, then
limit \eqref{circcurtoratio} follows from \eqref{biglambdacirc}.
\endproof
\begin{remk}
One natural hypothesis that implies 
condition \eqref{hypothesiscirc} holds is that 
 $s_k=K_k(\lambda_k^2)$. In this case,
$H(\lambda_k^2)\lambda_k^{6k}=1$. Then, \eqref{hypothesiscirc}
follows from  Proposition \ref{comparison2}.
\end{remk}
Since
$$Q_m=\frac{\beta^{-m}-(-1)^m\beta^m}{\sqrt{5}}$$
for all $m$, it follows that $m_0(n)\sim [\log_{\beta^{-1}}(n)]+1$ 
as  $n\rightarrow\infty$.

In the following result, we show that the Lyapunov condition
\ref{lyapcon} holds for all orbits starting at points of the form
$x=f^l(0)$, $l\in\mathbb N$.

\begin{thm}\label{lyapconorbit0circ}
\begin{itemize}
\item[1.] If $\{n_i\}$ is an increasing sequence of integers
that satisfy \eqref{laccirc}, then
\beqnl{ratiocirclac}
\lim\limits_{i\rightarrow\infty}\sup_{x\in I_{m_0(n_i+l)}}
\frac{\Lambda_p(f^l(x),n_i)}{\left(\Lambda_2(f^l(x),n_i)\right)^{p/2}}=0
\eeqnl
for all $p>2$.
\item[2.] Under condition \eqref{hypothesiscirc}
\beqnl{ratiocirc}
\lim\limits_{i\rightarrow\infty}\sup_{x\in I_{m_0(n_i+l)}}
\frac{\Lambda_p(f^l(x),n_i)}{\left(\Lambda_2(f^l(x),n_i)\right)^{p/2}}=0
\eeqnl
for all $p>2$.
\end{itemize}
\end{thm}
\proof
We will use a method  similar to the one used in  Theorem \ref{kurtnx}.
since the proof of (1) and (2) are very similar, 
we will only prove (2).

 Let $l\in \mathbb N$
be fixed. For any  $n\in\mathbb N$, consider 
$$n+l = Q_{m_0}+\cdots Q_{m_r}$$
where $\{m_j\}$ is as in \eqref{fibodecomp}. For
$x\in I_{m_0(n+l)}$, define the sequence
of returns $\{\upsilon_j\}_j$ with $\upsilon_{-1}=x$ as
in \eqref{retcirck}. 
Notice that
\begin{multline}
\Lambda_p(f^l(x),n)=|\Psi_{n+l,0}(x)|^p A_p(x,Q_{m_0})+
\label{almostfibo}\\
\sum_{i=1}^{r}|\Psi_{n+l,i}(x)\|^p\Lambda_p(\upsilon_{i-1},Q_{m_i})
\end{multline}
where
\beqnl{differencecirc}
\Lambda_p(x,Q_{m_0})-A_p(x,Q_{m_0})=
\sum_{j=1}^{l}\left|\left(f^{Q_{m_0}-j}\right)\circ f^j(x)\right|^p
\eeqnl
From the identity
\beqn
\left(f^{Q_{m_0}-j}\right)^\prime\circ f^{j}(x)=
\frac{\left(f^{Q_{m_0}}\right)^\prime(x)}{\left(f^j\right)^\prime(x)}
\eeqn
and Lemma \ref{growthcircle}, it follows that 
\beqnl{goodboundcirc}
\frac{c_j}{|\lambda_k|^{2km_0}}
\leq \left|\left(f^{Q_{m_0}-j}\right)^\prime\circ f^j(x)\right|
\leq \frac{d_j}{|\lambda_k|^{2km_0}}
\eeqnl
where  $c_j$ and $d_j$ depend only  on $j$. 
Since\eqref{differencecirc} contains $l$ terms, Theorem
 \eqref{comparison2} and Corollary \ref{kurtoQnresult} imply that
\beqn
\frac{\Lambda_p(x,Q_{m_0})-A_p(x,Q_{m_0})}
{\Lambda_p(x,Q_{m_0})}\leq C_{l,p}\frac{1}{(\lambda_k^{2kp}\rho_p)^{m_0}}
\rightarrow0
\eeqn
Therefore,
\beqnl{blockasympcirc}
\lim_{n\rightarrow}\sup_{x\in I_{m_0(n+l)}}
\frac{\Lambda_p(f^l(x),n)}{\Lambda_p(x,n+l)}=1
\eeqnl
and the limit \eqref{ratiocirc} follows from Proposition 
\eqref{kurto0ncirc}.
\endproof

\subsubsection{ Proof of Theorems \ref{circclt} and \ref{circclt2}}
Recall that for large integers $n$, 
$m_0(n)\sim[\log_{\beta^{-1}}(n)]+1$ and 
$r_n+1$ equals the number of terms in the Fibonacci expansion 
of $n$.

By Theorem \ref{main} and Theorem \ref{lyapconorbit0circ},
 it suffices to choose $\sigma_n$ that satisfies \eqref{weaknoise}.

Let $c$ and $d$ be as in Lemma \ref{growthcircle}. 
From \eqref{blockasympcirc} and Proposition
\ref{kurto0ncirc}, for each
$p>0$ there  there are constants $C_{p,l}$ and $D_{p,l}$
such that for all for all $x\in I_{m_0(n+l)}$,
\beqn
C_{p,l}\left(\frac{(c\, \lambda_k^{6k})^{r_{n+l}}}
{\lambda_k^{(2k-1)m_{r_{n+l}}}}\right)^p 
\leq 
\frac{\Lambda_p(f^l(x),n)}
{(\lambda_k^{2kp}\rho_p)^{m_0(n+l)}}
 \leq  D_{p,l}\left(\frac{d^{r_{n+l}}}
{\lambda_k^{(2k-1)m_{r_{n+l}}}}\right)^p
\eeqn

 Therefore, we have that
\beqnl{sigoptcirc}
1\leq \frac{(\widehat{\Lambda}(f^l(x),n))^3}{\sqrt{\Lambda_2(f^l(x),n)}}
\leq C\,\left(\frac{d^3}{c\, \lambda_k^{2k}}\right)^{r_{n+l}}
\left(\frac{\rho_1^3}{\sqrt{\rho_2}}\right)^{m_0(n+l)}
\eeqnl
for some constant $C>0$.

If we consider an increasing sequence of iterations $\{n_k\}_k$ which
have  lacunar Fibonacci expansions, see \eqref{laccirc}.
Then,  from \eqref{sigoptcirc} it follows that  it suffices to consider
$\sigma_{n_k}=n_k^{-(\gamma+1)}$
to have a  central limit holds along $\{n_k\}$.

Otherwise, under condition \eqref{hypothesiscirc} we have the following
cases:

Case (a) If  
\beqnl{caseobscirc}
\left|\frac{d^3}{c\, \lambda_k^{2k}}\right|\leq 1
\eeqnl
then,  it suffices to consider 
$$\sigma_m=\frac{1}{n^{\delta+1}}$$
with $\delta> \log_{\beta^{-1}}(\rho_1^3/\sqrt{\rho_2})$.

Case (b) If
\beqnl{casevaccirc}
\left|\frac{d^3}{c\, \lambda_k^{2k}}\right|> 1,
\eeqnl
it is enough to consider 
$$\sigma_n=\frac{1}{n^{\delta^*+1}}$$
with 
$\delta^*> \log_{\beta^{-1}}(d^3\rho_1^3)+
\log_{\beta^{-1}}(c\,\lambda_k^{2k}/\sqrt{\rho_2})$.

\bigskip

The Berry--Esseen estimates follow by choosing 
$\sigma_n$ that satisfies  \eqref{weakernoise}
in  Theorem \ref{berry-esseen}.

If \eqref{caseobscirc} holds, then it 
suffices to consider
$$\sigma_n=\frac{1}{n^\tau}$$
with
$\tau> \log_{\beta^{-1}}(\rho_1^3/\sqrt{\rho_2})+
\log_{\beta^{-1}}(\rho_2^3/\rho_3)+1$.

In the case \eqref{casevaccirc}, it suffices to consider
$$\sigma_n=\frac{1}{n^{\tau^*}}$$
with 
$\tau^*> \log_{\beta^{-1}}(d^3 \rho_1^3)+
\log_2(c\,\lambda_k^{2k}/\sqrt{\rho_2})
+ \log_2(\rho_2^3/\rho_3)$ +1
$\hfill{\Box}$

\begin{remk}
We know from the theory of
critical circle maps that $\{f^l(0)\, \mod 1: l\in\mathbb N\}$ 
is dense in $\mathbb T^1$.
Numerics in \cite{DiazL} suggests that  the Lyapunov
condition holds for all points in the circle. We conjecture that
this is indeed the case and that 
that the effective noise
of orbits starting
 in the any arbitrary point of $\mathbb T^1$ approaches
a Gaussian. We expect the   the speed of convergence 
in the central limit theorem to be not uniform.
\end{remk}
\section{Possible extensions of the results}
In this section, we suggest 
several extensions of the results in this paper that are
presumably be accessible. 
\begin{itemize}
\item[1.] Assume that $(\xi_n)$ is a sequence of independent 
random variables with mean zero $p$ finite 
moments  such that
$$A_-n^{\alpha_-}\leq \|\xi_n\|_2\leq  \|\xi_n\|_p\leq A_+ n^{\alpha_+},$$
with some $\alpha_{\pm}$ in a small range. 
\item[2.] Assume that the  random variables $(\xi_n)$
are weakly correlated (e.g. Martingale approximations).
\end{itemize}
 These assumptions are  natural in dynamical
systems applications when the noise is generated by a discrete
process. That is
\begin{eqnarray*}
x_{n+1}=f(x_n)&+& \sigma \psi(y_n)\\
y_{n+1}=h(y_n)
\end{eqnarray*}
and $h$ is an expanding map or an Anosov system.

\begin{itemize}
\item[3.] 
Related to the central limit theorem (even in the
case independent random variables $(\xi_n)$
of comparable sizes), it also would be desirable to
obtain higher order asymptotic expansions in the 
convergence to Gaussian, namely Edgeworth expansions.
\item[4.] We note that the estimates for the asymptotic growth of
 the variance of the effective
noise (\eqref{lambdaasympat0} with $p=2$) for
systems at the accumulation of period doubling are obtained
in \cite{vul:feiguniv} using the Thermodynamic formalism.
We think that it would be very interesting to develop
analogues to the log--convexity properties of the Lindeberg--Lyapunov
operators  or the Edgeworth 
expansions with the  thermodynamic formalism.
\end{itemize}

\section*{Acknowledgments}
We are 
  grateful to Prof. \`Alex Haro his support and 
encouragement in the realization of this paper.
We  wish to acknowledge 
  Professors Hans Koch and Klaus Bichteler for the
invaluable discussions
and comments of the results of this paper. The authors 
enjoyed very interesting discussions with D. Khmelev.
We miss him as a mathematician and as a generous colleague.

We thank the Fields Institute for their invitation 
and support to present part of this work
during the Renormalization in dynamical systems workshop in
December 2005. Both authors were supported in part by NSF Grants.


\bibliographystyle{alpha}

\def\polhk#1{\setbox0=\hbox{#1}{\ooalign{\hidewidth
  \lower1.5ex\hbox{`}\hidewidth\crcr\unhbox0}}}
  \def\polhk#1{\setbox0=\hbox{#1}{\ooalign{\hidewidth
  \lower1.5ex\hbox{`}\hidewidth\crcr\unhbox0}}}
  \def\polhk#1{\setbox0=\hbox{#1}{\ooalign{\hidewidth
  \lower1.5ex\hbox{`}\hidewidth\crcr\unhbox0}}}
  \def\polhk#1{\setbox0=\hbox{#1}{\ooalign{\hidewidth
  \lower1.5ex\hbox{`}\hidewidth\crcr\unhbox0}}}
  \def\polhk#1{\setbox0=\hbox{#1}{\ooalign{\hidewidth
  \lower1.5ex\hbox{`}\hidewidth\crcr\unhbox0}}} \def\cprime{$'$}

\end{document}